%% file: main.tex
\documentclass[preprint,12pt]{elsarticle}

\input{Packages_and_Command}

\begin{document}

\begin{frontmatter}

\title{Minimum-norm interpolation for unknown surface reconstruction}

\author{Alex Shiu Lun Chu\corref{cor1}\fnref{label1}} 
\ead{222482318@life.hkbu.edu.hk}
\cortext[cor1]{Corresponding author}

\author{Leevan Ling\fnref{label1}}
\ead{lling@hkbu.edu.hk}

\author{Ka Chun Cheung\fnref{label2}}
\ead{chcheung@nvidia.com}

\affiliation[label1]{organization={Department of Mathematics, Hong Kong Baptist University},city={Kowloon Tong, Kowloon},country={Hong Kong}}

\affiliation[label2]{organization={NVIDIA AI Technology Center NVAITC, NVIDIA}, country={USA}}

\begin{abstract}
We study algorithms to estimate geometric properties of raw point cloud data through implicit surface representations. Given that any level-set function with a constant level set corresponding to the surface can be used for such estimations, numerical methods need not specify a unique target function for these domain-type interpolation problems.
In this paper, we focus on kernel-based interpolation by radial basis functions (RBF) and reformulate the uniquely solvable interpolation problem into a constrained optimization model. This model minimizes some user-defined norm while enforcing all interpolation conditions.
To enable nontrivial feasible solutions, we propose to enhance the trial space with 1D kernel basis functions inspired by Kolmogorov--Arnold Networks (KANs).
Numerical experiments demonstrate that our proposed mixed-dimensional trial space significantly improves surface reconstruction from raw point clouds. This is particularly evident in the precise estimation of surface normals, outperforming traditional RBF trial spaces including the one for Hermite interpolation.
This framework not only enhances the processing of raw point cloud data but also shows potential for further contributions to computational geometry. We demonstrate this with a point cloud processing example.
\end{abstract}

\end{frontmatter}

\section{Introduction}\label{sec:intro}
Surface reconstruction from raw point cloud data is a fundamental task in computer vision and image processing, with applications spanning rendering \cite{Shape-implicitSurf}, shading \cite{surface-shading}, filtering \cite{ptcld-filtering}, and point cloud registration \cite{ptcld-registration}. Over the years, various reconstruction approaches have emerged that differ primarily in their reliance on surface normals. Some methods, such as Poisson surface reconstruction \cite{Poisson-Surfrecon} and HRBF-based techniques \cite{HRBF-surfrecon}, require accurate normal vectors. Other methods, such as tangent plane estimation \cite{surface-reconstruction-PCA}, Moving Least Squares (MLS) \cite{movinglsq}, Modified Gauss Formula \cite{gauss-formula} and radial basis function (RBF) approaches \cite{meshless-MATLAB, reconstruct-3d-rbf, SurfRecon-RBF-GPU, Implicit-CurlFreeRBF, Implicit-GRBF} fit surfaces directly from point locations, thereby reducing but not entirely eliminating the dependency on precise normal data. In these approaches, normal vectors are often approximated using methods such as Principal Component Analysis (PCA) which is efficient but typically yields lower accuracy. Lastly, there are methods that do not require any normal vectors such as the level set method \cite{FastSurfRecon-levelset} or Delaunay triangulation \cite{surfrecon-delaunay}. Although these traditional approaches yield acceptable results on smooth surfaces, they often struggle with complex geometries.

In contrast, recent advances in GPU technologies have propelled the development of deep learning methods. Many deep learning methods have been introduced in recent years, and numerous studies have compared these methods based on various criteria. According to one of the criteria proposed in \cite{surfrecon-survey}, deep learning approaches can be divided into two groups. One group comprises methods that require surface normals during training such as ParSeNet \cite{ParSeNet}, Occupancy Network \cite{OccNet}, DeepSDF \cite{DeepSDF} and LIG \cite{LocalImpliciGrid}. The other group includes methods that do not require surface normals during training such as PCPNet \cite{PCPNet}, DeepFit \cite{Deepfit}, Nesti-Net \cite{Nesti-Net}, NeAF \cite{NeAF} and HSurf \cite{HSurf}. These methods rely on extensive training datasets, feature extraction and complex architecture optimizations to capture fine geometric details which significantly enhance accuracy (especially for surfaces with sharp edges and corners). However, they tend to be computationally demanding and often lack interpretability. In this work, we introduce a novel approach that enhances the precision of surface normal estimation while integrating effectively with these advanced deep learning frameworks to support the ongoing evolution in point-to-surface analysis rather than competing with existing methods.

In this paper, we explore enhanced algorithms for estimating geometric properties by constructing local implicit surface representations directly from raw point cloud data without requiring pre-computed normal vectors. Our primary focus is on refining the kernel-based interpolation framework using RBF methods, which have traditionally been employed for this purpose. Initially, we reformulate the interpolation problem as a constrained optimization model that minimizes the native space norm subject to appropriate interpolation conditions. This reformulation guarantees that the problem remains uniquely solvable, much like the original interpolation formulation but yields markedly more accurate normal estimates. Consequently, our method can serve as an effective pre-processing step for any surface reconstruction technique that requires pre-computed normal vectors, providing higher quality normals that improve downstream reconstruction fidelity.

Moreover, our framework can also serve as an analysis tool for studying deep learning-based normal estimation. To expand the set of feasible solutions within the optimization search space, we propose two expansions of the trial space where one incorporates the traditional Hermite radial basis function (HRBF) and another one employs 1D kernel basis functions, which can be viewed as forming a sub-optimal variant of Kolmogorov--Arnold Networks (KANs) \cite{KANs}. These sub‐optimal KANs bridge classical kernel methods and deep learning approaches, offering insight into the capacities and limitations of data‐driven approaches. This strategic enhancement results in a mixed-dimensional trial space that substantially improves the accuracy of surface reconstruction from raw point clouds, particularly in terms of surface normal estimation.

The remainder of this paper is organized as follows: \autoref{sec2} provides an overview of implicit surface reconstruction and describes the traditional RBF approach used to address these problems. In \autoref{sec3}, we introduce our optimization method based on the RBF technique which is designed to derive a well-defined implicit function for the interpolation problem in the embedding space. \autoref{sec4} discusses strategies for expanding the trial space.
Numerical results are presented in \autoref{sec5}, including comparison of different configurations of the proposed approach, convergence study on torus and sphube, surface normal estimation on a flattening sphube and a pillar, as well as a denoising example of a cube. Finally, \autoref{sec6} concludes the paper with remarks and potential directions for future research.

\section{Notations and Preliminaries} \label{sec2}
Given a smooth codimension-1 manifold \( \M \) embedded in the ambient space \( \R^3 \), an implicit function representation of \( \M \) can be provided by any smooth scalar function \( F_\M: \Omega_\M \to \R \) for some embedded domain \( \M \subset \Omega_\M \subseteq \R^3 \) such that
\begin{equation}
    \M = \{ \x \mid F_\M(\x) = C, \text{ for some }C \in \R \},
    \label{eq:implicit_surface}
\end{equation}
and \( \nabla F_\M(\x) \neq 0 \) for all \( \x \in \M \). We explicitly keep a constant \( C \) here instead of the conventional choice of \( C = 0 \) for a numerical reason later.
By the fact that  $\M$ is the $C$-level set of $F_\M$, we can then compute the surface normal vector at any $\x=[x,y,z]^T\in\M$:
\begin{equation}\label{eq:surface_normal}
    \n(\x)=\frac{\nabla F_\M(\x)}{||\nabla F_\M(\x) ||}.
\end{equation}
The principal curvatures and directions are intrinsic geometric properties derived from the shape operator (or Weingarten map), which itself is computed from the Hessian matrix of $F_\M$:
\begin{equation}
    \nabla^2F_\M = \begin{bmatrix}
        \frac{\partial^2 F_\M}{\partial x^2} & \frac{\partial^2 F_\M}{\partial x \partial y} & \frac{\partial^2 F_\M}{\partial x \partial z} \\
        \frac{\partial^2 F_\M}{\partial y \partial x} & \frac{\partial^2 F_\M}{\partial y^2} & \frac{\partial^2 F_\M}{\partial y \partial z} \\
        \frac{\partial^2 F_\M}{\partial z \partial x} & \frac{\partial^2 F_\M}{\partial z \partial y} & \frac{\partial^2 F_\M}{\partial z^2}
    \end{bmatrix}.
\end{equation}
We also define the projection operator  $P_\M = I - \n \otimes \n$ for all $\x \in \M$, where $I$ is the $3\times3$ identity matrix. The shape operator $S: T_{\mathbf{x}}\M \to T_{\mathbf{x}}\M$ for any vector $\mathbf{v} \in T_{\mathbf{x}}\M$ is defined by $S\mathbf{v} = -D_{\mathbf{v}} \mathbf{n}$ that can be expressed as:
\begin{equation}\label{S=PHP}
  S = -\frac{1}{||\nabla F_\M  ||}P_\M\,\nabla^2F_\M\, P_\M.
\end{equation}
Our definition results in a $3\times3$ shape operator when using the standard basis. This is consistent with the approach outlined in \cite{multi-scale-ptcld, stable-efficient-ptcld}, which derives a  $2\times2$ shape operator. In that work, the operator is constructed using a transfer matrix to express tangent vectors in terms of an orthonormal basis of the tangent plane.
The nonzero eigenvalues of the shape operator $S$ in \eqref{S=PHP} represent the principal curvatures $\kappa_1$ and $\kappa_2$, and their corresponding eigenvectors are the principal directions, which are tangent to the surface at $\x$.
The Gaussian curvature is the product of the principal curvatures \( \kappa_1 \kappa_2 \) providing insight into the intrinsic curvature of a surface at a point, whereas the mean curvature  is the average of these principal curvatures \(  \frac1{2}(\kappa_1 + \kappa_2) \)  measuring how a surface bends in the embedding space. Closed formulas for these curvatures can be found in \cite{curvature-formula}.

\subsection{Implicit surface reconstruction from point cloud}
\begin{figure}
    \centering
    \subfloat[Point cloud and a local stencil
    ]{\includegraphics[width=.5\textwidth]{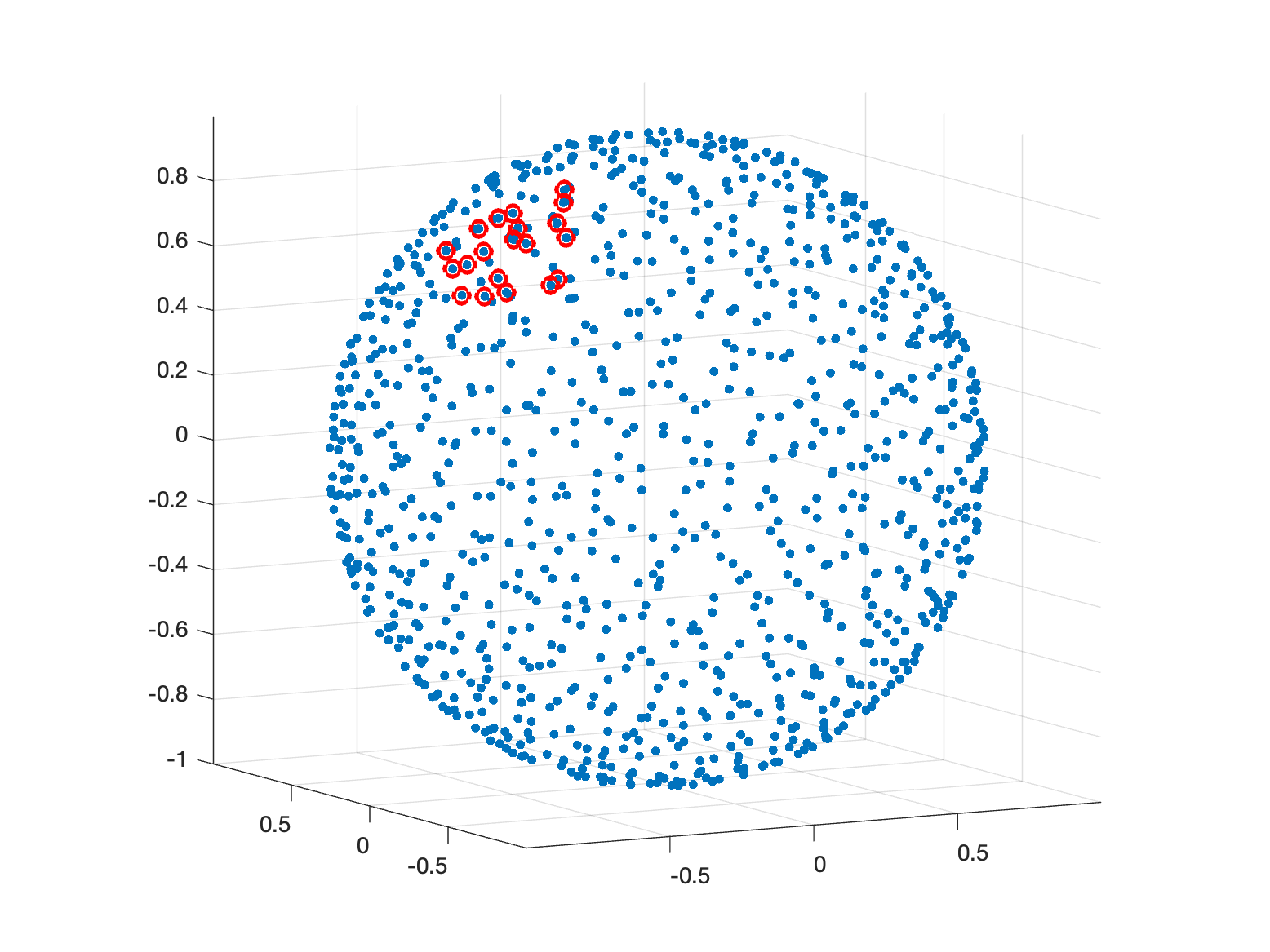}}
    \subfloat[Local stencil and ghost points
    ]{\includegraphics[width=.5\textwidth]{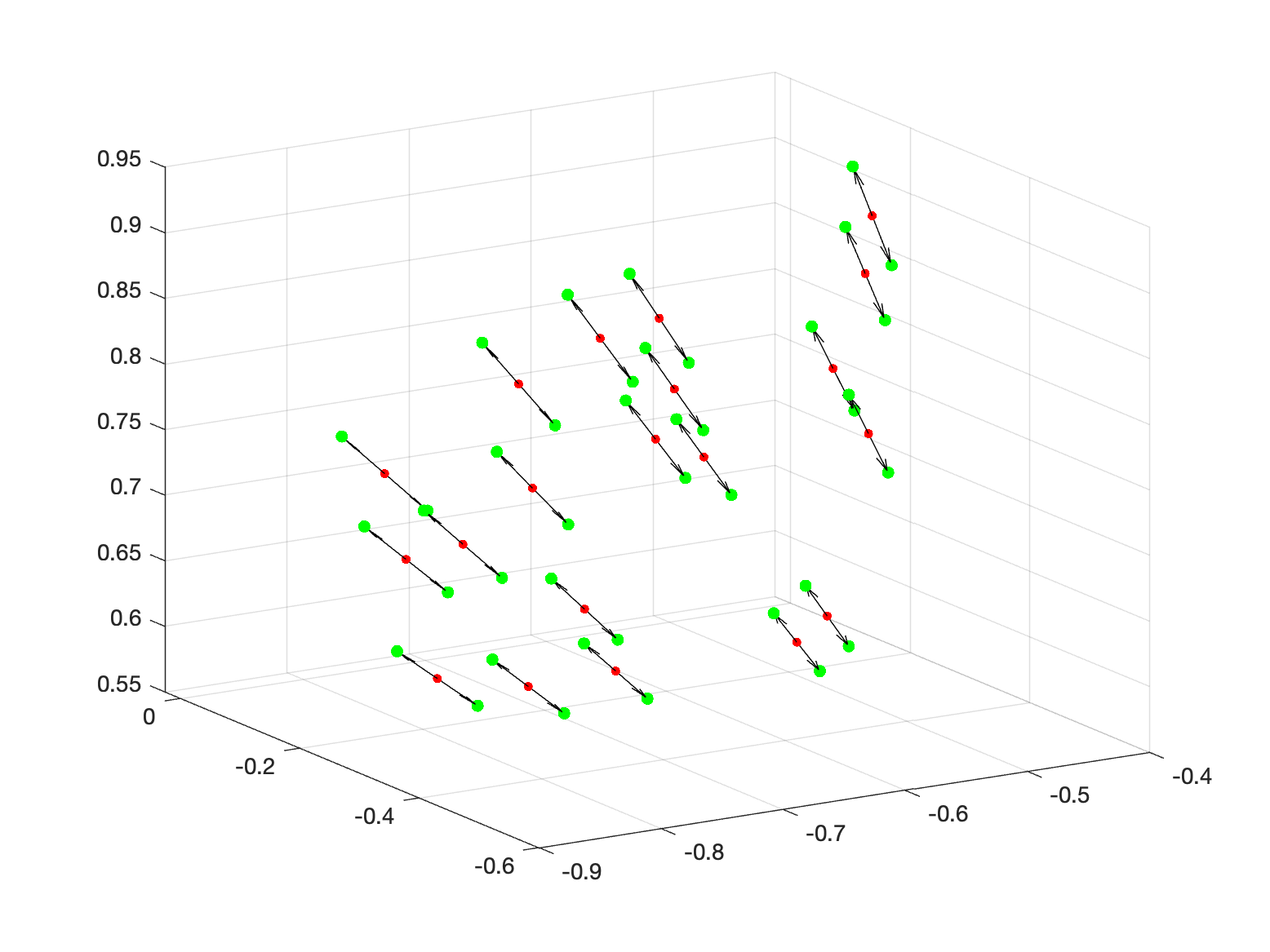}}
\caption{Local Stencils and Ghost Points in Point Cloud Data. (a) Illustrates the point cloud data $P$ shown in blue, with one local stencil highlighted in red. (b) Displays the same local stencil in red, augmented with off-surface (ghost) points $P_\pm$, depicted in green, which are generated using the offset parameter $\new{\delta}$ and estimated normals $\n$.}
\label{fig:Pointcloud_stencil}
\end{figure}

The goal in reconstructing a surface from a raw point cloud is to approximate an unknown manifold $\mathcal{M} \subset \mathbb{R}^3$.
Let  $P = \{\p_1, \ldots, \p_n\} \subset \M$ denote the given  point cloud dataset containing $n=|P|$ distinct surface points of $\M$.
We seek an implicit function $F_P$, which must be sufficiently smooth to approximate a function $F_{\mathcal{M}}$ that describes the surface $\M$. The function $F_P$ is defined such that it satisfies the following interpolation condition at each point in $P$:
\begin{equation}\label{eq:implicit_surface_ptcloud}
F_P(\p_i) = C \quad \text{for each } \p_i \in P,
\end{equation}
where $C$ is a constant, assumed nonzero unless additional conditions specify otherwise. This condition ensures that $F_P$ closely models the underlying surface $\M$.
Using a local stencil reduces the computational complexity and thus makes it easier to handle large datasets and the local geometric properties of the point cloud can also be approximated.
Note that the numerical approaches and procedures described below can be effectively applied both globally to the entire point cloud and locally to specific segments of $\M$. In cases where the focus is local, $P=P_s$ and $n=n_s$ should be understood as referring to local stencils and their respective sizes, see \autoref{fig:Pointcloud_stencil}(a). For simplicity, we use the same notation for both global and local contexts; the subscript-$s$ is omitted when it does not lead to confusion.

It is common practice in surface reconstruction methods \cite{meshless-MATLAB, reconstruct-3d-rbf, SurfRecon-RBF-GPU, Implicit-CurlFreeRBF} to include off-surface points in addition to interpolation conditions in \eqref{eq:implicit_surface_ptcloud}. This approach helps define the interior and exterior of the surface, avoids trivial solutions, and enhances stability. Notably, constructing these off-surface points requires normals; however, instead of needing highly accurate normal vectors, it is sufficient to use approximated normals provided they are correctly oriented. In cases where the surface normals are unavailable, estimates obtained via principal component analysis (PCA) \cite{surface-reconstruction-PCA,PCA-normal} or least squares plane fitting \cite{Implicit-GRBF} are generally adequate. By relying on properly oriented approximations, the reconstruction function becomes well-defined without the strict dependence on high-fidelity normals that other surface reconstruction techniques often require. Nevertheless, obtaining a complete set of off-surface points with consistently correctly oriented normals remains challenging, and algorithms such as SNO \cite{SNO} and ODP \cite{ODP} may still be needed to assign the proper direction to the normals.

To enhance the interpolation defined in \eqref{eq:implicit_surface_ptcloud}, one can integrate $2n$ off-surface (ghost) points \cite{levelset-shape,Shape-levelset} by expanding the point cloud data $P$ with an estimated normal $\n$ at each point to an $N=3n$ point set:
\begin{equation}\label{Ppm}
 P \cup P_\pm := P \cup \{\p \pm \new{\delta}\n(\p) \mid \p \in P\},
\end{equation}
where $\delta$ is an offset parameter that must be carefully selected to avoid crossing of the level set \cite{RBF-mesh-repair}, see \autoref{fig:Pointcloud_stencil}(b). Assigning level set values $(C, C\pm \new{\delta})$ to $P$ and $P_\pm$ ensures that the target function in \eqref{eq:implicit_surface} is the (uniquely defined) \textit{signed distance function}.

In algorithms that use a single stencil (local method) to approximate geometric quantities at a specific surface point, employing the entire set $P_\pm$ may be excessive. Alternatively, one could define just a pair of ghost points at the point of interest $\p_0 \in P$ within a local stencil, thereby reducing the size of the interpolation problem to $N=n+2$ points:
\begin{equation}\label{P+2}
\tilde{P}_0 := P \cup \{\p_0 \pm \new{\delta}\n(\p_0)\}.
\end{equation}
This approach results in a manageable $N\times N$ interpolation problem \reviewB{that reduces the per-stencil cost compared with augmenting all $2n$ points as in \eqref{Ppm},  while still avoiding the trivial constant-level solution and enforcing a signed-distance-like behavior near $p_0$.} However, the target unknown function of the interpolation remains undefined and is dependent on how the chosen numerical method extrapolates the approximate function away from $\p_0$ and $P$.
\reviewB{
\begin{remark}
Since off-surface information is imposed only at $p_0$, the interpolant can be less constrained away from $p_0$ when sampling is highly nonuniform, noise is present, or the stencil crosses sharp features. This may reduce robustness and degrade the accuracy of the reconstructed implicit function and its derivatives. In such cases, one may increase the local stencil density and/or add ghost points for a small subset of stencil points to provide additional off-surface constraints.
\end{remark}
}

\subsection{Kernel-based RBF Interpolation Method}

To resolve the interpolation problems described previously (with or without ghost points), we select a set of basis functions to form a trial space and seek an interpolant in this space. A standard choice is RBF interpolation method. In this work, we consider translation-invariant symmetric kernels $\Phi:\mathbb{R}^3\times\mathbb{R}^3\to\mathbb{R}$ that are strictly positive definite (SPD). For any set of $N:=|\Xi|$ distinct interpolation points $\Xi\subset\Omega_\M$, either $\Xi=P$ or its extensions as defined in \eqref{Ppm} or \eqref{P+2}, we define the trial space with centers in $\Xi$ \cite{Discrete-lsq-rbf}:
\begin{equation}
\mathcal{U}_{\Xi,\Phi} := \text{span} \{\Phi(\|\cdot - \p_j\|) : \p_j \in \Xi\}.
\label{eq:RBF_trial_space}
\end{equation}
We then seek an implicit function in \eqref{eq:implicit_surface_ptcloud} of the form
\begin{equation}\label{eq:RBF_ansatz}
F_{\Xi,\mathcal{U}_{\Xi,\Phi}}
= \sum_{\p_j \in \Xi}\lambda_j\Phi(\|\cdot - \p_j\|) \in \mathcal{U}_{\Xi,\Phi},
\end{equation}
where the coefficients $\lambda_j$ are found by enforcing the interpolation conditions. 
\reviewA{For SPD kernels and distinct points, the associated kernel matrix $\Phi(\Xi,\Xi)$ is symmetric positive definite and hence invertible, so the interpolation problem has a unique solution for any prescribed data values \cite{meshless-MATLAB,scatter-data,sobolev-error-estimate}.}
Using the estimated signed distances to the $C$-level set of $\M$ as function values in \eqref{eq:implicit_surface_ptcloud}, the interpolation conditions give the linear system
\begin{equation}\label{RBF Ax=b}
  F_{\Xi,\mathcal{U}_{\Xi,\Phi}}(\Xi) 
  = \left[\sum_{\p_j \in \Xi}\lambda_j\Phi(\|\p_i - \p_j\|)\right]_{i=1}^N
  = C+\text{dist}_\M(\Xi):=
  \left\{
  \begin{array}{ll}
    C   & \text{if } \p_i \in P, \\
    C+\new{\delta} & \text{if } \p_i \in P_+, \\
    C-\new{\delta} & \text{if } \p_i \in P_-. \\
  \end{array}
  \right.
\end{equation}
Although we use the standard interpolation procedure in \eqref{RBF Ax=b} for $F_{\Xi,\Phi} = I_{\Xi,\Phi}F_\M$, except in the case of \(N=3n\) data points in \eqref{Ppm}, we do not know the form of the target \(F_\M\). To address this issue, we propose a constrained optimization framework that selects an interpolant by minimizing a user-defined norm while enforcing the interpolation constraints. Details are given in \autoref{sec3}.

\section{Minimum-norm Interpolation} \label{sec3}
Given the nonuniqueness nature of the implicit function $F_\M$, it is logical to seek a function with the smallest Sobolev norm. Theoretically, this approach is supported by the scattered zeros lemma in $\Omega_\M$ \cite{sobolev-bound} and on $\M$ \cite{sobolev-error-estimate}. This lemma asserts that if a smooth function vanishes at a set of scattered points, the Sobolev norms of the function are controlled by the distribution of these zeros.

Consider the detailed analysis of \cite[Theorem 1.1]{sobolev-bound} and apply it to the difference $I_{\Xi,\Phi}F_\M-F_\M$, where the zeros are scattered at $\Xi \subset \Omega_\M$. Define the fill distance as
    \begin{equation}
        h = h_{\Xi,\Omega_\M} := \sup_{\x \in \Omega_\M} \min_{\p_j \in \Xi} \|\x - \p_j\|_2.
    \end{equation}
This leads to a simplified (with $p=2=q$ in the original lemma there) fractional Sobolev spaces estimate:
    \begin{equation} \label{eq:errbnd_sobolev}
        \big|I_{\Xi,\Phi}F_\M - F_\M\big|_{W_2^{|\alpha|}(\Omega_\M)} \leq c h^{\tau - |\alpha|} \big|I_{\Xi,\Phi}F_\M - F_\M\big|_{W_2^{\tau}(\Omega_\M)},
    \end{equation}
for some sufficiently large $\tau \in \R^+$ with respect to the multi-index $\alpha$, and some constant $c$ independent of $F_\M$ and $h$.
This gives an interpolation error bound to any interpolation process with property
$ |I_{\Xi,\Phi}F_\M  |_{W_2^{\tau}(\Omega_\M)}\leq  | F_\M |_{W_2^{\tau}(\Omega_\M)}$.

The native space $ \mathcal{N}_\Phi(\Omega_\M) \supset \mathcal{U}_{\Xi,\Phi}$ for all sets $\Xi\subset\Omega_\M$ of trial centers naturally arises in RBF interpolation methods. For any trial function $u\in\mathcal{U}_{\Xi,\Phi}$ in the form in \eqref{eq:RBF_ansatz}, its native space norm is defined by
\begin{equation}
    \|u\|_{\mathcal{N}_\Phi(\Omega_\M)} := \sum_{\p_i,\p_j\in \Xi}\lambda_i\lambda_j\Phi(\p_i,\p_j) = \boldsymbol\lambda^T\Phi(\Xi,\Xi)\boldsymbol\lambda.
    \label{eq:native_space}
\end{equation}
For any $f \in \mathcal{N}_\Phi$, the interpolant to $f$ on $\Xi$ from $\mathcal{U}_{\Xi,\Phi}$ satisfies $\|f - I_{\Xi,\Phi}f \|_{\mathcal{N}_\Phi} =\min_{u\in\mathcal{U}_{\Xi,\Phi}} \|f - u \|_{\mathcal{N}_\Phi}  \leq \| f \|_{\mathcal{N}_\Phi}$, similar to the upper bound in \eqref{eq:errbnd_sobolev}.
The aim of having $\mathcal{N}_\Phi$ be norm-equivalent \cite{nativespace1, nativespace2, ReproducingKernel} to $W_2^{\tau}(\Omega_\M)$, i.e.,
there exist constants $c_1$ and $c_2$ such that:
        \begin{equation}
            c_1 \|\cdot\|_{W_2^\tau(\Omega_\M)} \leq \|\cdot\|_{\mathcal{N}_\Phi(\Omega_\M)} \leq c_2 \|\cdot\|_{W_2^\tau(\Omega_\M)},
        \end{equation}
motivates the use of Sobolev space reproducing kernels, whose Fourier transform decays as:
\begin{equation}
    \hat{\Phi}(\omega) = (1 + \|\omega\|^2)^{-\tau}.
    \label{eq:ms_FT}
\end{equation}
The native space of such kernels is defined as:
        \begin{equation}\label{eq:sobolev_space}
            \mathcal{N}_\Phi(\mathbb{R}^d) = \{f \in L_2(\mathbb{R}^d): \widehat{f}(\cdot)(1 + \|\cdot\|_2^2)^{\tau/2} \in L_2(\mathbb{R}^d) \text{ with } \tau > d/2
            \}.
        \end{equation}
This coincides with the Sobolev space $W_2^\tau(\mathbb{R}^d)$, see \cite{sobolev-space}.
Assuming further that $\Omega_\M$ is a bounded Lipschitz domain, then $\mathcal{N}_\Phi(\Omega_\M) \approx W_2^\tau(\Omega_\M)$, also showing norm equivalence.
A commonly used example is the Mat\'ern Sobolev kernel:
\begin{equation}\label{eq:ms_kernel}
    \Phi(r) = \Phi_{\tau,d}(r) = \mathcal{K}_{\tau-d/2}(r)r^{\tau-d/2} : \mathbb{R}^d \times \mathbb{R}^d \to \mathbb{R}, \text{ for } \tau > d/2,
\end{equation}
where $\mathcal{K}_{\tau-d/2}$ denotes the modified Bessel function of the second kind.

\subsection{Conversion to Constrained Optimization}

If we can choose the target unknown function $F_\M$, it is theoretically beneficial to select one with a minimal Sobolev norm in \eqref{eq:errbnd_sobolev}. This strategy  is logical because smoother target functions generally allow for better approximation accuracy. Although directly computing the Sobolev norm could be computationally costly, employing the native space norm offers a more cost-effective method. An even faster option is the $\ell^2$-norm on the RBF coefficients.
Let $\|\cdot\|_U$ denote a user-selected norm to be minimized for the selection of trial functions from some trial space $\mathcal{U}$.
We now proceed to define the optimization problem by incorporating the interpolation conditions in \eqref{RBF Ax=b} as constraints.
For a sufficiently large $\tau$ and a given set of interpolation data points $\Xi \in \Omega_\M$, we define the approximate implicit function as follows:
\begin{equation}\label{eq:opt_problem}
F_{\Xi, \mathcal{U}} = \underset{f \in \mathcal{U}}{\text{minimize}} \quad \|f\|_U
\quad \text{subject to} \quad f(\Xi) = C + \text{dist}_{\M}(\Xi).
\end{equation}
The associated implicit function $F_\M$ is then defined by
\begin{equation}\label{FM}
F_\M = \underset{f \in W_2^{\tau}(\Omega_\M)}{\text{minimize}} \quad \|f\|_{W_2^{\tau}(\Omega_\M)}
\quad \text{s.t.} \quad f(\M) = C, \ \text{and}\ f(\Xi\setminus\M) = C + \text{dist}_\M(\Xi\setminus\M).
\end{equation}
In this formulation, $F_\M$ captures $\M$ as the $C$-level surface and interpolates the signed distance from all points $\Xi\setminus\M$ outside the surface contained in $\Xi$. This comprehensive approach ensures that the function $F_\M$ meets the specific requirements set forth by the interpolation conditions and norm minimization.
\reviewA{In the discrete setting discussed later in \eqref{eq:KA_RBF_K_norm} and \eqref{eq:KA_RBF_K_matrix}, we write $f$ in the chosen basis so that the objective becomes $\boldsymbol\lambda^T K\boldsymbol\lambda$, where $K$ is the square matrix induced by the selected trial space. For SPD kernels and distinct trial centers, $K$ is symmetric positive definite, so the quadratic objective is strictly convex. With the interpolation conditions giving linear equality constraints, the resulting quadratic program has a unique minimizer whenever the constraints are feasible \cite{Boyd_Vandenberghe_2004}. In our setting with surface points prescribed at $C$ and two ghost points per stencil at $C\pm\delta$ with a consistent normal orientation, the linear constraints are feasible.}

While surface points alone can define the interpolation constraints, our numerical experiments show that adding ghost points improves stability. Ghost points provide off-surface constraints that locally mimic a signed distance function. For the two-ghost point setting in \eqref{P+2}, we take $\Xi=\tilde{P}_0$, so $\Xi\setminus\M$ consists of the two ghost points and the constraint in \eqref{FM} remains unchanged.

While the choice of norm affects smoothness and stability, accuracy is more strongly influenced by the choice of the trial space $\mathcal{U}$.
Expanded trial spaces that are effective in capturing geometric features are discussed in \autoref{sec4}.

\section{Expanded Trial Spaces} \label{sec4}
The selection of a trial space is pivotal as it defines the function space in which the solution is sought, playing a crucial role in the accuracy of the results. Building on the logical progression from the previous section, our objective is to identify a way to expand the trial spaces $\mathcal{U}$ for the optimization problem \eqref{eq:opt_problem} that still effectively approximate $W_2^\tau(\Omega_\M)$.

\subsection{KAN-inspired Trial Space}\label{sec:KA-RBF}
\begin{figure}
    \centering
    \input{sub_opt_KANs_1.tex}
\caption{Schematic representation of the KAN architecture for the anisotropic RBF Neural Network, featuring the newly inspired KAN-RBF basis functions.} 
\label{fig:sub_opt_KANs}
\end{figure}
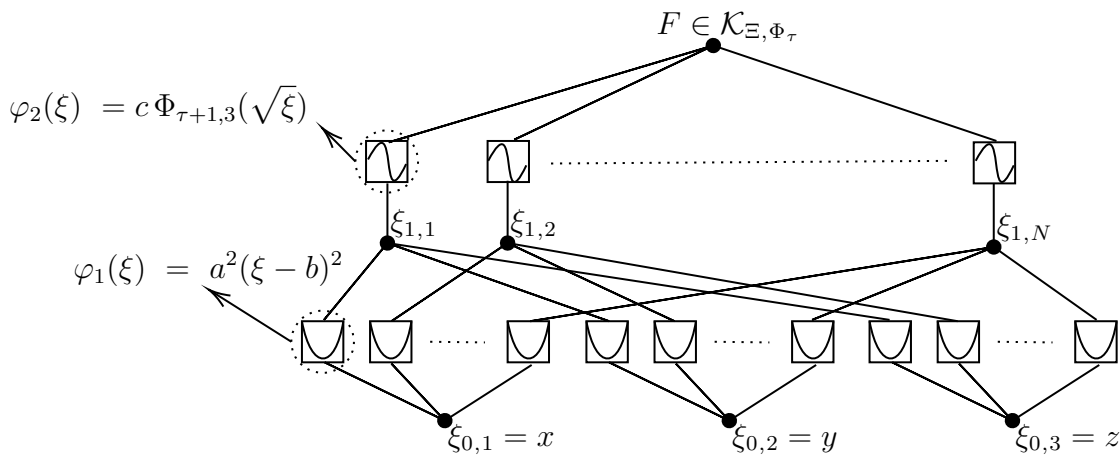

The proposed trial space for minimum-norm  interpolation draws inspiration from the Kolmogorov--Arnold representation theorem, which states that any multivariate continuous function can be expressed as a superposition of continuous single-variable functions \cite{KAtheorem}. This theorem underscores the flexibility and potential of using such structured decompositions in numerical approximations.
Conventionally, RBF interpolation methods employ a trial space comprising exactly $N$ basis functions for $N$ data points. To broaden the feasible set in the optimization problem \eqref{eq:opt_problem}, we propose a trial space by the KAN framework:

\begin{definition}[KAN-inspired trial space]
The KAN-inspired trial space for a distinct point set $\Xi \subset \Omega_\M \subset \R^d$ includes functions generated by the KAN in Figure~\ref{fig:sub_opt_KANs}. This space is formally defined by:
\begin{equation}\label{eq:KA_RBF_trial_space}
    \mathcal{K}_{\Xi,\Phi_{\tau,1}} := 
    \text{\rm span}\Big\{
    \Phi_{\tau,1}\big(|(\cdot-\p_j)^T\mathbf{e}_1 |\big),
    \Phi_{\tau,1}\big(|(\cdot-\p_j)^T\mathbf{e}_2 |\big),
    \Phi_{\tau,1}\big(|(\cdot-\p_j)^T\mathbf{e}_3 |\big): \p_j \in \Xi\Big\}.
\end{equation}
\end{definition}

Figure~\ref{fig:sub_opt_KANs} shows a KAN architecture designed for processing 3D data like a radial basis function Neural Network (RBFNN). The network accepts a 3D input $\p=(x, y, z)$ and begins with a first layer comprising $N$ nodes. Each node employs an activation function $\varphi_1(\xi)=a^2\,(\xi-b)^2$. For the $j$-th node ($j=1,\ldots, N$), the post-activation output is expressed as:
\[
\xi_{1,j} = a_{0,1,j}^2 (x- b_{0,1,j})^2 + a_{0,2,j}^2 (y- b_{0,2,j})^2 + a_{0,3,j}^2 (z- b_{0,3,j})^2 =: \| \p-\mathbf{b}_{0,j}\|_{W_{0,j}}^2,
\]
which is a weighted distance-squared function for the point $\mathbf{b}_{0,j}=[b_{0,1,j},b_{0,2,j},b_{0,3,j}]^T$ with anisotropic weight $ W_{0,j} = \text{diag}[a_{0,1,j}^2,a_{0,2,j}^2,a_{0,3,j}^2]$. When $W_{0,j} = \varepsilon I$ for some $\varepsilon$ and all $j=1,\ldots,N$, this corresponds to traditional RBF with $\varepsilon$ being the shape parameter.
The next activation function  $\varphi_2(\xi) = c\,\Phi_{\tau+1,3}( \sqrt{\xi} )$ completes the configuration of this anisotropic RBFNN. In a fully adaptive KAN, the parameters $(a, b, c)$ would be optimized based on the training data. However, these parameters are fixed for our analysis to yield a sub-optimal network. This fixed configuration allows us to use the trial space as an effective analytical tool, providing a convenient means to compare our approach with conventional methods and to assess the full network's potential.

In our expanded trial space $\mathcal{K}_{\Xi,\Phi_{\tau,1}}$, we set the weights $W_{0,j}=\text{diag}[1,0,0]$ \reviewB{in x-direction} to degenerate the distance function and align the centers with data points in the set $\Xi$ (i.e. $\mathbf{b}_{0,1}=\p_1,\cdots,\mathbf{b}_{0,N}=\p_N$). Thus, we derive a set of 3D RBFs that solely acts on the $x$-coordinates. Note that the Mat\'ern Sobolev kernel in \eqref{eq:ms_kernel} has the property that $\Phi_{\tau+1,3}=\Phi_{\tau,1}$, leading to our KAN-inspired 1D kernel basis $\Phi_{\tau,1}\big(|(\cdot-\p_j)^T\mathbf{e}_1 |\big)$.
By repeating this procedure for the $y$ and $z$ coordinates, we generate a total of $3N$ new basis functions from the span defined in \eqref{eq:KA_RBF_trial_space}.

This trial space, denoted as $\mathcal{K}_{\Xi,\Phi_{\tau,1}}$, may also be interpreted as a degeneration of an anisotropic kernel \cite{beatson2010error}. Here, anisotropic scaling is used to modify the distance computation along each coordinate axis rather than the Euclidean distance \reviewB{and we set the weights $W_{0,j} = \mathbf{e}_i\mathbf{e}_i^T$ with $i\in\{1,2,3\}.$} As the scaling factors for all but one coordinate approach zero, the kernel settings collapse towards a singular limit, effectively converging to a 1D kernel. \reviewB{This adds 1D kernel features to the trial space and increases its approximation power under the same interpolation constraints. In particular, the enriched space better captures neighborhoods with strong one-dimensional structure, such as locally flat regions and plate-like parts.}

\begin{figure}
    \centering
    \begin{tabular}{cc}
    \subfloat[Original 3D dataset
    ]{\includegraphics[width=.4\textwidth]{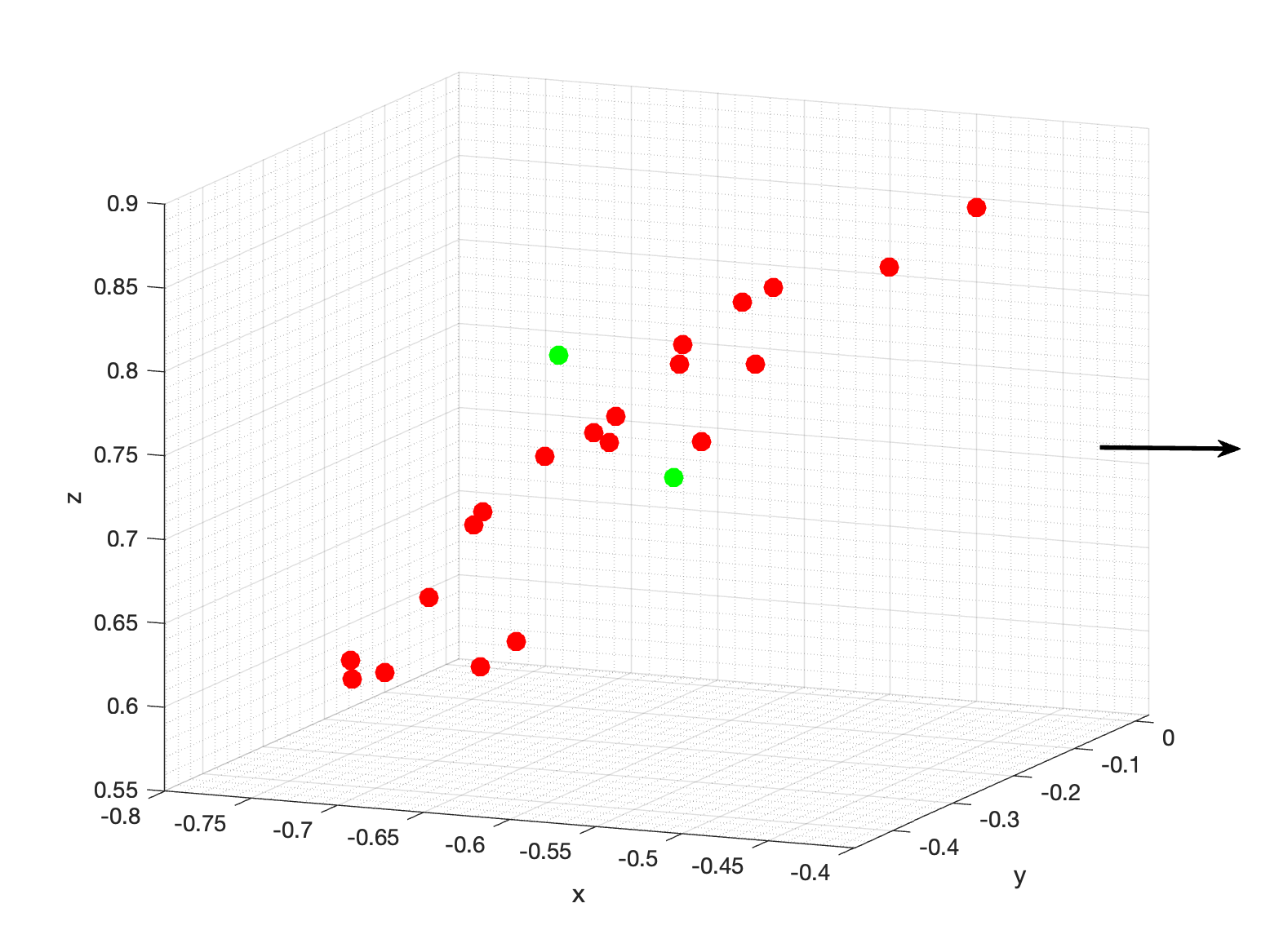}} &
    \subfloat[1D dataset of x data]{\includegraphics[width=.4\textwidth]{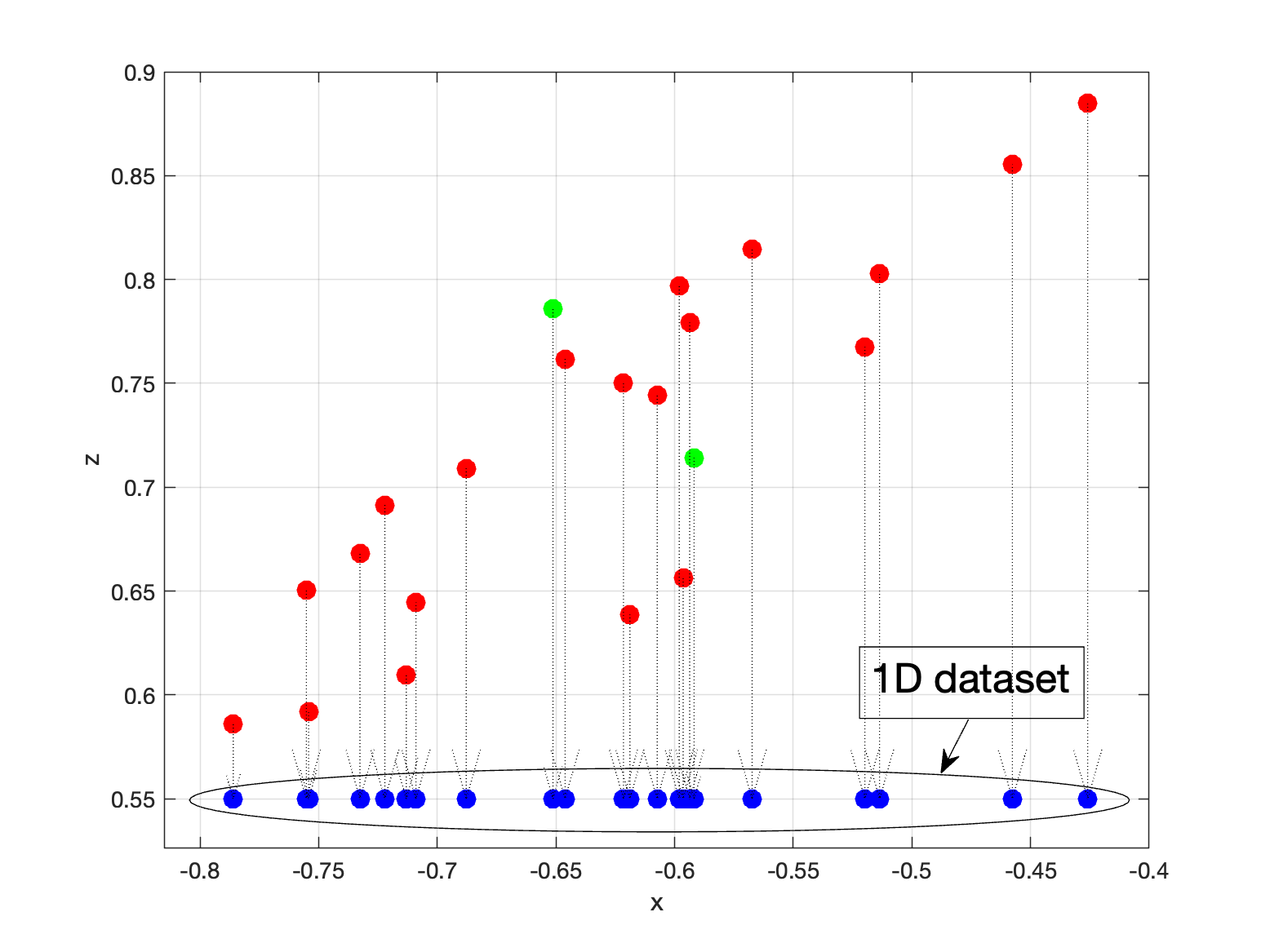}} \\
    \subfloat[1D dataset of y data]{\includegraphics[width=.4\textwidth]{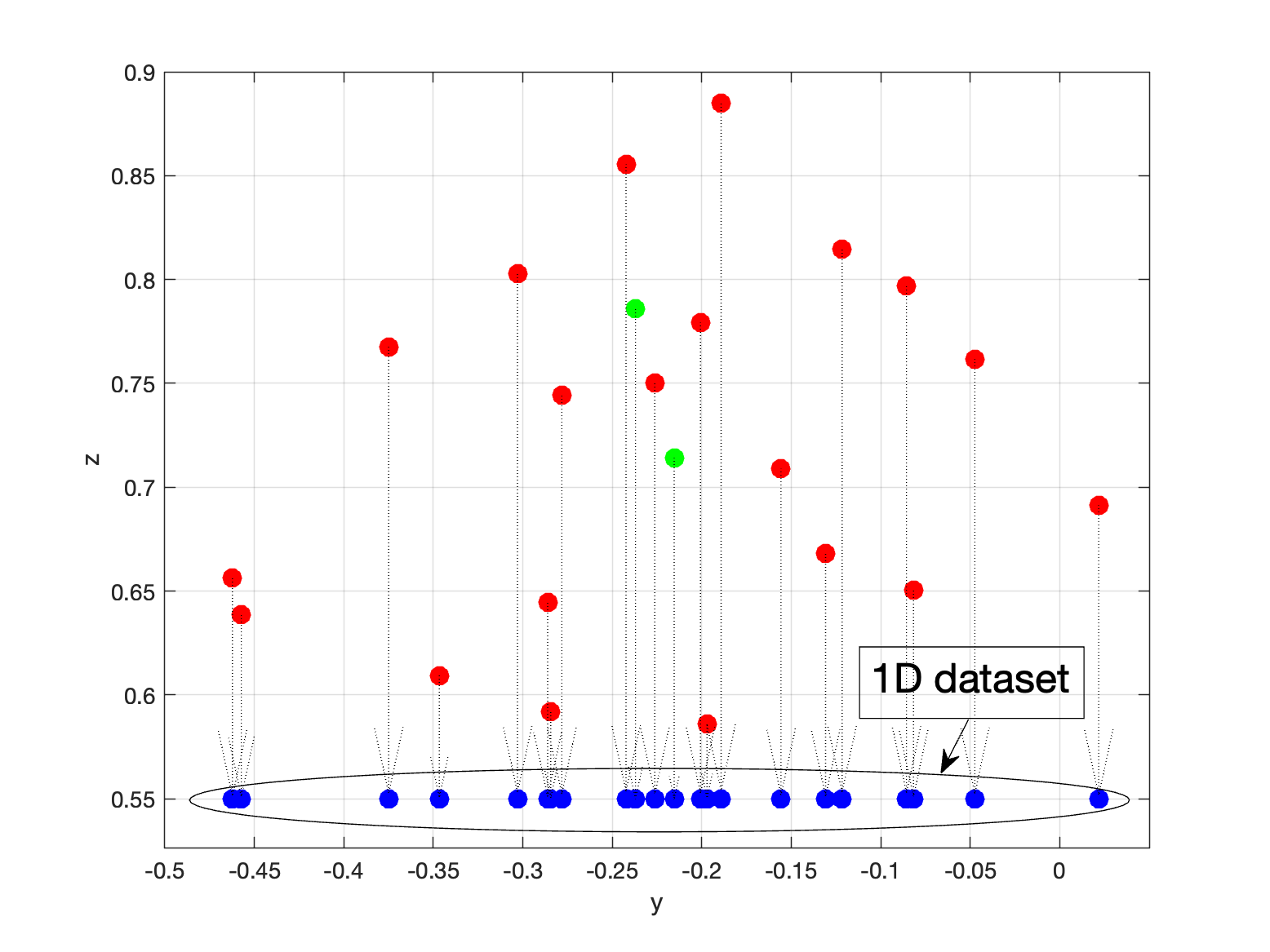}} &
    \subfloat[1D dataset of z data]{\includegraphics[width=.4\textwidth]{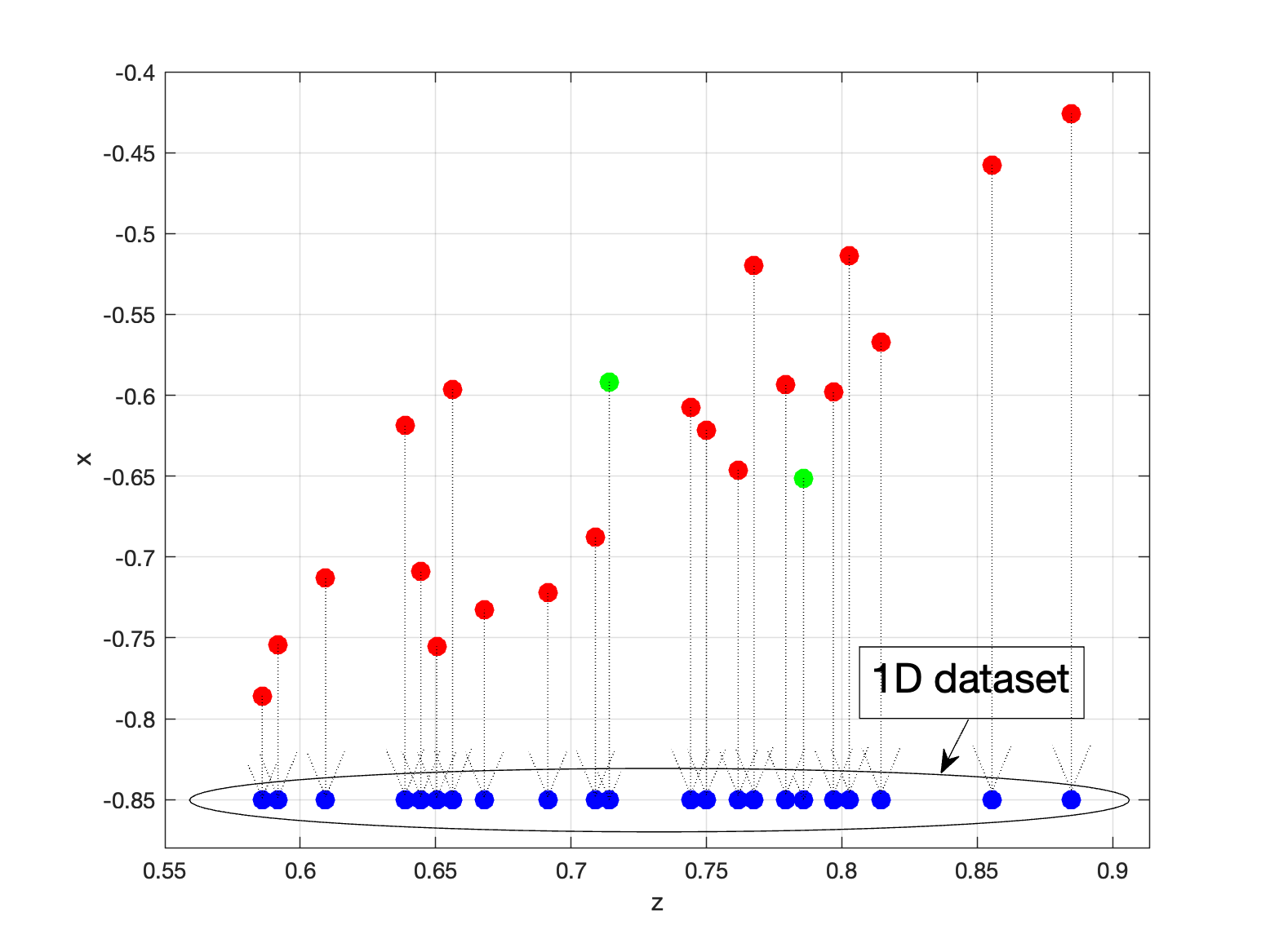}}
    \end{tabular}
    \caption{Decomposition of a 3D data set into 1D projections: (a) shows the original data in the xyz-axis. (b) through (d) illustrate projections onto the xz, yz, and zx planes, respectively, with further projections onto the x, y, and z axes.}
    \label{fig:eg_KAN}
\end{figure}

\autoref{fig:eg_KAN} demonstrates the segmentation of a 3-dimensional data set $\Xi$ into three 1-dimensional data sets. This process exemplifies the concept of degenerate anisotropic scaling, where each spatial dimension is independently handled, allowing the 3D data set to yield more basis functions.

Now, we use the sum of function spaces \( \mathcal{U}_{\text{KAN}} = \mathcal{U}_{\Xi,\Phi_{\tau,3}} + \mathcal{K}_{\Xi,\Phi_{\tau,1}} \) as the trial space for the optimization problem described in \eqref{eq:opt_problem}. The coordinates are denoted by \( \Xi = [X, Y, Z] \).
Any trial function in \( \mathcal{U}_{\text{KAN}} \) can be expressed as an \( N \times 4N \) underdetermined matrix-vector product:
\begin{equation}\label{eq:KA_RBF_K_trial_function}
F = \begin{pmatrix}
\Phi_{\tau,3}(\cdot,\Xi) &
\Phi_{\tau,1}(\cdot^T\mathbf{e}_1,X) &
\Phi_{\tau,1}(\cdot^T\mathbf{e}_2,Y) &
\Phi_{\tau,1}(\cdot^T\mathbf{e}_3,Z)
\end{pmatrix} \boldsymbol{\lambda} \in \mathcal{U}_{\text{KAN}},
\end{equation}

for any \( \boldsymbol{\lambda} \in \mathbb{R}^{4N} \). We assume unique decomposition in the sum spaces and equipped $F$ with a sum of native space norms defined by:
\begin{equation}\label{eq:KA_RBF_K_norm}
\|F\|_{\mathcal{N}_{\text{KAN}}} = \boldsymbol{\lambda}^T K\boldsymbol{\lambda} :=
\boldsymbol{\lambda}^T \text{diag} \begin{pmatrix}
\Phi_{\tau,3}(\Xi,\Xi) & \Phi_{\tau,1}(X,X) & \Phi_{\tau,1}(Y,Y) & \Phi_{\tau,1}(Z,Z)
\end{pmatrix} \boldsymbol{\lambda}.
\end{equation}
The interpolation matrix  in the constraint is given by:
\begin{equation}\label{eq:KA_RBF_K_matrix}
 A\boldsymbol{\lambda} :=
 \begin{bmatrix}
\Phi_{\tau,3}(\Xi,\Xi) & \Phi_{\tau,1}(X,X) & \Phi_{\tau,1}(Y,Y) & \Phi_{\tau,1}(Z,Z)
\end{bmatrix} \boldsymbol{\lambda}=\text{RHS}.
\end{equation}

\subsection{Hermite Radial Basis Function (HRBF)} \label{sec:HRBF}
To facilitate a numerical comparison, we introduce new basis functions linked to Hermite interpolation.
\begin{definition}[HRBF trial space]
The HRBF trial space comprises the span of the RBF and its first derivatives:
\begin{equation}
    \mathcal{H}_{\Xi,\Phi_{\tau,3}} := \text{\rm span}\Big\{\Phi_{\tau,3}(\cdot,\p_j), \nabla^{\p}\Phi_{\tau,3}(\cdot,\p_j):\p_j\in \Xi\Big\}.
    \label{eq:HRBF_trial_space}
\end{equation}
Note that each component function of $\nabla^{\p}\Phi_{\tau,3}(\cdot,\p_j)$ is used to define a basis function and the superscript $\p$ indicates the fact that the gradient operator acts upon the second variable of the kernel.
\end{definition}
A trial function in \( \mathcal{H}_{\Xi,\Phi_{\tau,3}} \) is represented as an \( N \times 4N \) underdetermined matrix-vector product:
\begin{equation}\label{eq:HRBF_K_trial_function}
F = \begin{pmatrix}
\Phi_{\tau,3}(\cdot,\Xi) & \nabla^{\p}\Phi_{\tau,3}(\cdot,\Xi)
\end{pmatrix} \boldsymbol{\lambda} \in \mathcal{H}_{\Xi,\Phi_{\tau,3}},
\end{equation}
where \( \boldsymbol{\lambda} \) is any vector in \( \mathbb{R}^{4N} \), and comes with a standard native space norm:
\begin{equation}\label{eq:HRBF_K_norm}
\|F\|_{\mathcal{N}_{\Phi_{\tau,3}}} =\boldsymbol{\lambda}^T K\boldsymbol{\lambda} :=
\boldsymbol{\lambda}^T  \begin{pmatrix}
\Phi_{\tau,3}(\Xi,\Xi) &  \nabla^{\p}\Phi_{\tau,3}(\Xi,\Xi) \\
\nabla\Phi_{\tau,3}(\Xi,\Xi)^T & \nabla\cdot \nabla^{\p}\Phi_{\tau,3}(\Xi,\Xi)
\end{pmatrix} \boldsymbol{\lambda}.
\end{equation}
The interpolation matrix in  the constraint is as follows:
\begin{equation}\label{eq:HRBF_K_matrix}
 A\boldsymbol{\lambda} :=
 \begin{pmatrix}
\Phi_{\tau,3}(\Xi,\Xi) & \nabla^{\p}\Phi_{\tau,3}(\Xi,\Xi)
\end{pmatrix} \boldsymbol{\lambda}=\text{RHS}.
\end{equation}

\section{Numerical examples} \label{sec5}

\begin{figure}
    \centering
    \subfloat[Ellipsoid]{\includegraphics[width=.33\textwidth]{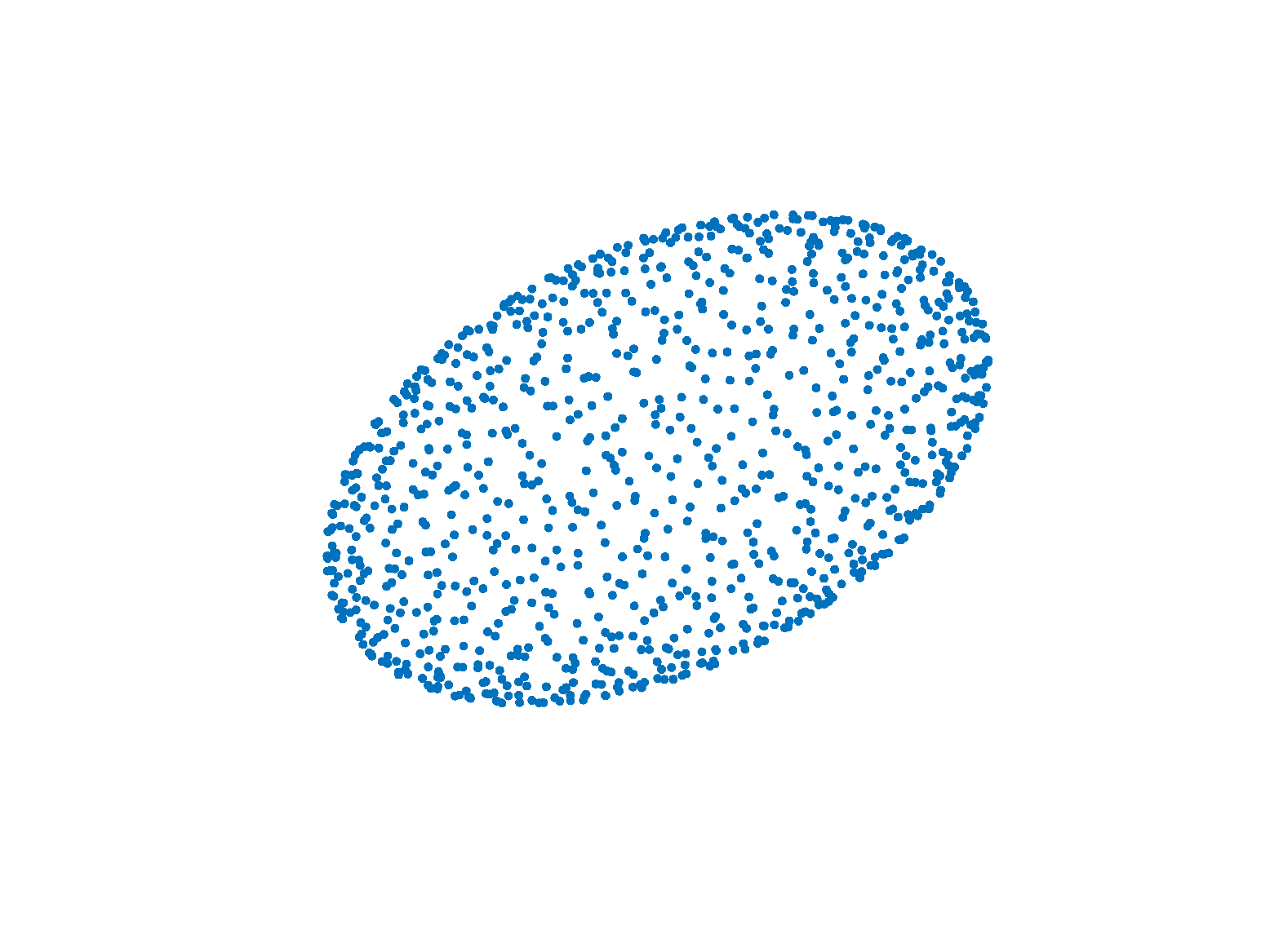}}
    \subfloat[Torus]{\includegraphics[width=.33\textwidth]{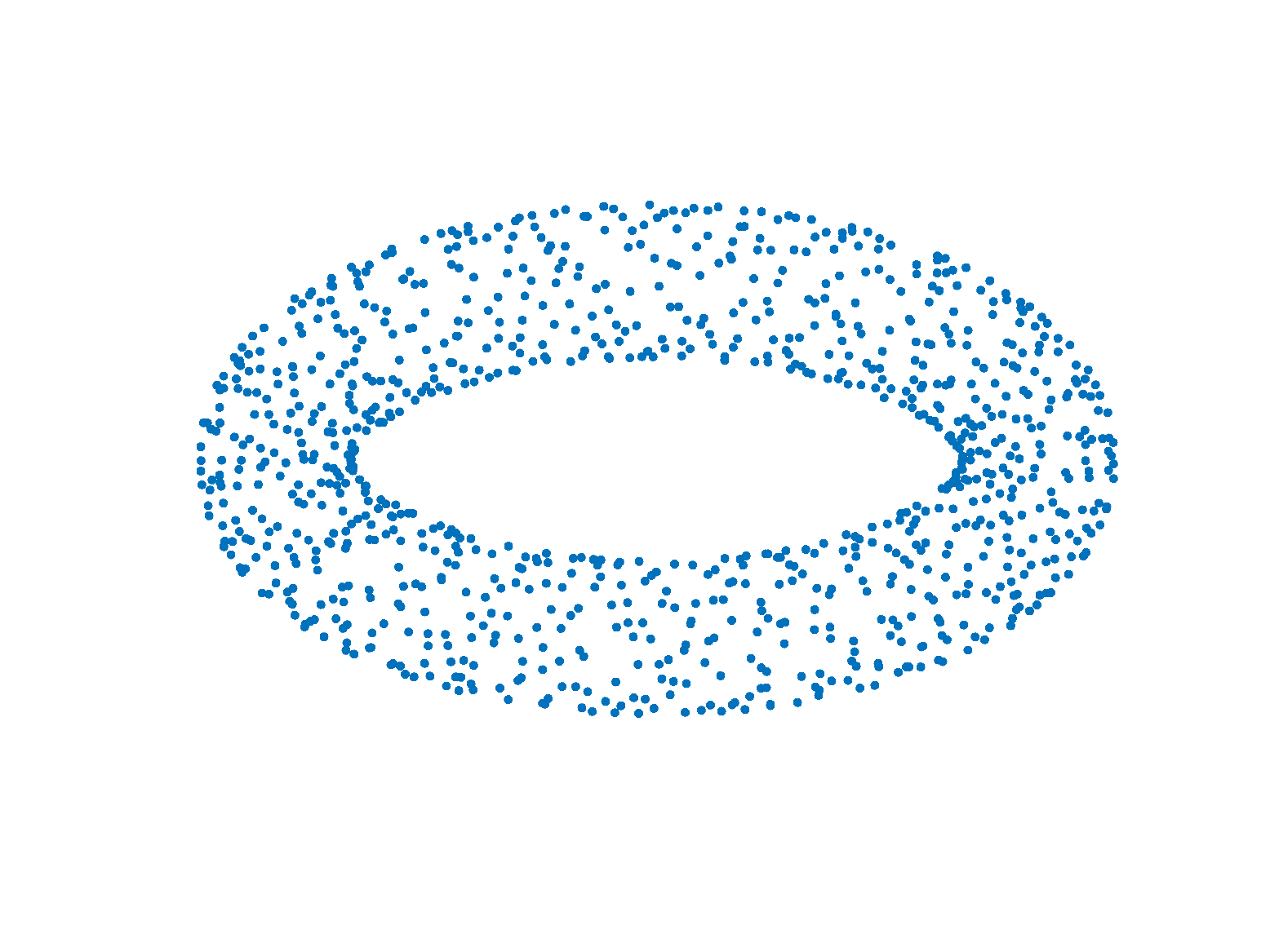}}
    \subfloat[Sphube]{\includegraphics[width=.33\textwidth]{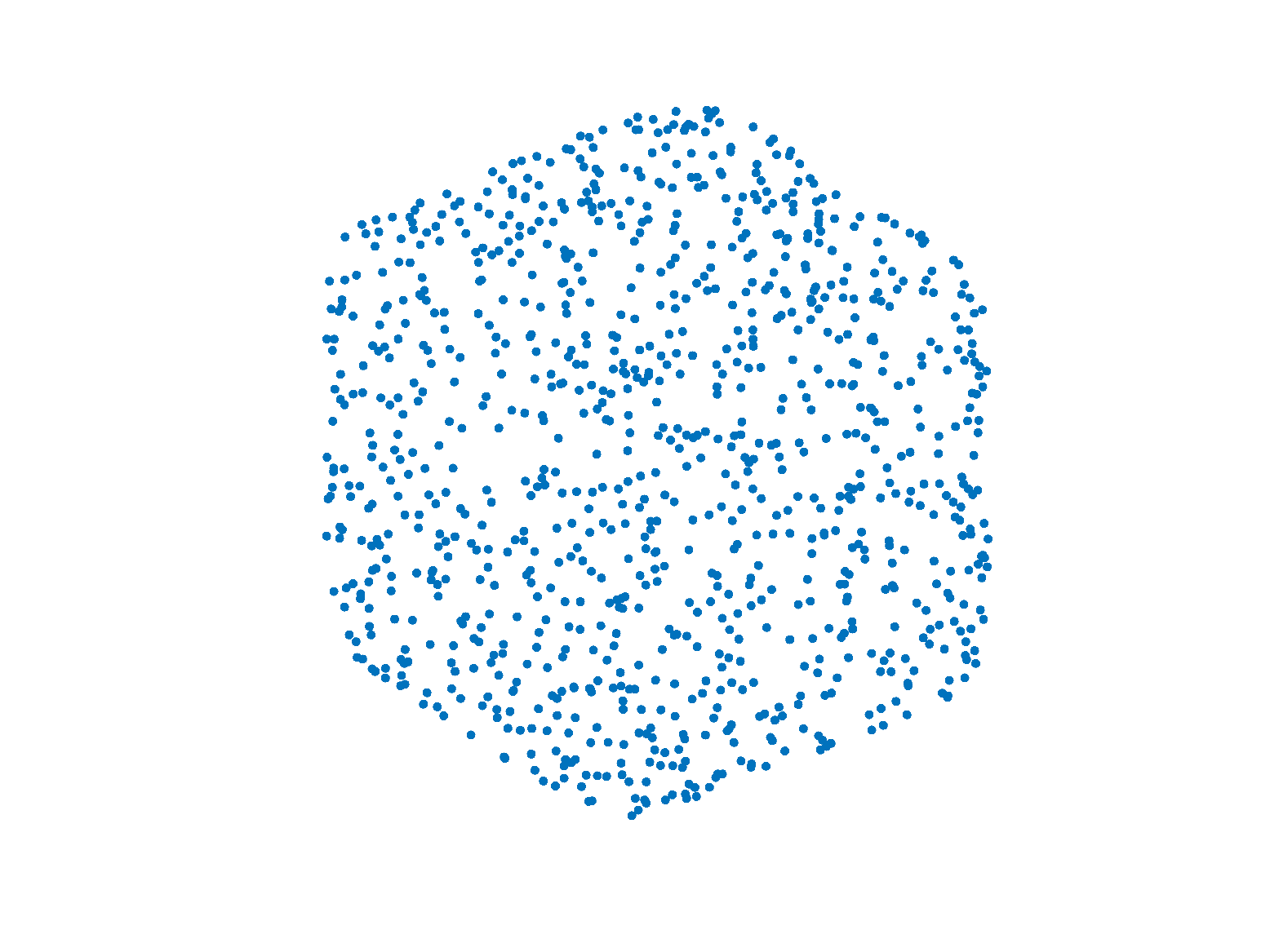}}
    \caption{Point clouds of $N=1000$ points for an ellipsoid, a torus, and a sphube. These point clouds are generated using parametric equations with points allocated through a Halton distribution to ensure low discrepancy.}
    \label{fig:eg_ptcld}
\end{figure}

This section presents a series of numerical experiments designed to assess the accuracy of our proposed interpolation method by the constrained optimization in \eqref{eq:opt_problem} posed on the sum of the traditional RBF trial space and the KAN-inspired space in \eqref{eq:KA_RBF_trial_space}. The experiments are structured to address specific aspects of the methodology, including the impact of trial center placements within the KAN layer, the computational trade-offs between using different norms in the constraints, and the effects of kernel smoothness, stencil sizes, and increased curvature on solution accuracy. Additionally, we investigate the performance of point cloud processing, comparing curvature estimations derived from our proposed method and the Tikhonov regularization method within the same trial space.
Three types of point clouds are utilized in our numerical experiments: an ellipsoid $\mathbb{E}$\footnote{$\mathbb{E} = \left\{(x,y,z)\in\mathbb{R}^3 \mid \frac{x^2}{a^2}+\frac{y^2}{b^2}+\frac{z^2}{c^2}-1=0\right\},$}, a torus $\mathbb{T}$\footnote{$\mathbb{T} = \left\{(x,y,z)\in\mathbb{R}^3 \mid (R-\sqrt{x^2+y^2})^2+z^2-r^2=0\right\},$}, and a sphube $\mathbb{S}_{s,r}$\footnote{$\mathbb{S}_{s,r} = \left\{(x,y,z)\in\mathbb{R}^3 \mid x^2+y^2+z^2-\frac{s^2}{r^2}x^2y^2-\frac{s^2}{r^2}y^2z^2-\frac{s^2}{r^2}x^2z^2 + \frac{s^4}{r^4}x^2y^2z^2 - r^2 = 0\right\},$}. The parameters $a$, $b$, and $c$ represent the lengths of the semi-axes of $\mathbb{E}$, $R$ and $r$ is the major and minor radius of $\mathbb{T}$, $s$ is the squareness parameter of $\mathbb{S}_{s,r}$.
While $\mathbb{S}_{0,r}$ resembles a sphere with radius $r$, $\mathbb{S}_{1,r}$ approaches the shape of a cube with an edge length of $2r$.
The exact surface normals for each shape can be obtained by computing the gradient of their respective implicit equations. Examples of a point cloud for an ellipsoid, a torus, and a sphube, each containing $N = 1000$ points, are depicted in \autoref{fig:eg_ptcld}. These point clouds are constructed using a Halton distribution in the parameter space to allocate points such that the points in the dataset exhibit low discrepancy.

For reproducibility, we adopt the method of using one stencil per point of interest to negate the impact of patch geometry variations. Unless specified differently, we employ the Mat\'ern Sobolev kernel with a smoothness order of $\tau=5$ and a local stencil size of $N_s=40$. The interpolation data $\Xi=\tilde{P}_0$ includes two ghost points as detailed in \eqref{P+2}.
The error function for estimating the normal vector is defined using the Euclidean norm of the difference between the exact  and the approximated normal, i.e.,
\begin{equation}
\mathcal{E}_{\n}(\p) := \|\vec{n}^*(\p) - \vec{n}(\p)\|_2
\quad\text{for }\p\in\M,
\label{eq:error_normal}
\end{equation}
and we report the maximum error $\| \mathcal{E}_{\n}\|_{\infty}$.

\subsection*{Example 1: Distribution of Trial Centers within KAN's First Layer}
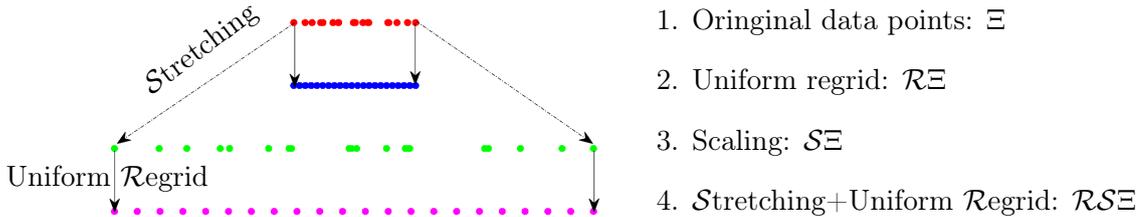
\begin{figure}
    \centering
    \hspace*{-4cm}
    \begin{overpic}[width=0.5\textwidth, trim=100 100 0 100, clip=true,tics=10]{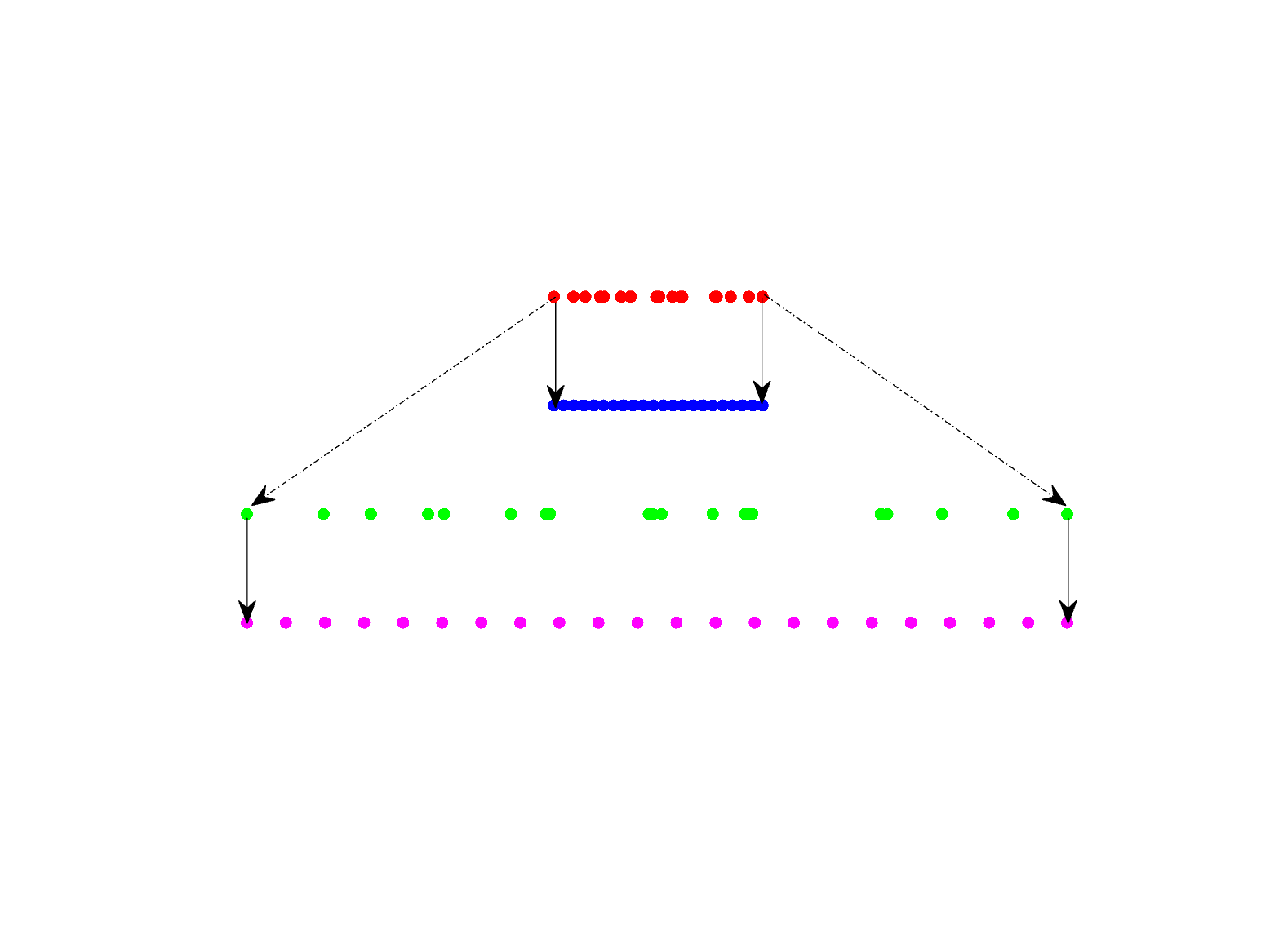}
    \put(6,29){\rotatebox{37}{\small{$\mathcal{S}$tretching }}}
    \put(-16,14){\small{Uniform $\mathcal{R}$egrid}}
    \put(91,40){\small{1. Original data points: $\Xi$}}
    \put(91,30){\small{2. Uniform regrid: $\mathcal{R}\Xi$}}
    \put(91,20){\small{3. Scaling: $\mathcal{S}\Xi$ }}
    \put(91,10){\small{4. $\mathcal{S}$tretching+Uniform $\mathcal{R}$egrid: $\mathcal{R}\mathcal{S}\Xi$}}
    \end{overpic}
    \caption{Example 1: Trial Center Distributions in KAN's First Layer: Config 1 uses the original data points, $\Xi$, as trial centers. Config 2 applies a uniform regrid, transforming $\Xi$ to $\mathcal{R}\Xi$. Config 3 stretches $\Xi$ to $\mathcal{S}\Xi$, adjusting the interval around the data center. Config 4 combines both stretching and uniform regrid techniques.}
    \label{fig:local_stencil_distribution}
\end{figure}

In our proposed trial space, as discussed in \autoref{sec:KA-RBF}, the method requires projections of a 3D dataset into three distinct 1D datasets. This decomposition may  lead to clustering  within the data. To address this, we explore different configurations of the trial center distribution within our KAN-inspired trial space approach. This investigation aims to evaluate the impact on point cloud surface reconstruction performance. Four configurations are designed to assess both the robustness and accuracy of the method:
\begin{itemize}
  \item Config 1 utilizes the original data points, $\Xi$, as trial centers, following the definition provided in \eqref{eq:KA_RBF_trial_space}. This setup serves as a benchmark.
  \item Config 2 implements a uniform regrid, transforming $\Xi$ to $\mathcal{R}\Xi$. This modification redistributes the data points into a uniform distribution, aiming to maximize the minimum separating distance.
  \item Config 3 involves stretching $\Xi$ to $\mathcal{S}\Xi$, expanding the interval around the data center to a specific reference length. This configuration standardizes the range of trial centers across stencils.
  \item Config 4 integrates both the stretching and uniform regrid techniques.
\end{itemize}
See \autoref{fig:local_stencil_distribution} for a pictorial overview of these configurations.

The performance of different trial center distributions for the minimum-norm interpolation on $\mathcal{U}_{\text{KAN}}$ is evaluated using an ellipsoid with semi-axes $a=0.85$, $b=0.35$, and $c=0.5$. The tests involve estimating surface normals on this ellipsoid with datasets ranging from $N =100$ to 5000 points in total. The experiments are conducted for varying local stencil sizes  and the behavior is assessed relative to the smoothness order $2\leq \tau \leq 5$ of the kernel. The maximum error across stencil sizes $40\leq N_s\leq80$ is reported in \autoref{table:KnR_config}.

In terms of normal estimation accuracy, all configurations exhibit relatively similar performance across the tested range of stencil sizes and kernel smoothness orders from $\tau = 2$ to $\tau = 5$. Notable exceptions are highlighted in boldface in the table, where some configurations significantly outperform the benchmark Config 1 by an order of magnitude. This superior performance is usually observed for Config 4 when the total number of data points is large which is consistent with the effect of enforcing a common coordinate range across stencils and improving separation in the 1D projected centers. Based on these results, Config 4 is selected for the remainder of the experiments.
\reviewB{
\begin{remark}
    In practice, we suggest using a simple indicator such as the local mesh ratio $\rho:=h/q$ where $h$ is the fill distance and $q:=\min_{i\neq j}\|\p_i-\p_j\|_2$ is the separation distance to determine which configuration to use. In our experience, Config 4 is most reliable when $\rho$ is moderate so the stencil is not overly sparse or irregular. When $\rho$ is large, the benefit becomes less consistent and increasing the stencil density can improve its performance.
\end{remark}
}

\subsection*{Example 2: Selection of the Norm to be Optimized}

Solving the proposed minimum-norm interpolation problem requires addressing an optimization problem in matrix-vector form:
\[
        \text{minimize } \|K^{1/2}\boldsymbol\lambda\|_2 \quad \text{subject to } A\boldsymbol\lambda = b,
\]
see \eqref{eq:KA_RBF_K_norm} and \eqref{eq:HRBF_K_norm} for $K$, and \eqref{eq:KA_RBF_K_matrix} and \eqref{eq:HRBF_K_matrix} for $A$ of the KAN-inspired and HRBF trial spaces.
To solve this, one can first compute the square root or Cholesky factorization of $K$. Then, solve the modified system
\[
(AK^{-1/2})(K^{1/2}\boldsymbol\lambda) = b
\]
for $(K^{1/2}\boldsymbol\lambda)$ in the minimum-norm sense. Finally, solve
\[
K^{1/2}\boldsymbol\lambda = (K^{1/2}\boldsymbol\lambda)
\]
to recover $\boldsymbol\lambda$. The leading computational cost in this process is the factorization of $K$, especially when $K$ is close to singular. 
\reviewB{Compared with standard RBF interpolation, which solves a dense $N\times N$ system per stencil, the enriched trial spaces lead to an underdetermined constraint $A\boldsymbol\lambda=b$ with $A\in\mathbb{R}^{N\times 4N}$ and $\boldsymbol\lambda\in\mathbb{R}^{4N}$. As a result, the overall scaling remains cubic in $N$ but with a larger constant that depends on the chosen linear algebra routine. In our experiments, we apply a dense Cholesky factorization to $K\in\mathbb{R}^{4N\times 4N}$, so this step scales as $O((4N)^3)=64\,O(N^3)$.}
To reduce computational cost, one may instead minimize $\|\boldsymbol\lambda\|_{\ell^2}$ instead of the native space norm. In such cases, $K$ can be replaced by the identity matrix $I$, reducing the problem to the standard minimum-norm least-squares solution for a linear system.

In \autoref{table:opt_norm_choice}, we report the maximum error in normal estimation based on the minimum-norm interpolant for the KAN-inspired and HRBF trial spaces, considering both the native space norm and the $\|\boldsymbol\lambda\|_{\ell^2}$ norm. The setup is identical to Example 1, with tests conducted for varying local stencil sizes $40\leq N_s\leq80$, and the maximum error is presented.
The results indicate that using the KAN-inspired trial space consistently achieves higher accuracy compared to the HRBF trial space. For small $h$, the accuracy difference between the two trial spaces can reach 2 to 3 orders of magnitude. In contrast, the choice between minimizing the native space norm or the $\|\boldsymbol\lambda\|_{\ell^2}$ norm shows minimal impact on performance in this test.

\subsection*{Example 3: Convergence Studies on Torus and Sphube }
\begin{figure}
    \centering
    \subfloat[$\tau=3 $]{
    \begin{overpic}[width=0.3\textwidth, trim=0 0 0 0, clip=true,tics=10]{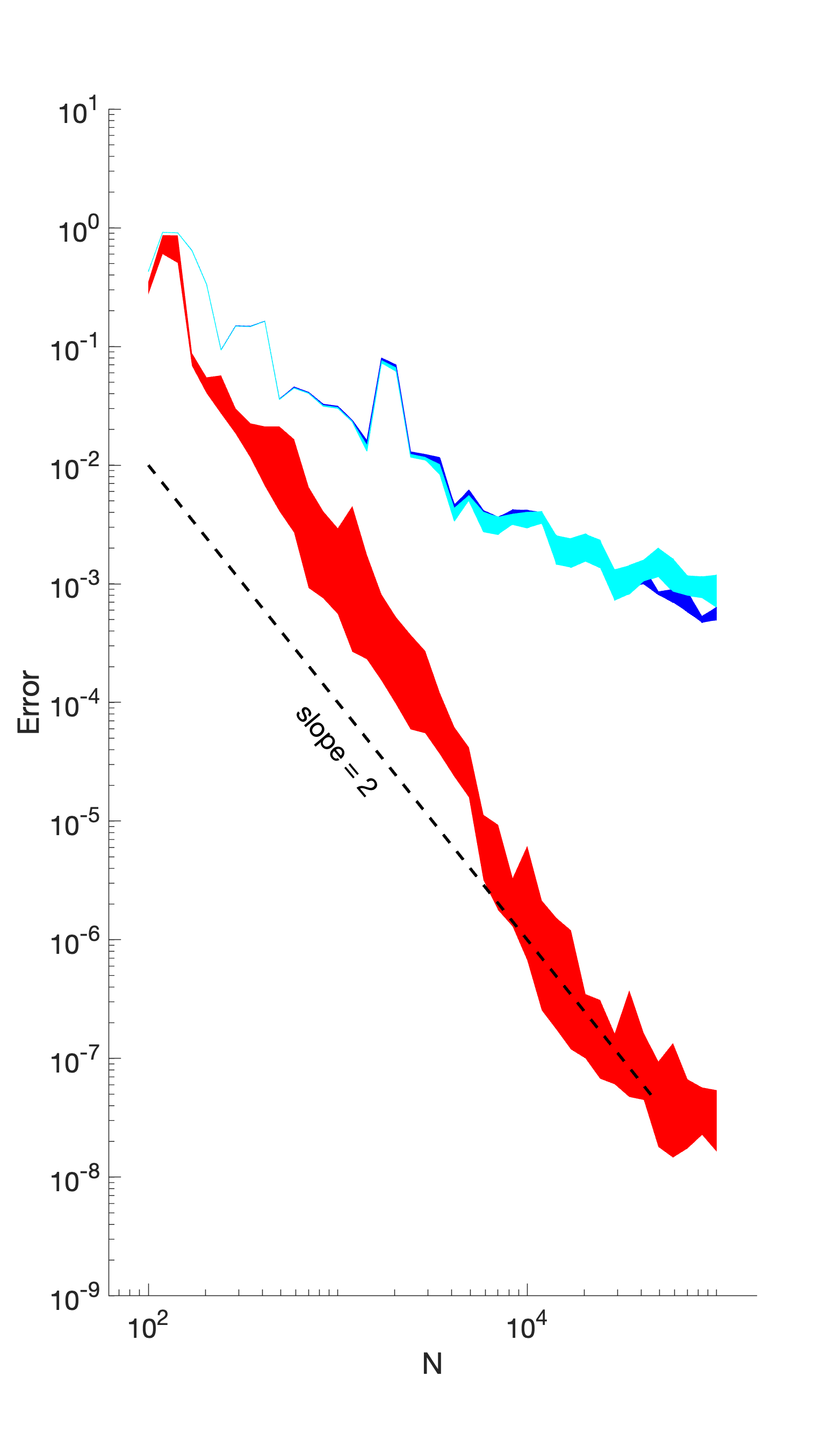}
        \put(35,85){\csq{blue}} \put(40,85){\tiny{RBF}}
        \put(35,82){\csq{Aqua}} \put(40,82){\tiny{HRBF}}
        \put(35,79){\csq{red}}  \put(40,79){\tiny{KRBF}}
        \put(34,78){\linethickness{0.25mm}\color{black}\polygon(0,0)(20,0)(20,10)(0,10)}
    \end{overpic}}
    \subfloat[$\tau=4 $]{
    \begin{overpic}[width=0.3\textwidth, trim=0 0 0 0, clip=true,tics=10]{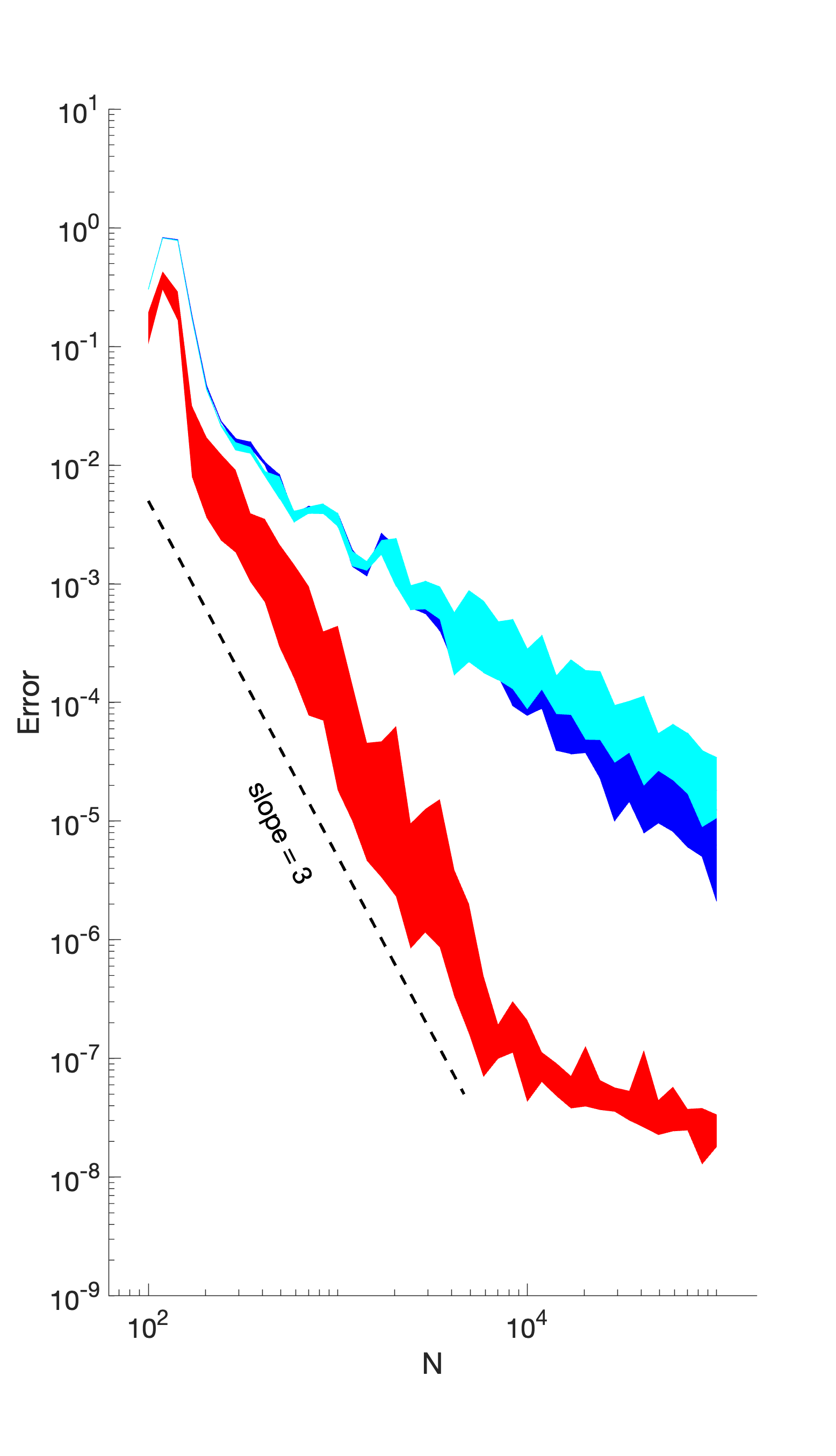}
        \put(35,85){\csq{blue}} \put(40,85){\tiny{RBF}}
        \put(35,82){\csq{Aqua}} \put(40,82){\tiny{HRBF}}
        \put(35,79){\csq{red}}  \put(40,79){\tiny{KRBF}}
        \put(34,78){\linethickness{0.25mm}\color{black}\polygon(0,0)(20,0)(20,10)(0,10)}
    \end{overpic}}
    \subfloat[$\tau=5  $]{
    \begin{overpic}[width=0.3\textwidth, trim=0 0 0 0, clip=true,tics=10]{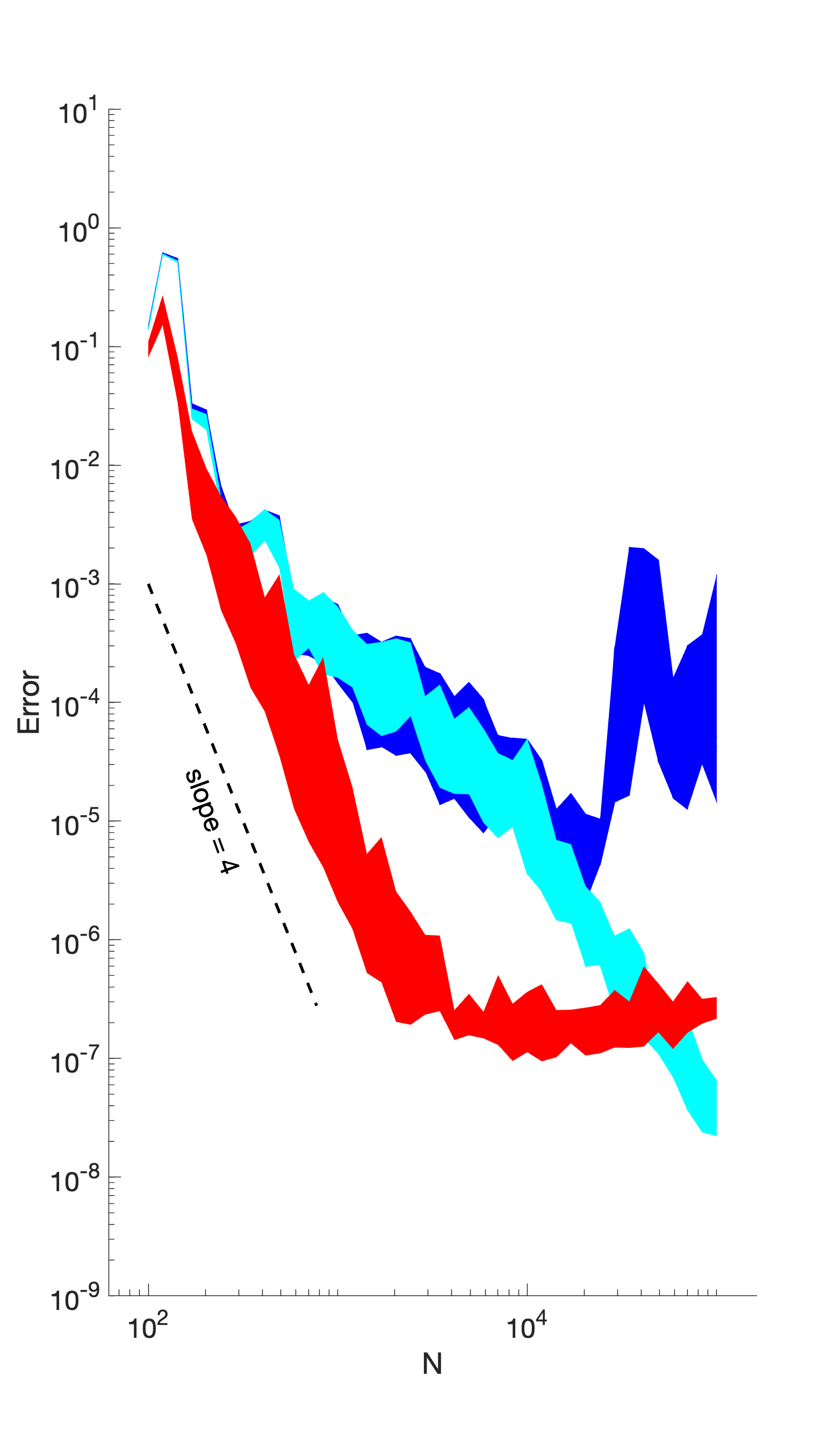}
        \put(35,85){\csq{blue}} \put(40,85){\tiny{RBF}}
        \put(35,82){\csq{Aqua}} \put(40,82){\tiny{HRBF}}
        \put(35,79){\csq{red}}  \put(40,79){\tiny{KRBF}}
        \put(34,78){\linethickness{0.25mm}\color{black}\polygon(0,0)(20,0)(20,10)(0,10)}
    \end{overpic}}
    \caption{Example 3: Convergence study on a torus for normal estimation using various methods with a Sobolev kernel of smoothness order ${\tau}$ and stencil sizes ranging from ${40 \leq N \leq 90}$.}
    \label{fig:convergence_Torus}
    \centering
    \subfloat[$\tau=3 $]{
    \begin{overpic}[width=0.3\textwidth, trim=0 0 0 0, clip=true,tics=10]{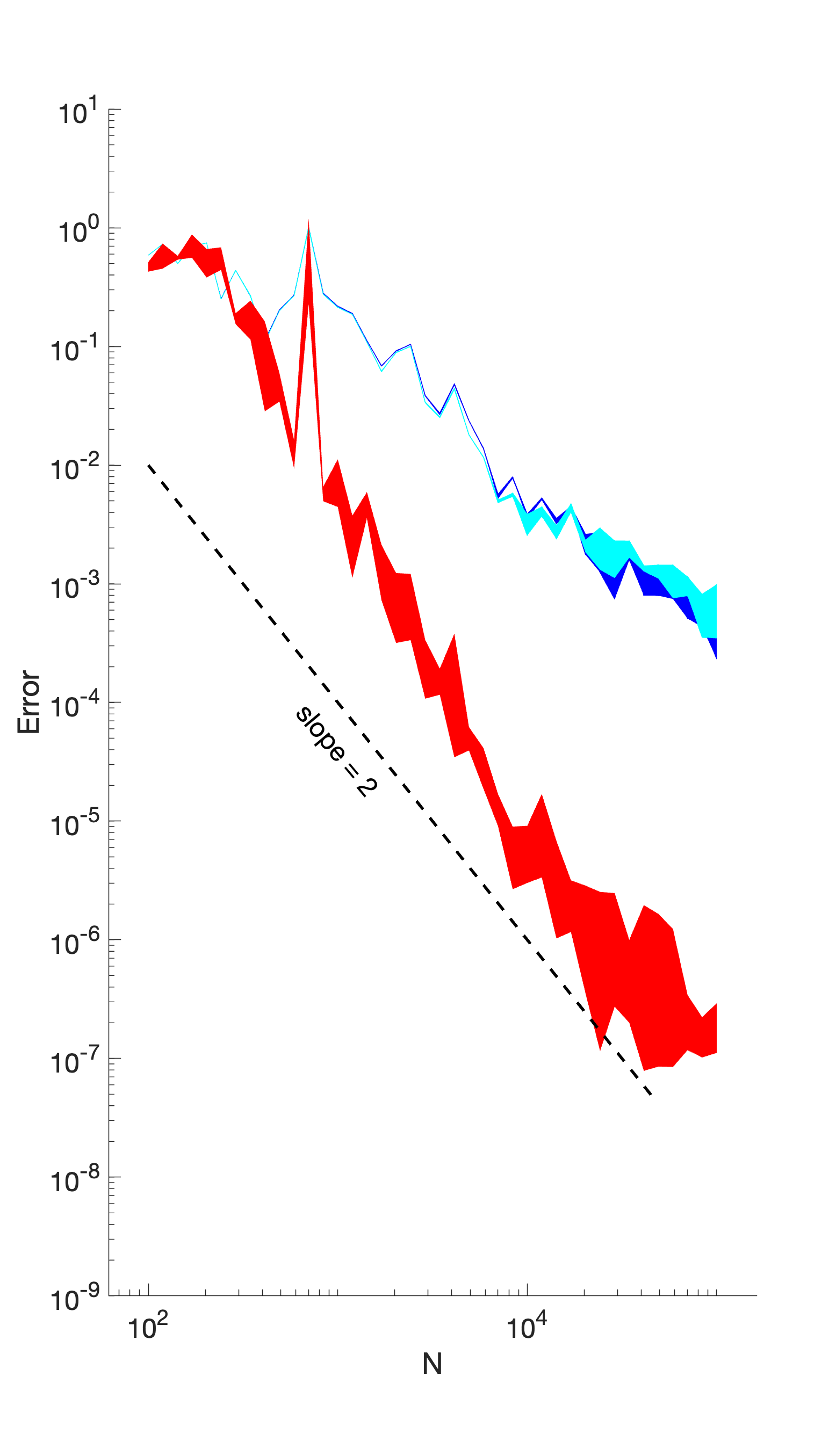}
        \put(35,85){\csq{blue}} \put(40,85){\tiny{RBF}}
        \put(35,82){\csq{Aqua}} \put(40,82){\tiny{HRBF}}
        \put(35,79){\csq{red}}  \put(40,79){\tiny{KRBF}}
        \put(34,78){\linethickness{0.25mm}\color{black}\polygon(0,0)(20,0)(20,10)(0,10)}
    \end{overpic}}
    \subfloat[$\tau=4 $]{
    \begin{overpic}[width=0.3\textwidth, trim=0 0 0 0, clip=true,tics=10]{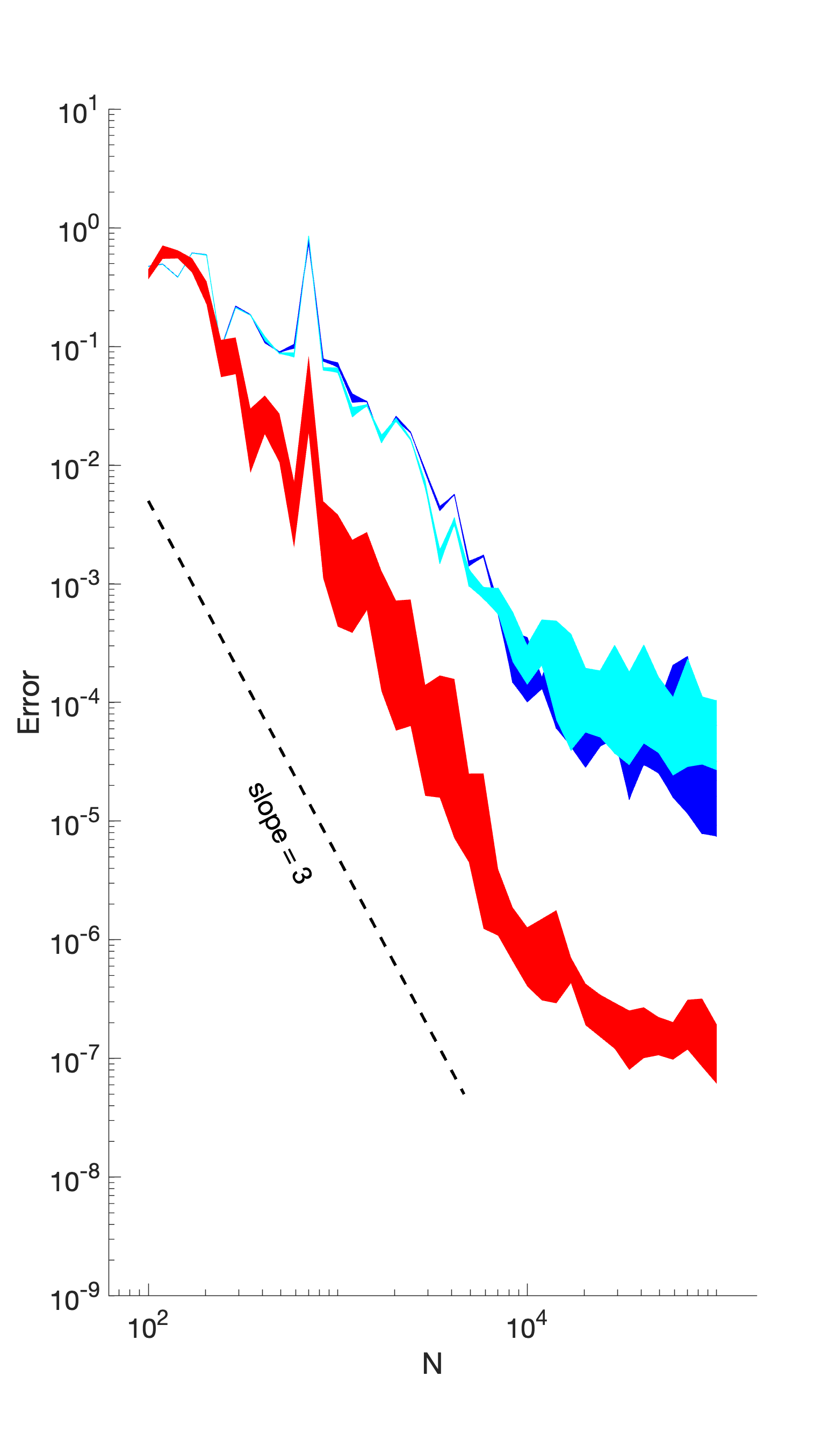}
        \put(35,85){\csq{blue}} \put(40,85){\tiny{RBF}}
        \put(35,82){\csq{Aqua}} \put(40,82){\tiny{HRBF}}
        \put(35,79){\csq{red}}  \put(40,79){\tiny{KRBF}}
        \put(34,78){\linethickness{0.25mm}\color{black}\polygon(0,0)(20,0)(20,10)(0,10)}
    \end{overpic}}
    \subfloat[$\tau=5 $]{
    \begin{overpic}[width=0.3\textwidth, trim=0 0 0 0, clip=true,tics=10]{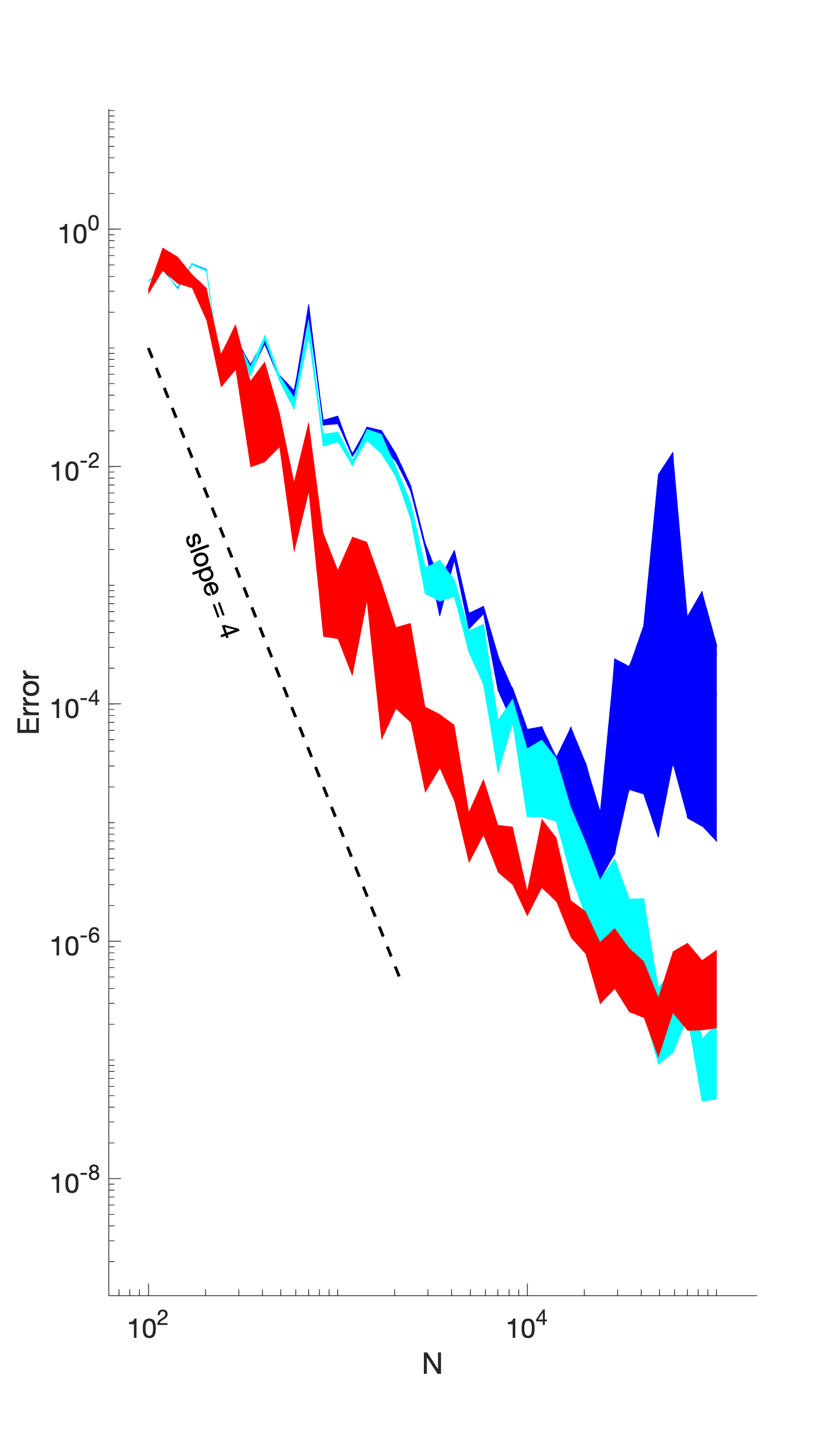}
        \put(35,85){\csq{blue}} \put(40,85){\tiny{RBF}}
        \put(35,82){\csq{Aqua}} \put(40,82){\tiny{HRBF}}
        \put(35,79){\csq{red}}  \put(40,79){\tiny{KRBF}}
        \put(34,78){\linethickness{0.25mm}\color{black}\polygon(0,0)(20,0)(20,10)(0,10)}
    \end{overpic}}
    \caption{Example 3: Convergence study on a sphube for normal estimation using various methods with a Sobolev kernel of smoothness order ${\tau}$ and stencil sizes ranging from ${40 \leq N \leq 90}$.}
    \label{fig:convergence_Sphube}
\end{figure}

In this example, we perform a convergence study for surface normal estimation using the following interpolants:
\begin{itemize}
  \item \textbf{RBF} is the standard interpolant from $\mathcal{U}_{\Xi,\Phi}$ in \eqref{eq:RBF_trial_space} by solving an $N_s\times N_s$ interpolation matrix system.
  \item \textbf{HRBF} is the min-norm interpolant in $ \mathcal{H}_{\Xi,\Phi_{\tau,3}}$ that minimizes native space norm in \eqref{eq:HRBF_K_norm}.
  \item \textbf{KRBF} is the min-norm interpolant in $ \mathcal{U}_{\text{KAN}} = \mathcal{U}_{\Xi,\Phi_{\tau,3}} + \mathcal{K}_{\Xi,\Phi_{\tau,1}}$ that minimizes \eqref{eq:KA_RBF_K_norm}.
\end{itemize}
We conducted experiments on a torus with major radius \( R=1 \) and minor radius \( r=0.2 \), and on a sphube with radius \( r=1 \) and squareness parameter \( s=0.9 \). Tests utilized Sobolev kernels with smoothness orders ranging from \( 3 \leq \tau \leq 5 \), and local stencils of size \( 40 \leq N_s \leq 90 \). A reference slope of \(\tau-1\) was included for comparative analysis. Results are presented in a collective format, displaying errors for all values of \(\tau\) and \(N\) as shaded areas in \autoref{fig:convergence_Torus} for the torus and \autoref{fig:convergence_Sphube} for the sphube.

The KRBF generally exhibited superior performance over both standard RBF and HRBF, with exceptions occurring at \(\tau = 5\) and larger \(N\) values (i.e., small \(h\)), where convergence stagnated because the Cholesky decomposition failed for nearly singular matrices \(K\). In such cases, we stabilize the computation by replacing $K$ with $K+\epsilon I$.  Notably, KRBF demonstrated its most significant advantages at the lower smoothness order \(\tau=3\), which shows the same order of convergence and accuracy as its \(\tau=4\) counterpart. Further analysis is needed to determine if the property \(\Phi_{\tau+1,3}=\Phi_{\tau,1}\) plays a role in this faster than expected convergence.

\subsection*{Example 4: Flattening Sphube with Increasing Curvature}

Previous examples primarily focused on relatively smooth surfaces. To further evaluate the limits of the proposed approach, we conduct tests on a sphube with radius \( r = 1 \) and squareness parameter \( s \) ranging from 0.1 to 0.9 in increments of 0.2.

\begin{figure}
    \centering
    \vspace{1mm}
    \begin{overpic}[width=0.19\textwidth, trim=60 60 120 40, clip=true,tics=10]{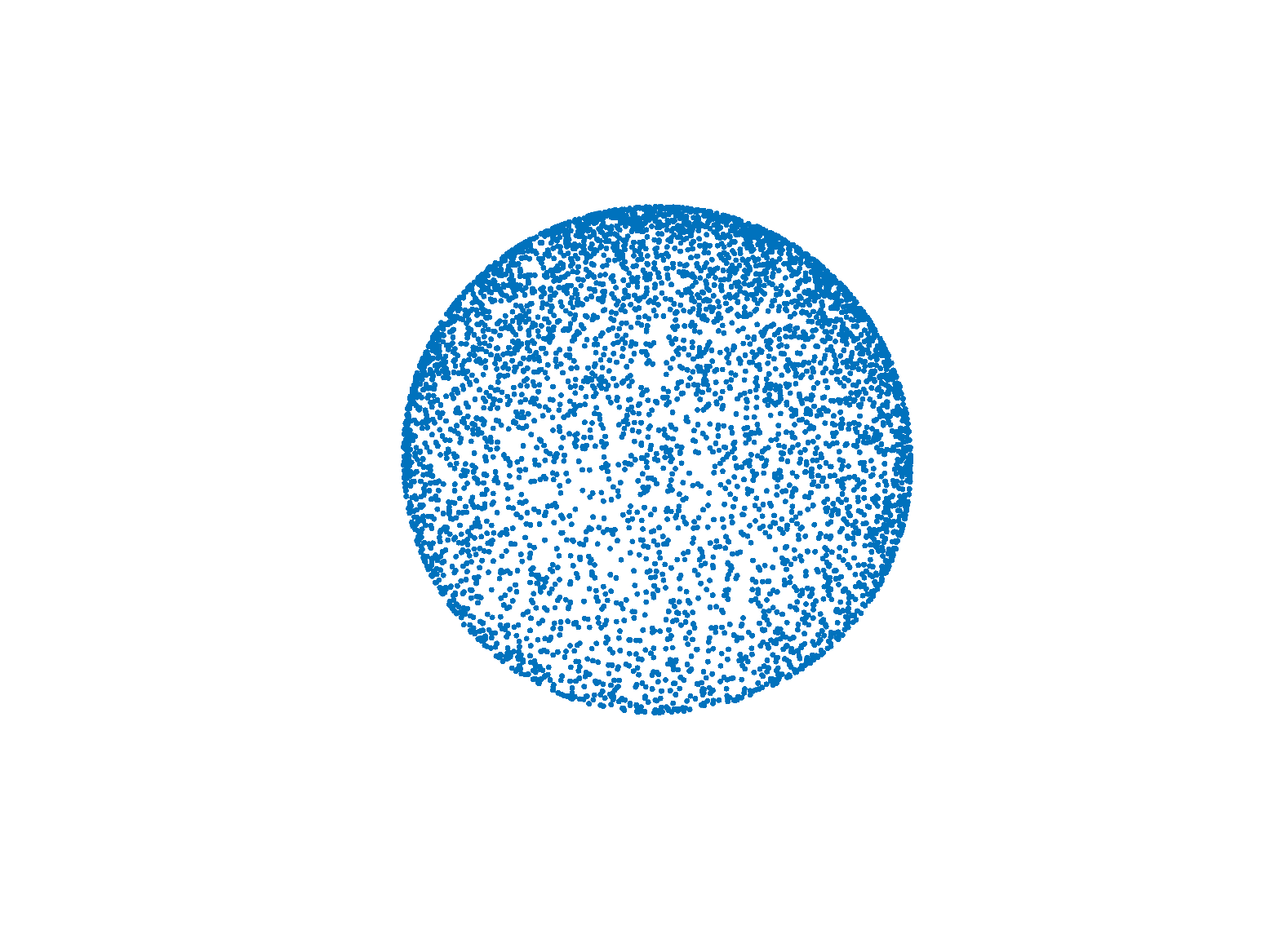}
    \put(38,86){$\bm{s = 0.1}$}
    \end{overpic}
    \vspace{1mm}
    \begin{overpic}[width=0.19\textwidth, trim=60 60 120 40, clip=true,tics=10]{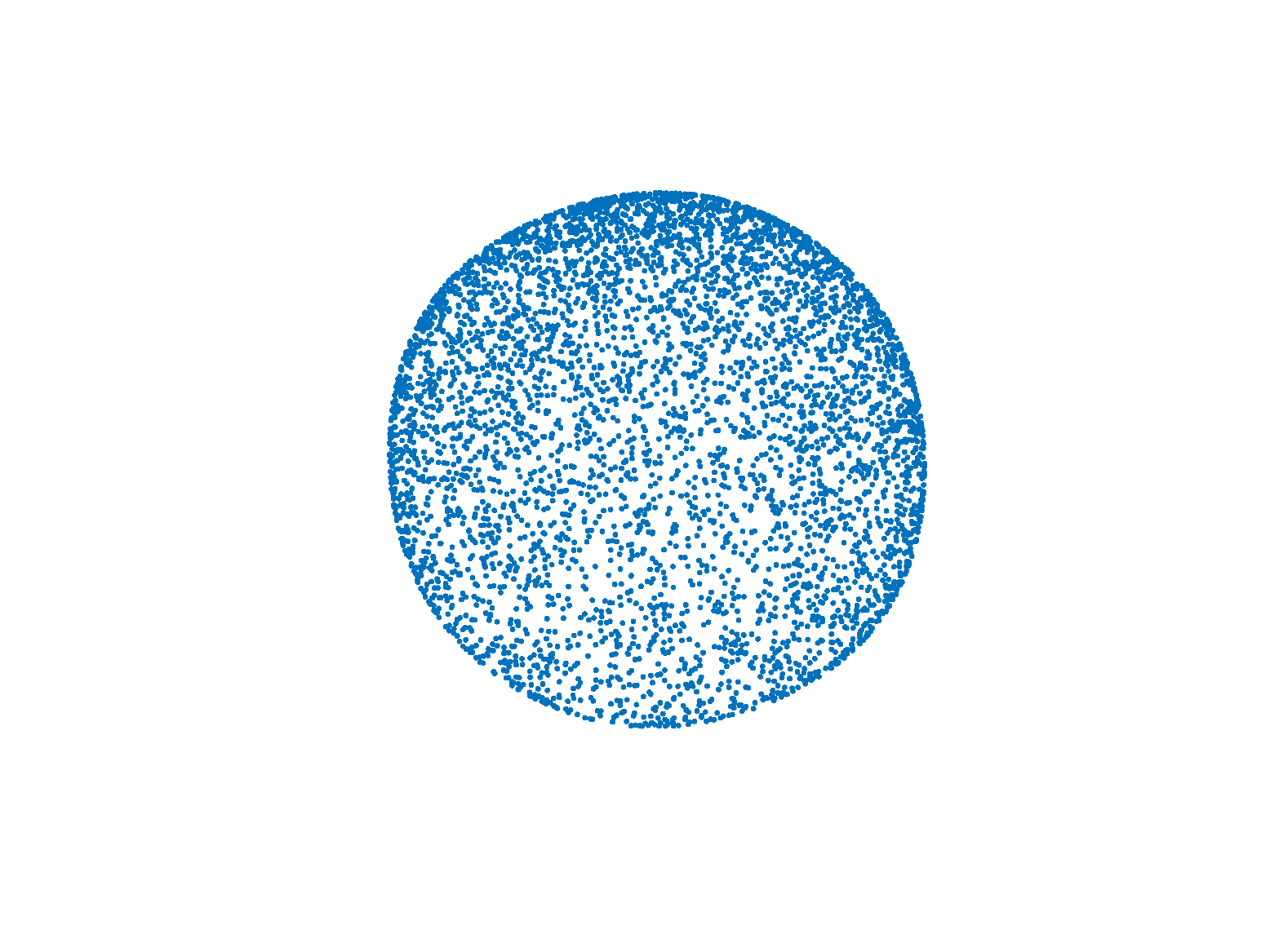}
    \put(38,86){$\bm{s = 0.3}$}
    \end{overpic}
    \vspace{1mm}
    \begin{overpic}[width=0.19\textwidth, trim=60 60 120 40, clip=true,tics=10]{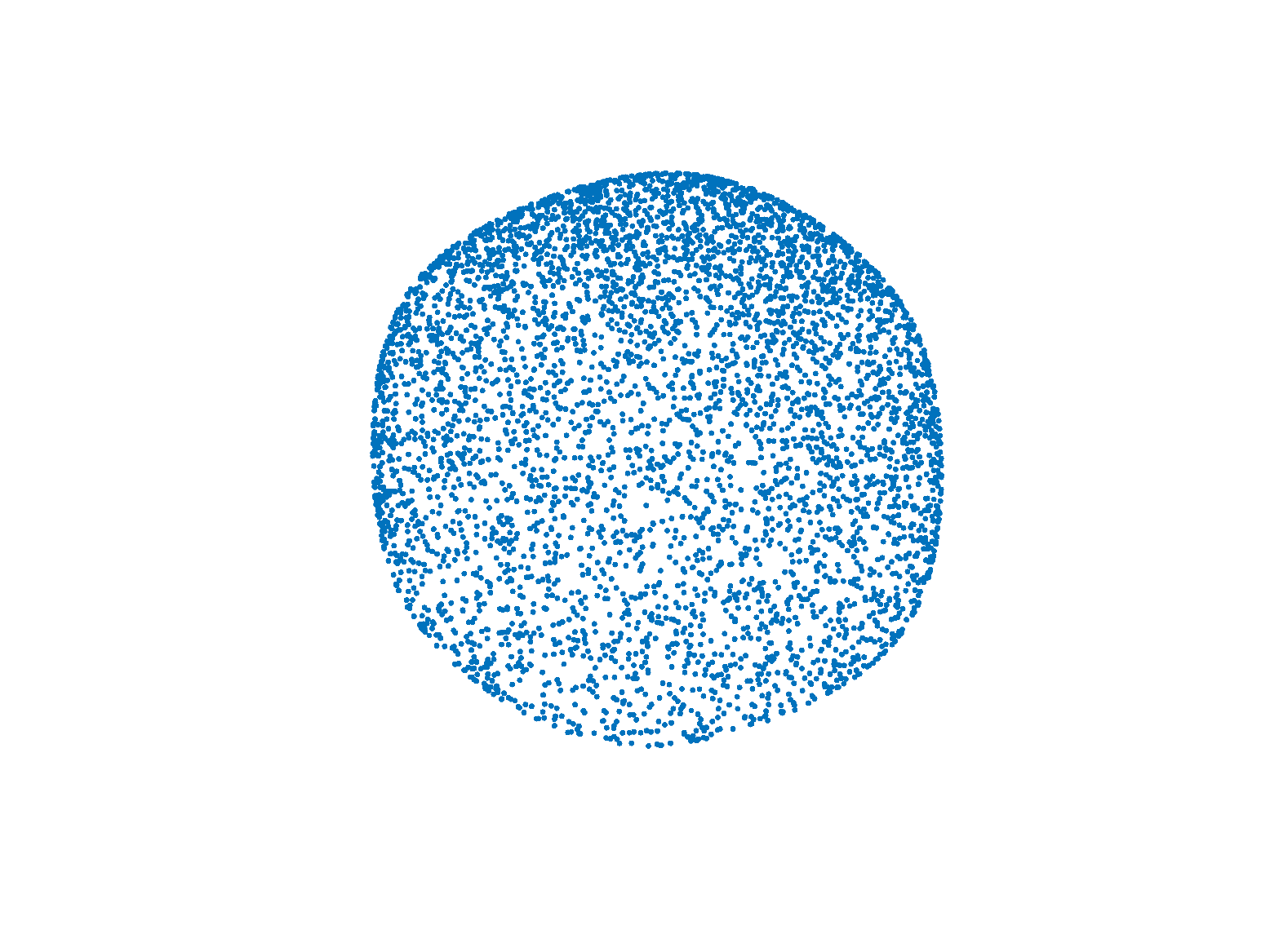}
    \put(38,86){$\bm{s = 0.5}$}
    \end{overpic}
    \vspace{1mm}
    \begin{overpic}[width=0.19\textwidth, trim=60 60 120 40, clip=true,tics=10]{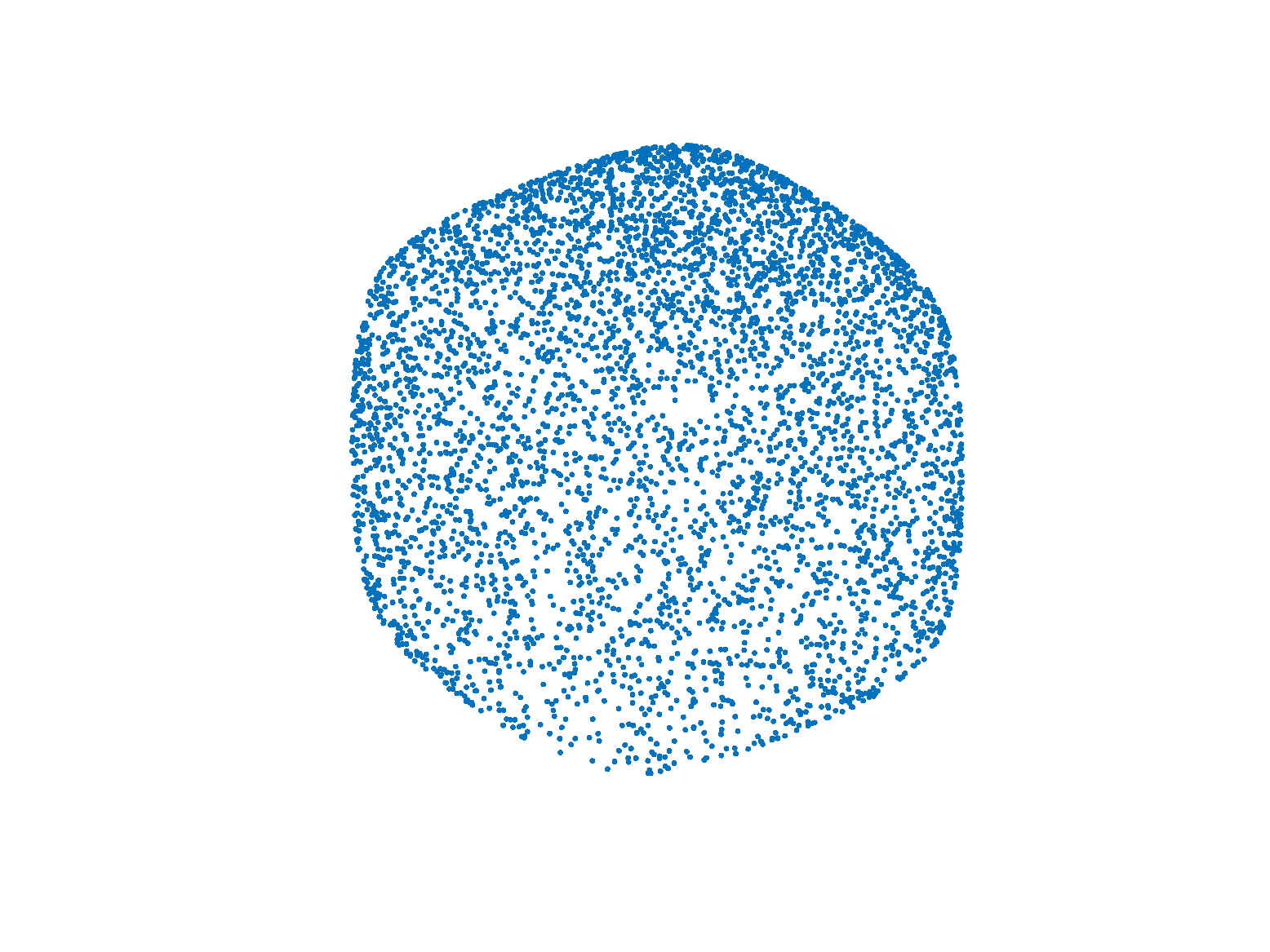}
    \put(38,86){$\bm{s = 0.7}$}
    \end{overpic}
    \vspace{1mm}
    \begin{overpic}[width=0.19\textwidth, trim=60 60 120 40, clip=true,tics=10]{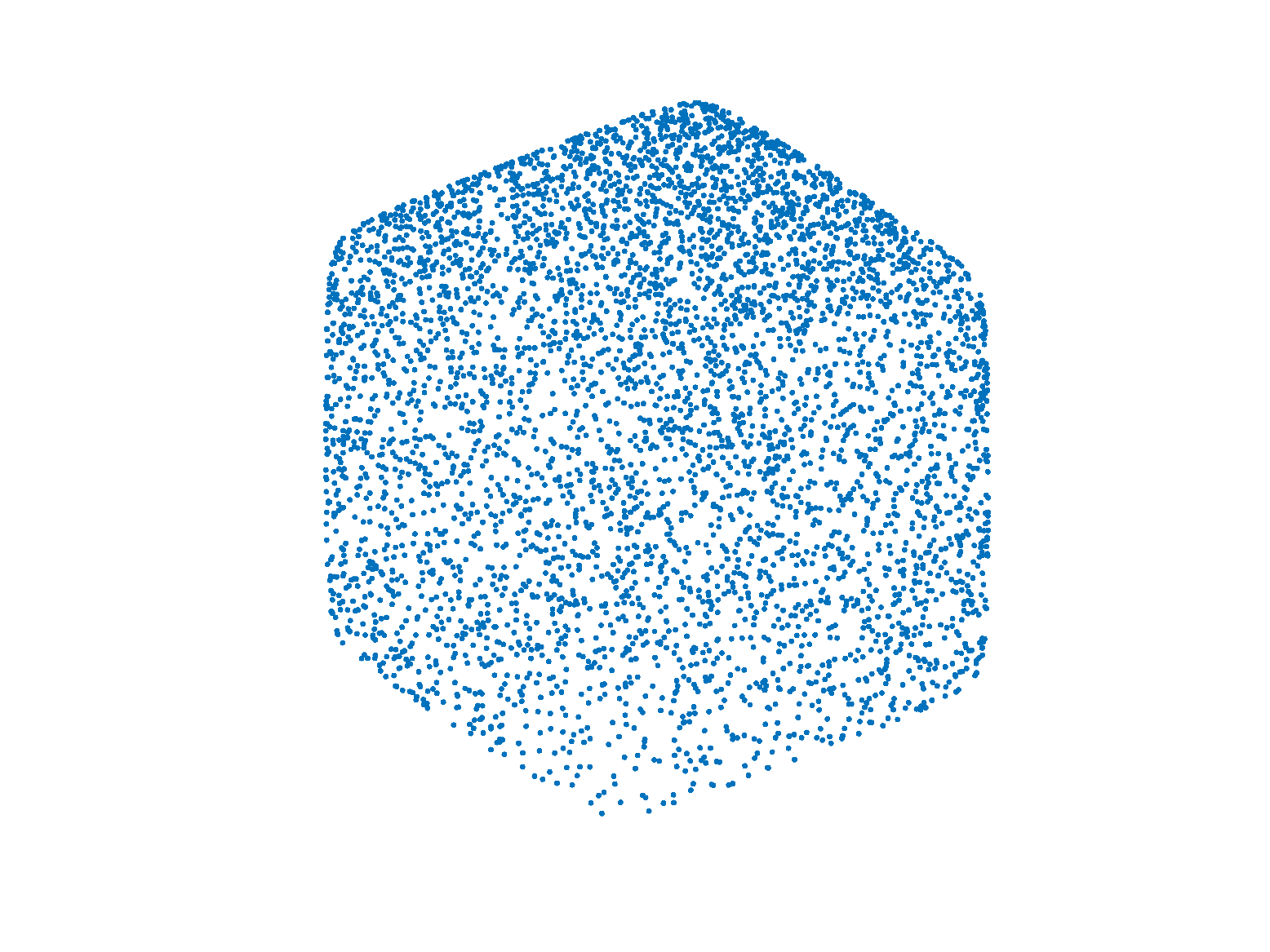}
    \put(38,86){$\bm{s = 0.9}$}
    \end{overpic}
    \\
    \begin{overpic}[width=0.19\textwidth, trim=60 60 120 40, clip=true,tics=10]{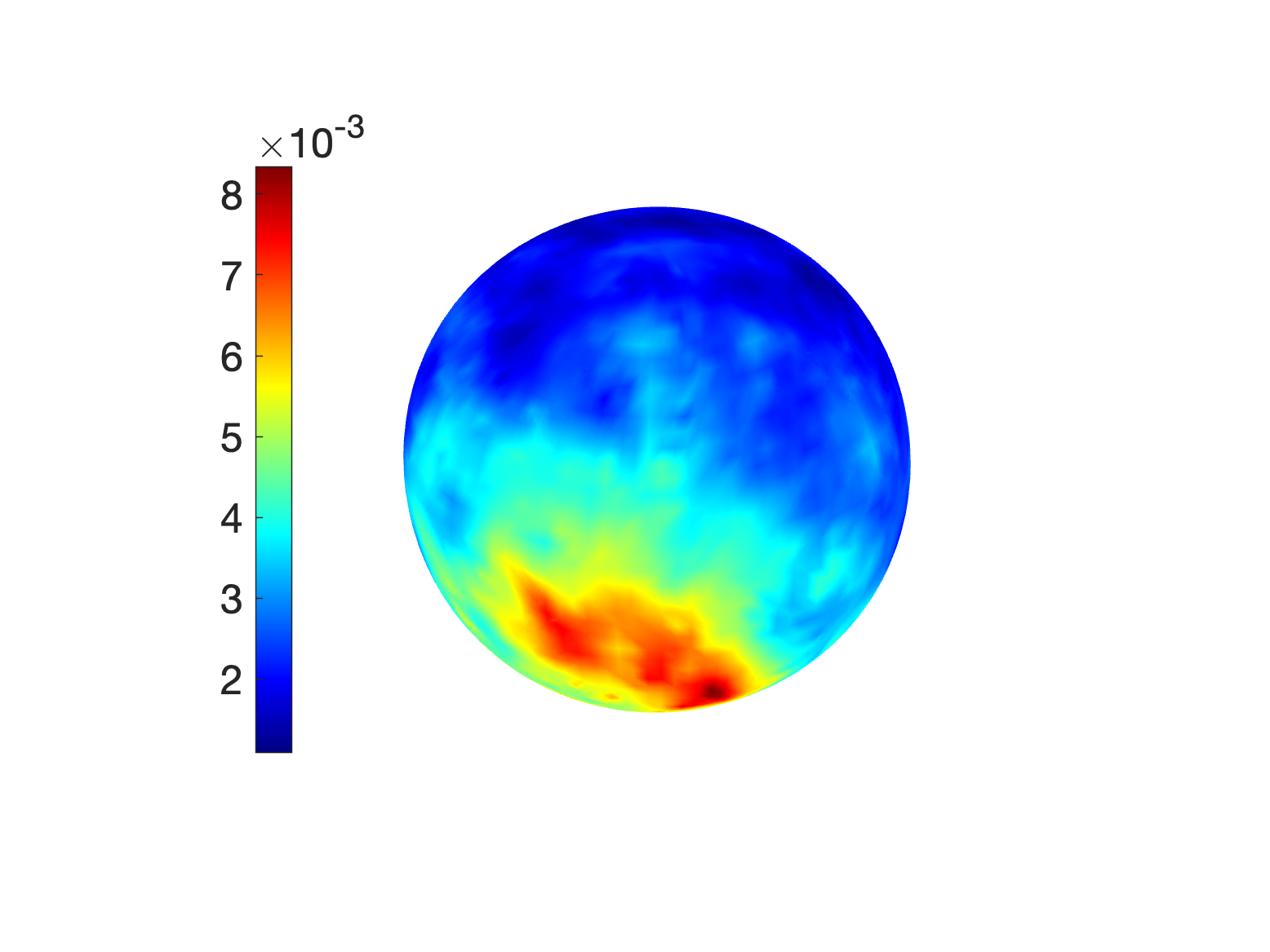}
    \end{overpic}
    \begin{overpic}[width=0.19\textwidth, trim=60 60 120 40, clip=true,tics=10]{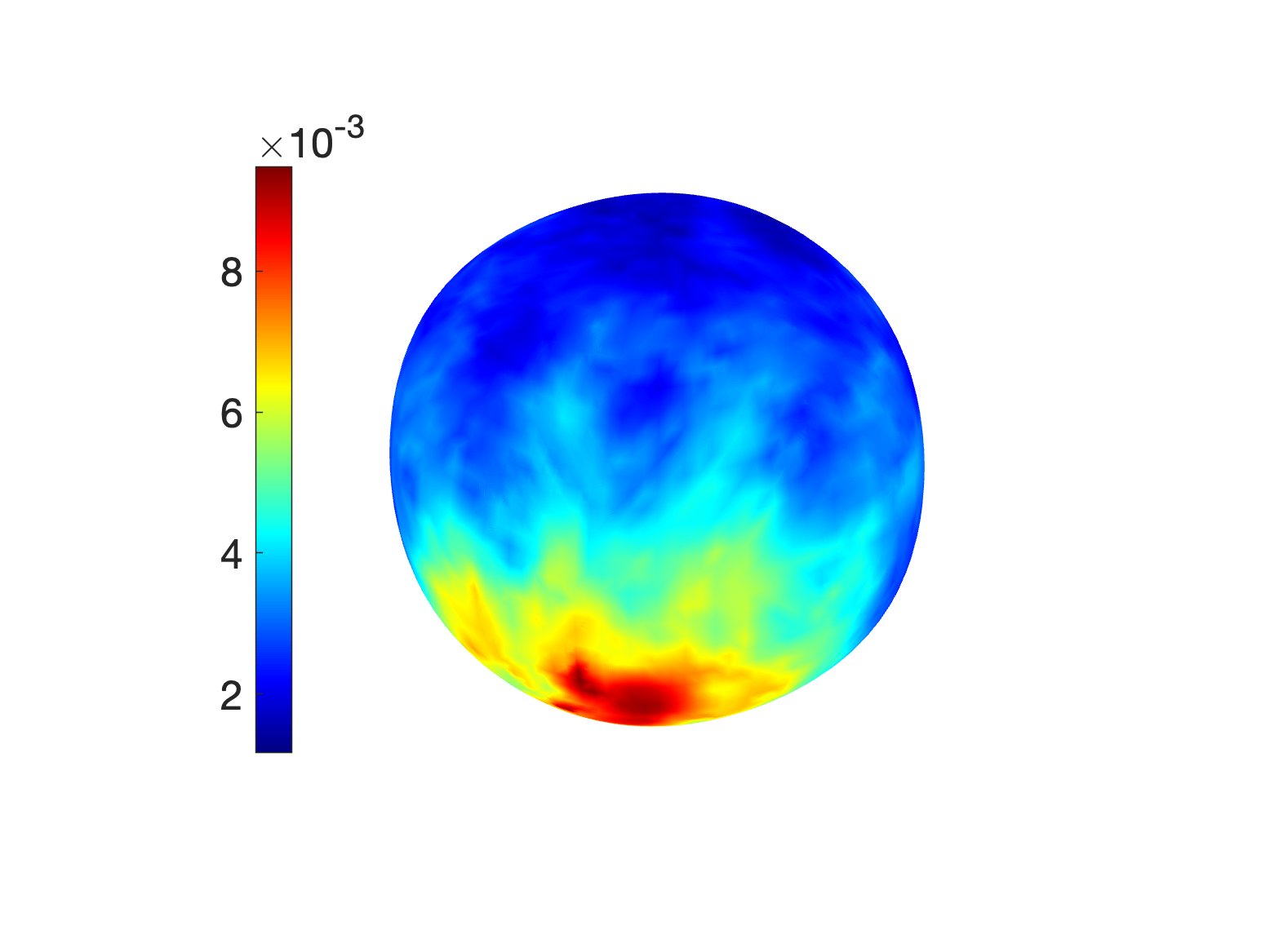}
    \end{overpic}
    \begin{overpic}[width=0.19\textwidth, trim=60 60 120 40, clip=true,tics=10]{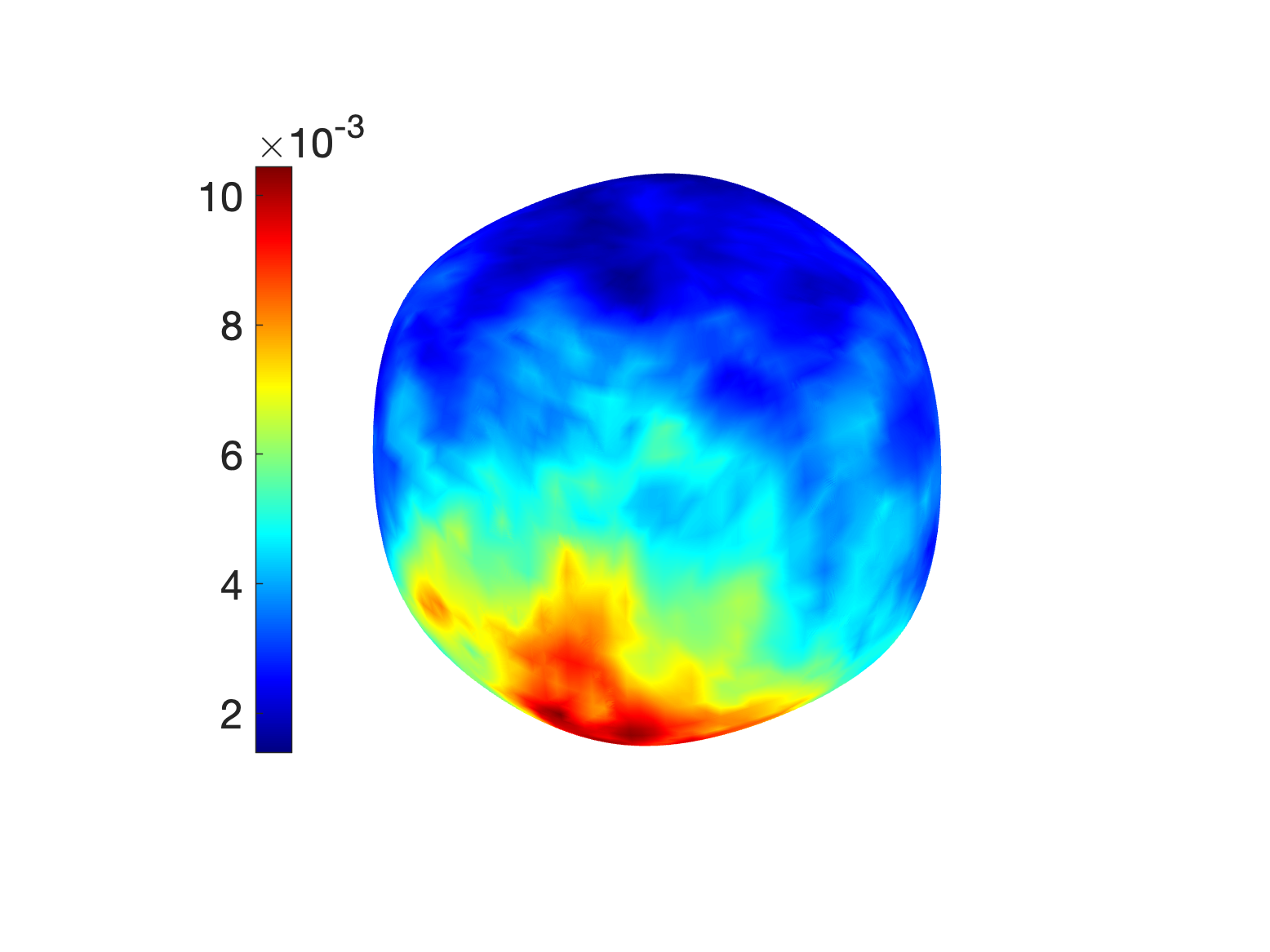}
    \end{overpic}
    \begin{overpic}[width=0.19\textwidth, trim=60 60 120 40, clip=true,tics=10]{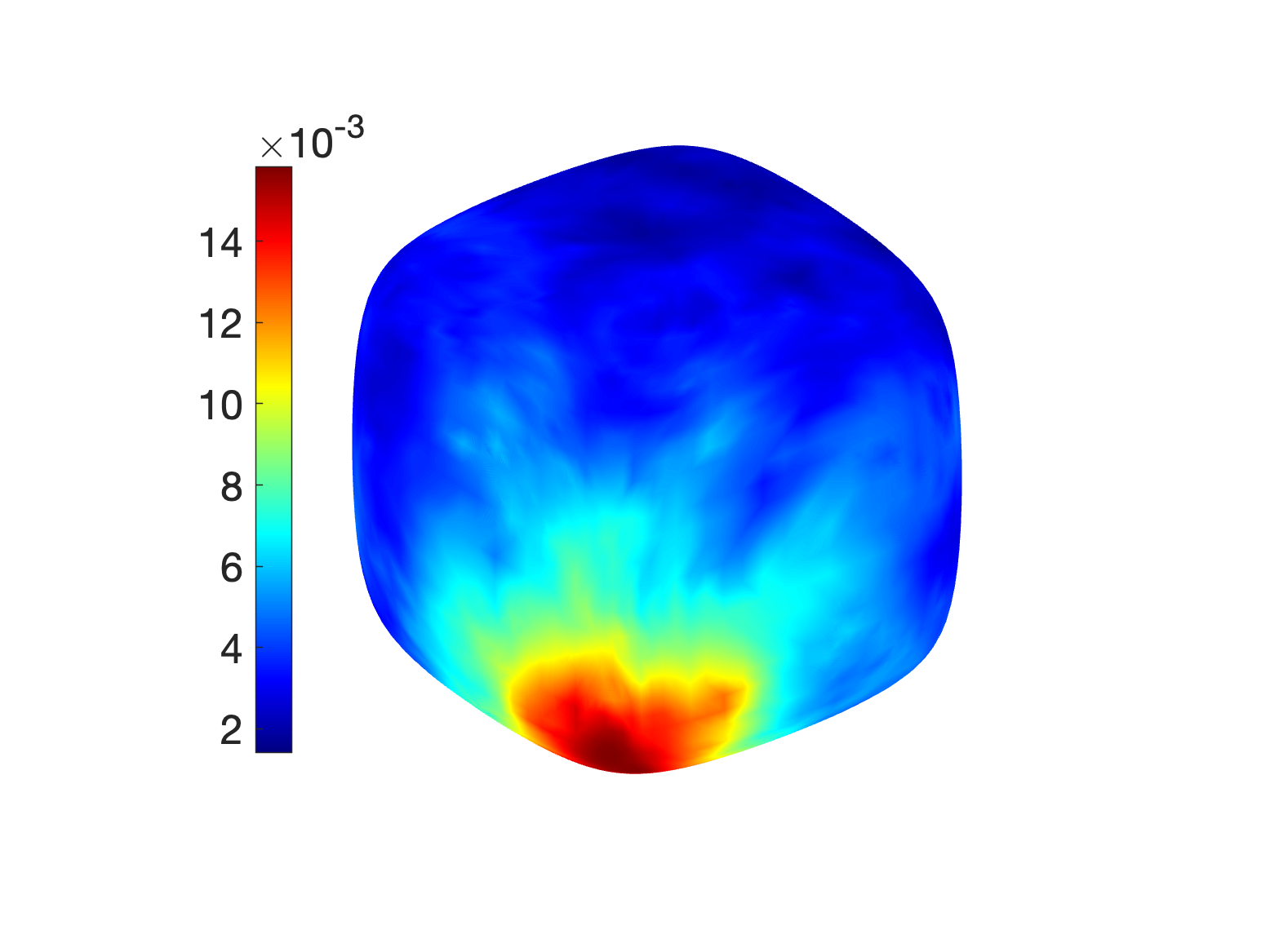}
    \end{overpic}
    \begin{overpic}[width=0.19\textwidth, trim=60 60 120 40, clip=true,tics=10]{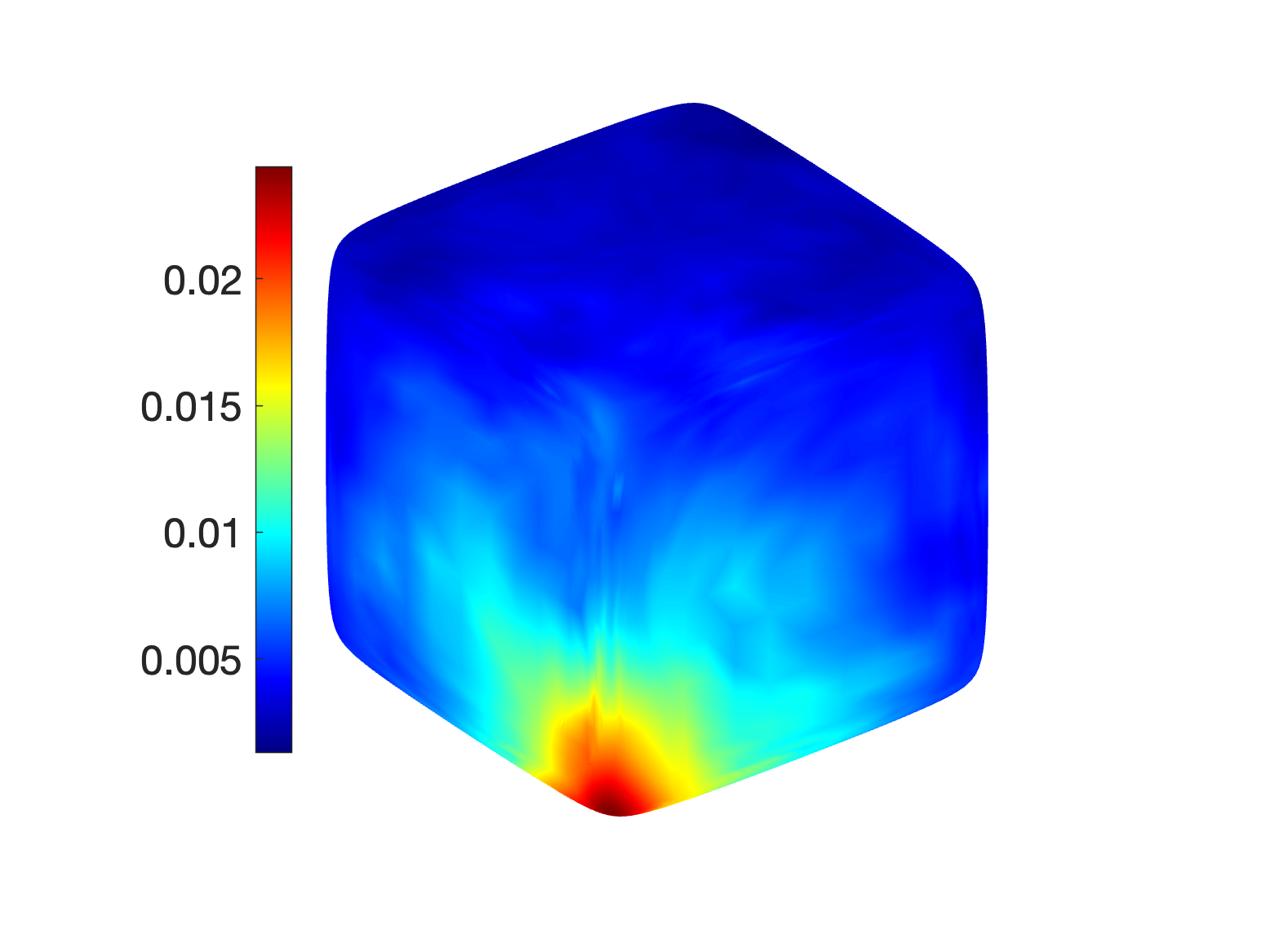}
    \end{overpic}
    \caption{Point distribution of a Flattening sphube of size $N = 5000$ with varying densities. The bottom row shows the local fill distance of the point set where red indicate larger distance (coarse region) and blue indicate smaller one (fine region).}
    \label{fig:flattening_sphube_ptset}
\end{figure}

\begin{figure}
    \centering
    \vspace{-2mm}
    \begin{overpic}[width=0.33\textwidth, trim=140 60 50 40, clip=true,tics=10]{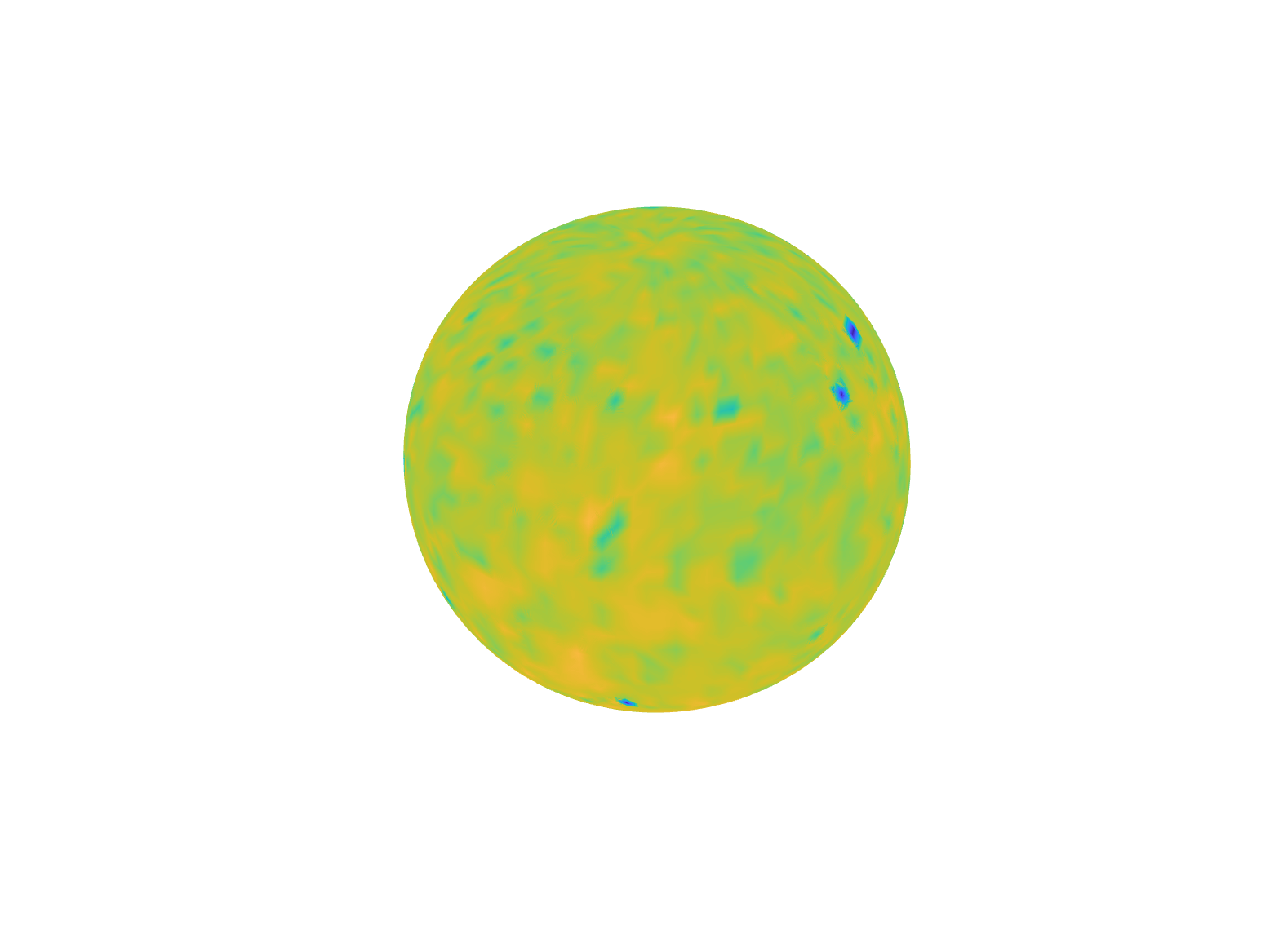}
    \put(-10,29){\rotatebox{90}{$s = 0.1$}}
    \put(34,86){\textbf{PCA}}
    \put(10,75){$\| \mathcal{E}_{\n}\|_{\infty}=7.4\times10^{-2}$}
    \end{overpic}
    \hspace{-7mm}\vspace{-2mm}
    \begin{overpic}[width=0.33\textwidth, trim=140 60 50 40, clip=true,tics=10]{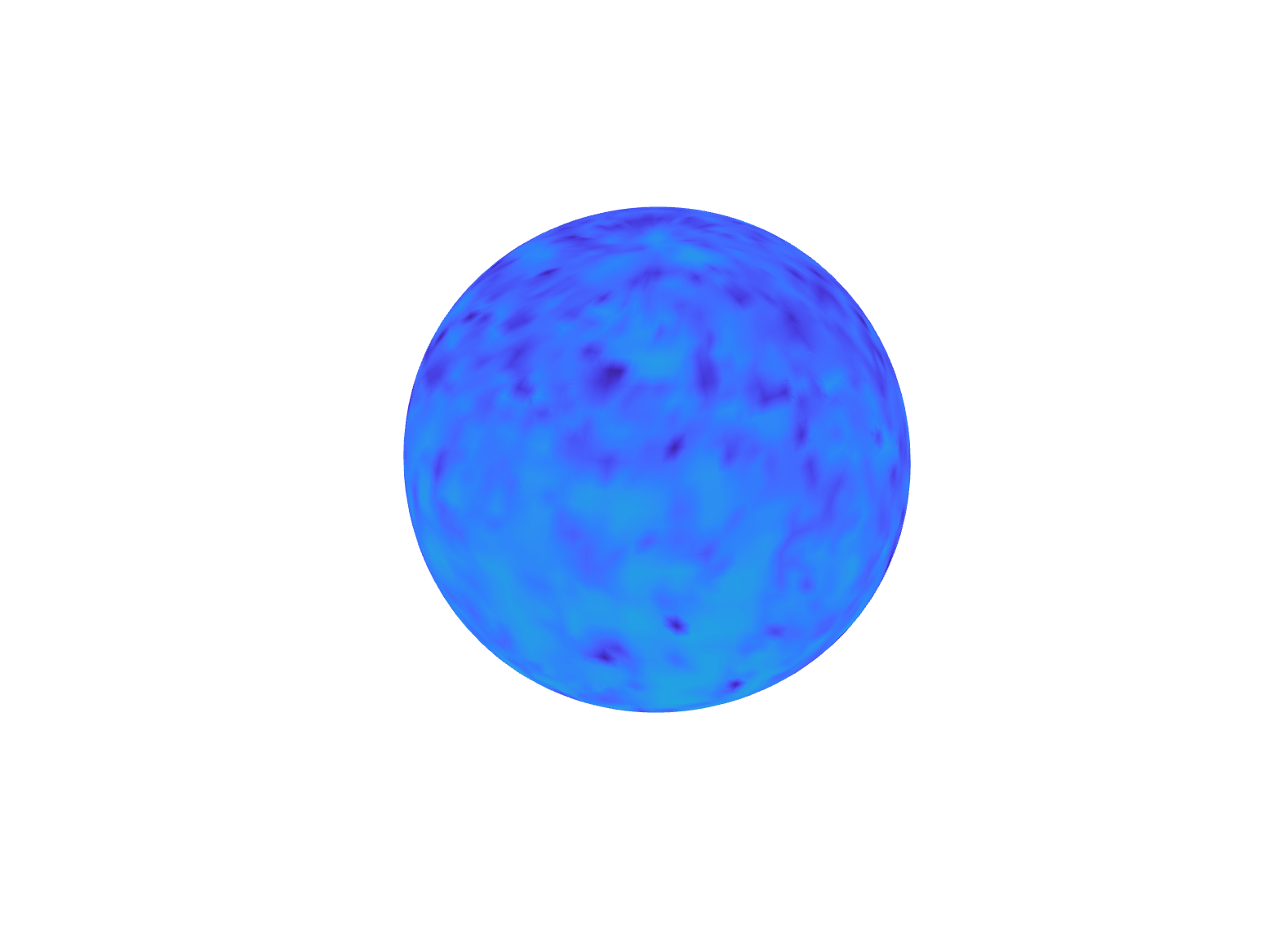}
    \put(8,86){\textbf{Best RBFs/HRBF}}
    \put(11,75){$\| \mathcal{E}_{\n}\|_{\infty}=4.7\times10^{-5}$}
    \end{overpic}
    \hspace{-7mm}\vspace{-2mm}
    \begin{overpic}[width=0.33\textwidth, trim=140 60 50 40, clip=true,tics=10]{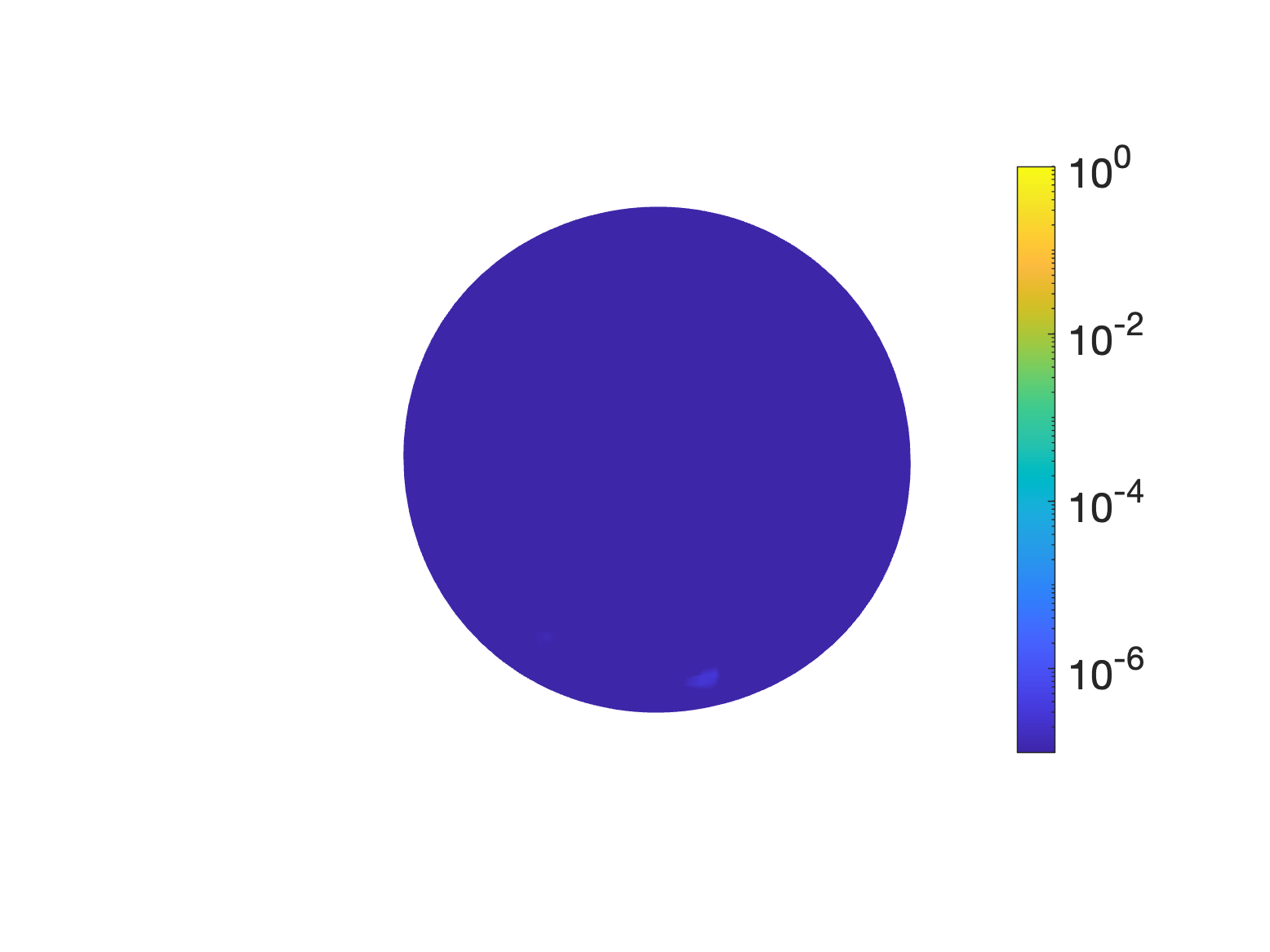}
    \put(32,86){\textbf{KRBF}}
    \put(12,75){$\| \mathcal{E}_{\n}\|_{\infty}=3.4\times10^{-7}$}
    \end{overpic}
    \\
    \vspace{-2mm}
    \begin{overpic}[width=0.33\textwidth, trim=140 60 50 40, clip=true,tics=10]{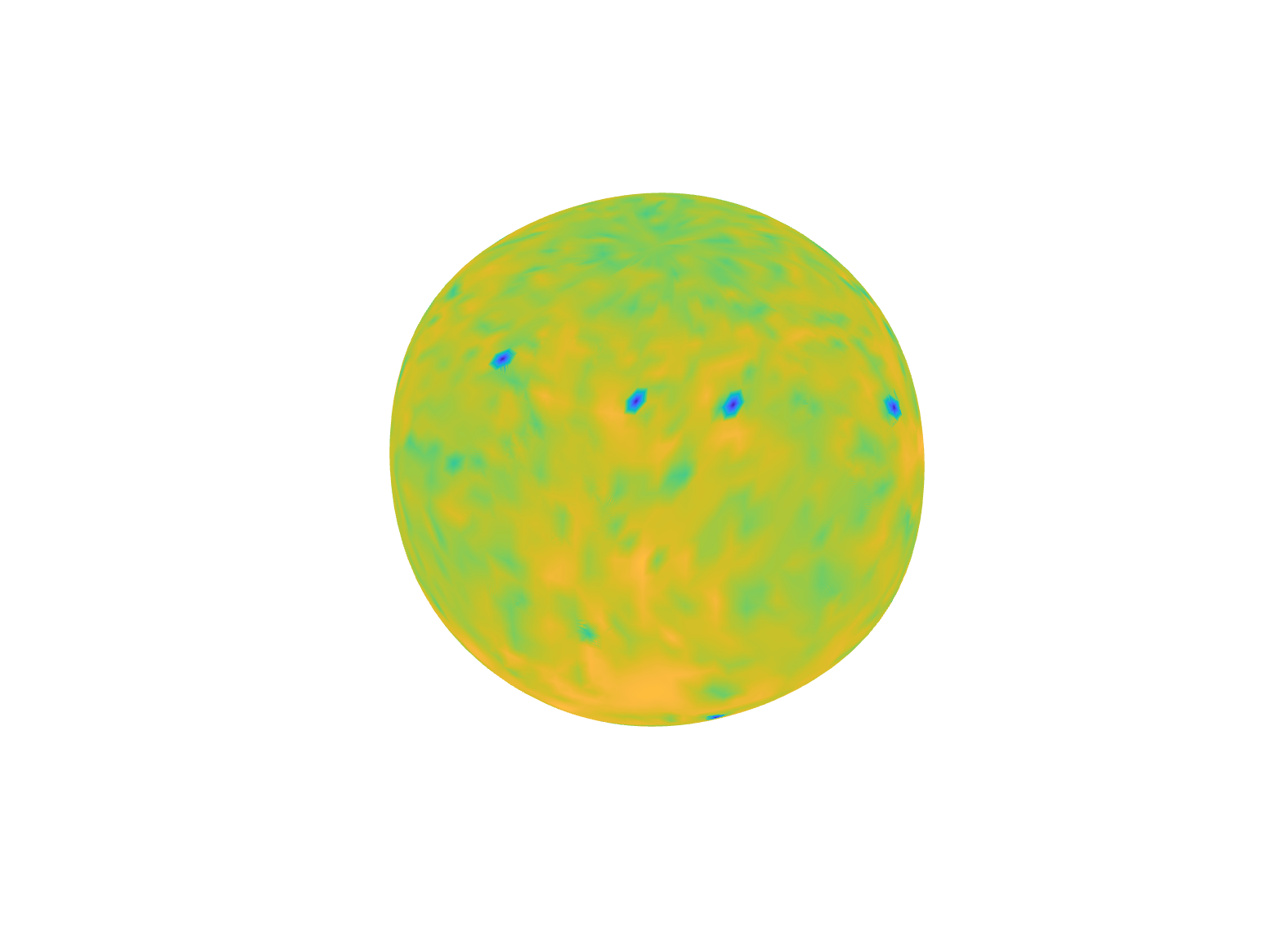}
    \put(-10,29){\rotatebox{90}{$s = 0.3$}}
    \put(10,77){$\| \mathcal{E}_{\n}\|_{\infty}=1.0\times10^{-1}$}
    \end{overpic}
    \hspace{-7mm}\vspace{-2mm}
    \begin{overpic}[width=0.33\textwidth, trim=140 60 50 40, clip=true,tics=10]{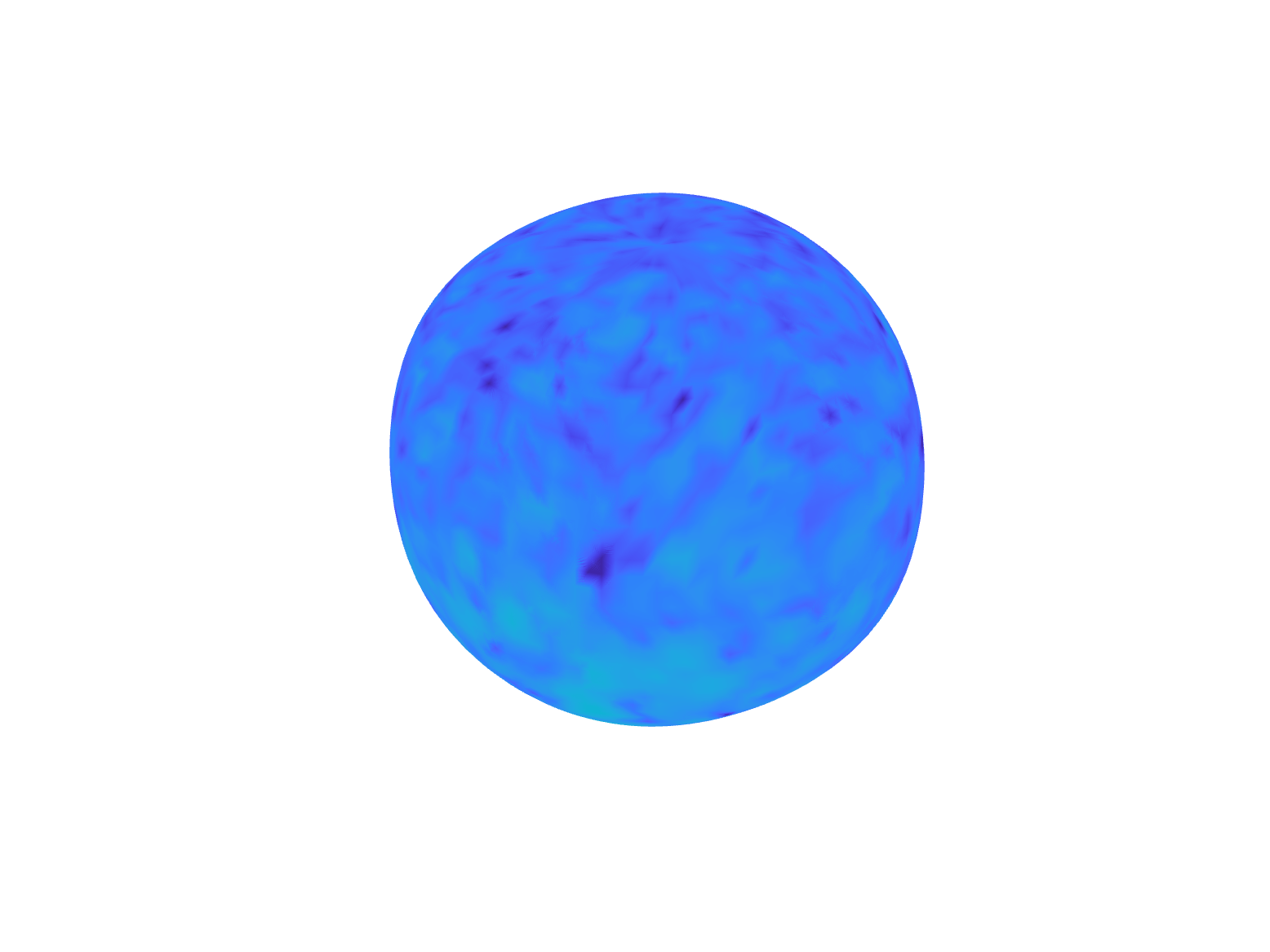}
    \put(11,77){$\| \mathcal{E}_{\n}\|_{\infty}=1.2\times10^{-4}$}
    \end{overpic}
    \hspace{-7mm}\vspace{-2mm}
    \begin{overpic}[width=0.33\textwidth, trim=140 60 50 40, clip=true,tics=10]{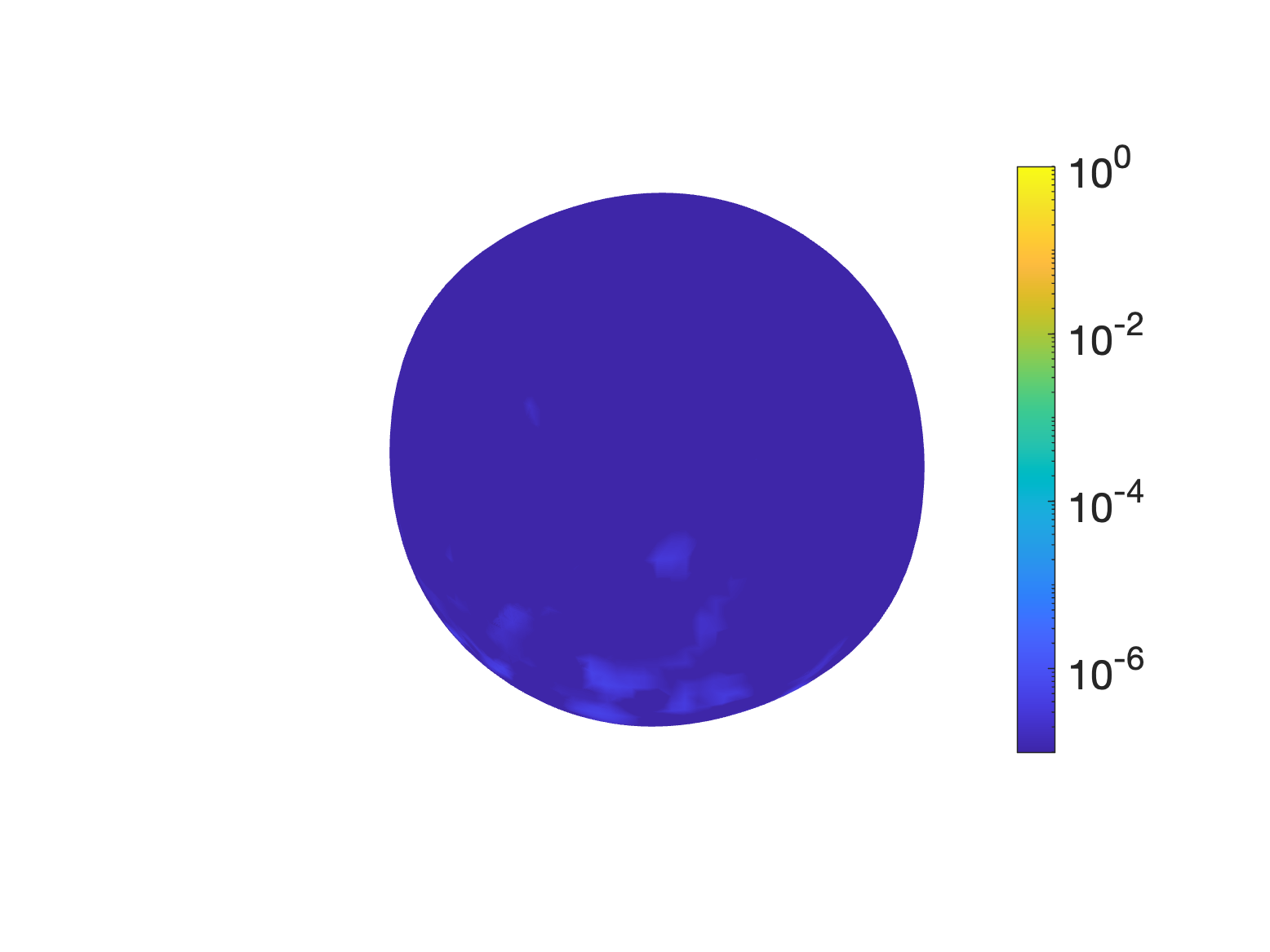}
    \put(12,77){$\| \mathcal{E}_{\n}\|_{\infty}=5.4\times10^{-7}$}
    \end{overpic}
    \\
    \begin{overpic}[width=0.33\textwidth, trim=140 60 50 40, clip=true,tics=10]{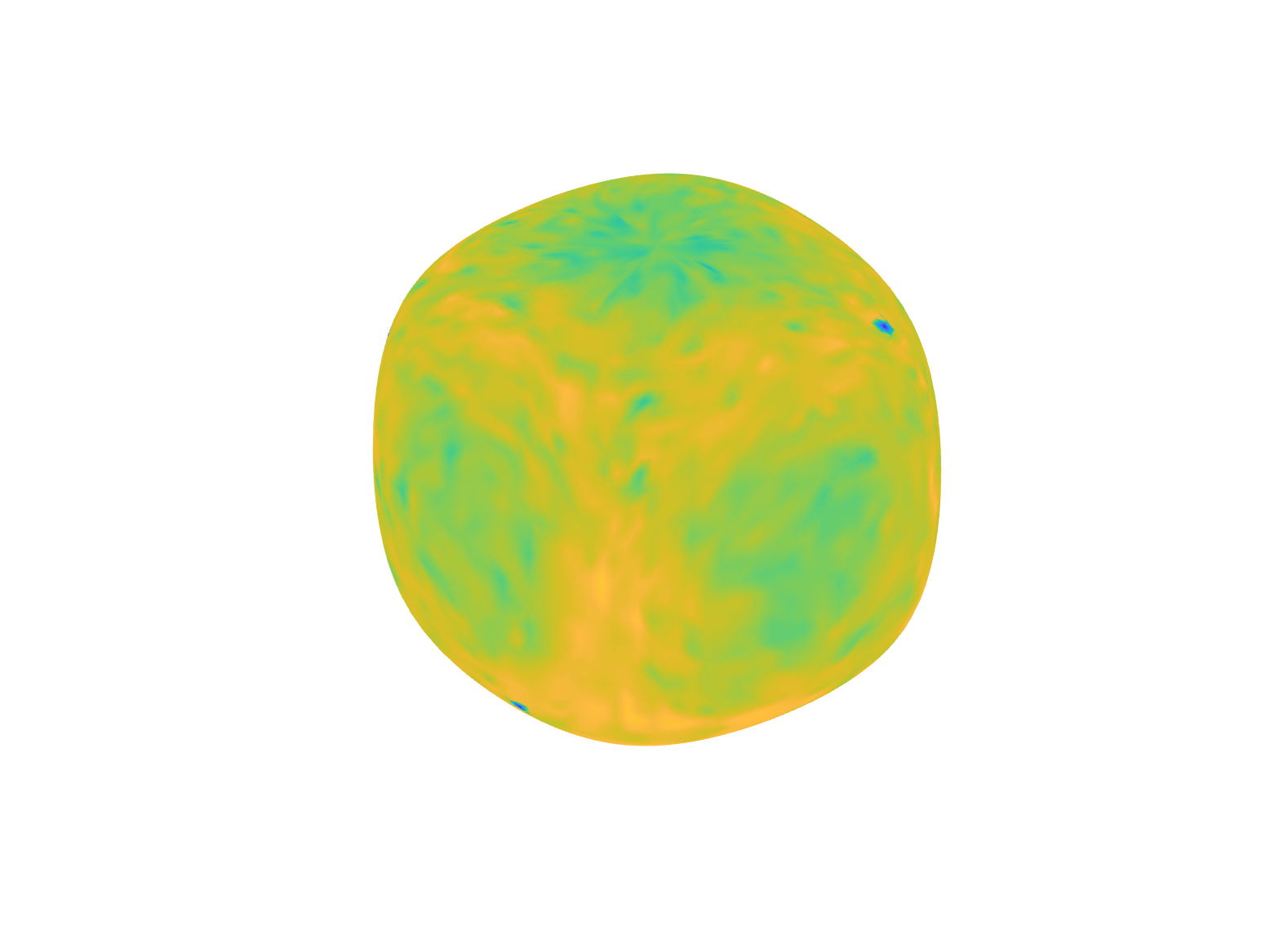}
    \put(-10,29){\rotatebox{90}{$s = 0.5$}}
    \put(10,81){$\| \mathcal{E}_{\n}\|_{\infty}=1.3\times10^{-1}$}
    \end{overpic}
    \hspace{-7mm}
    \begin{overpic}[width=0.33\textwidth, trim=140 60 50 40, clip=true,tics=10]{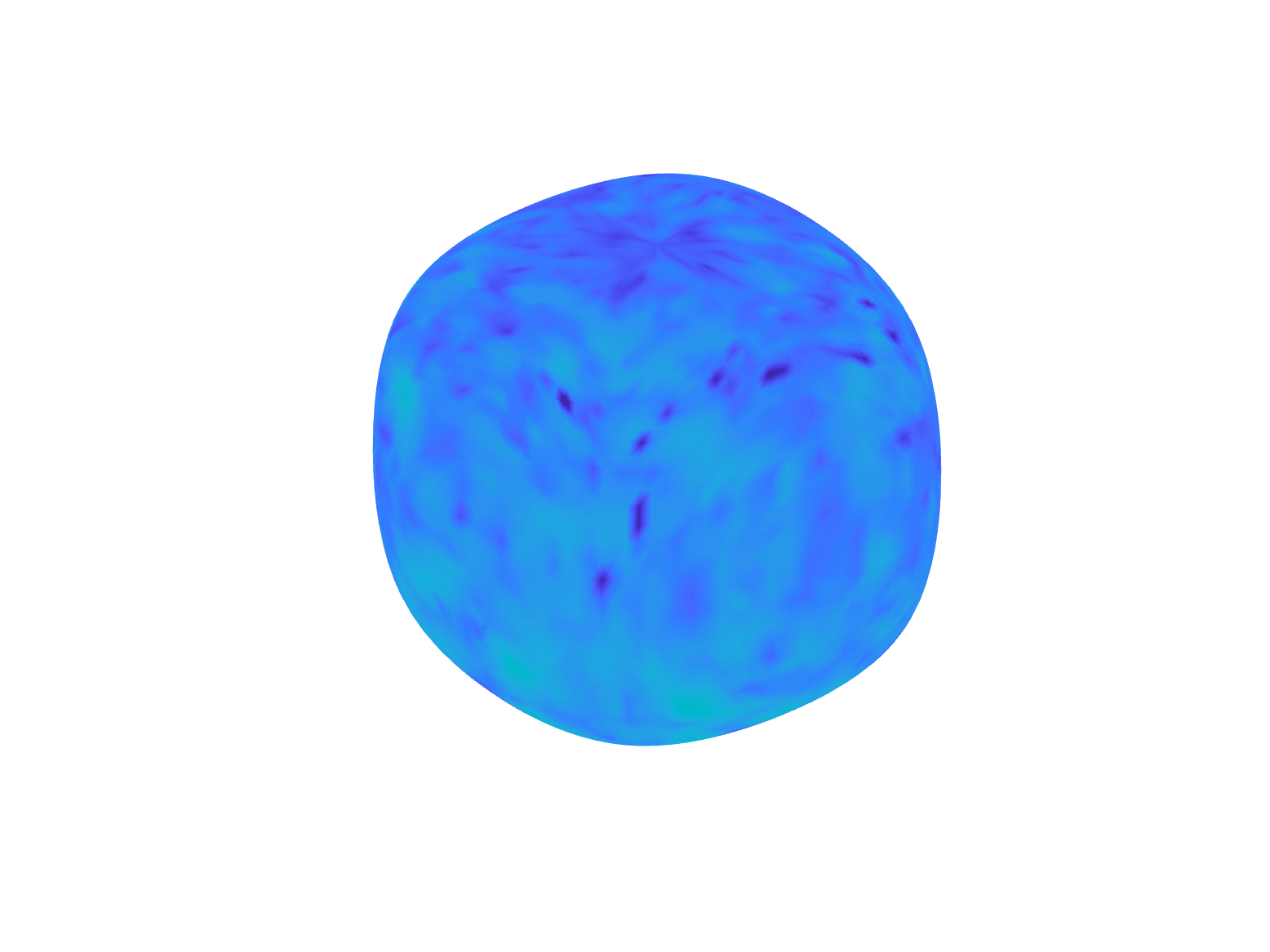}
    \put(11,81){$\| \mathcal{E}_{\n}\|_{\infty}=1.9\times10^{-4}$}
    \end{overpic}
    \hspace{-7mm}
    \begin{overpic}[width=0.33\textwidth, trim=140 60 50 40, clip=true,tics=10]{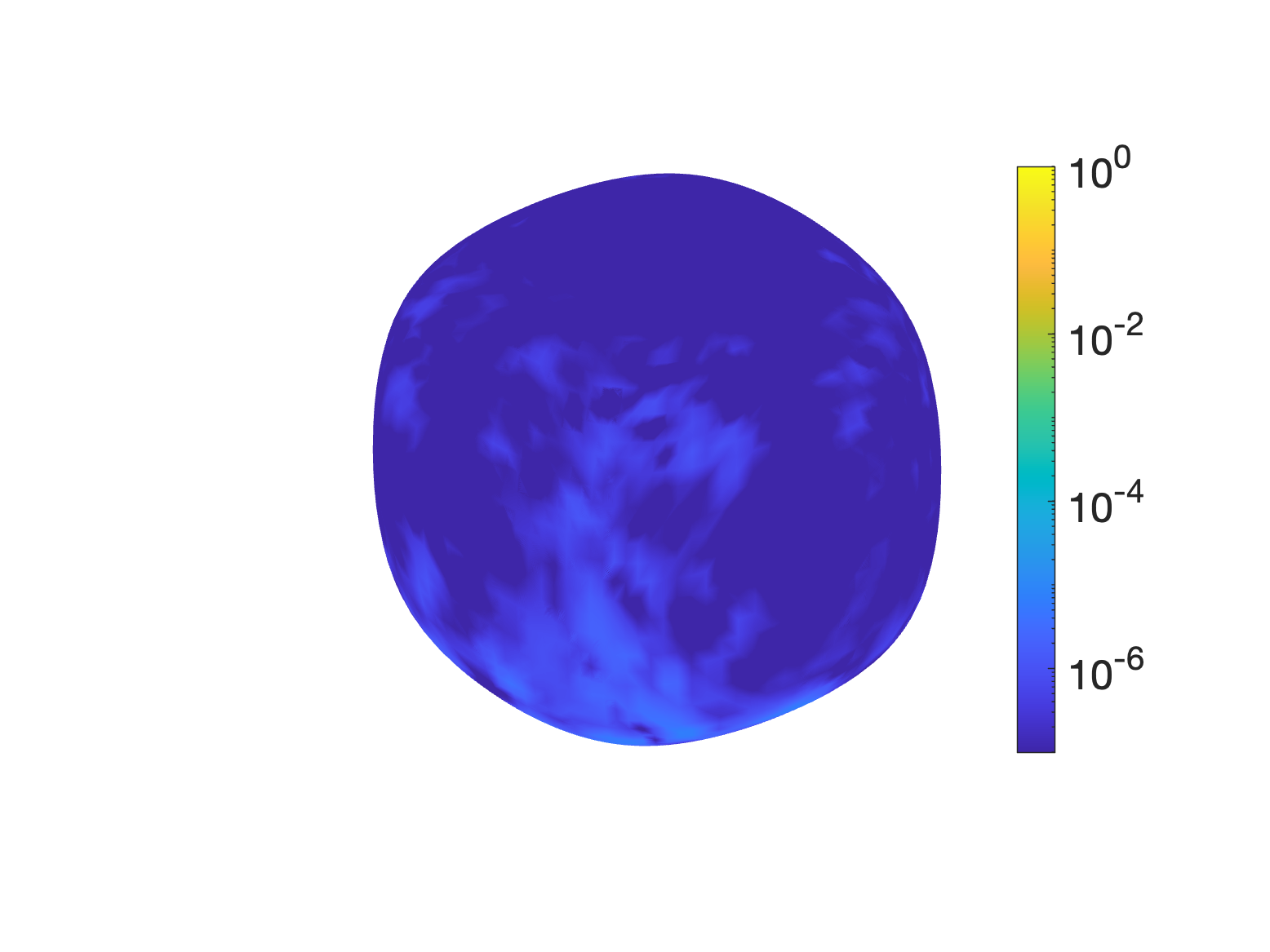}
    \put(12,81){$\| \mathcal{E}_{\n}\|_{\infty}=7.9\times10^{-6}$}
    \end{overpic}
    \\
    \begin{overpic}[width=0.33\textwidth, trim=140 60 50 40, clip=true,tics=10]{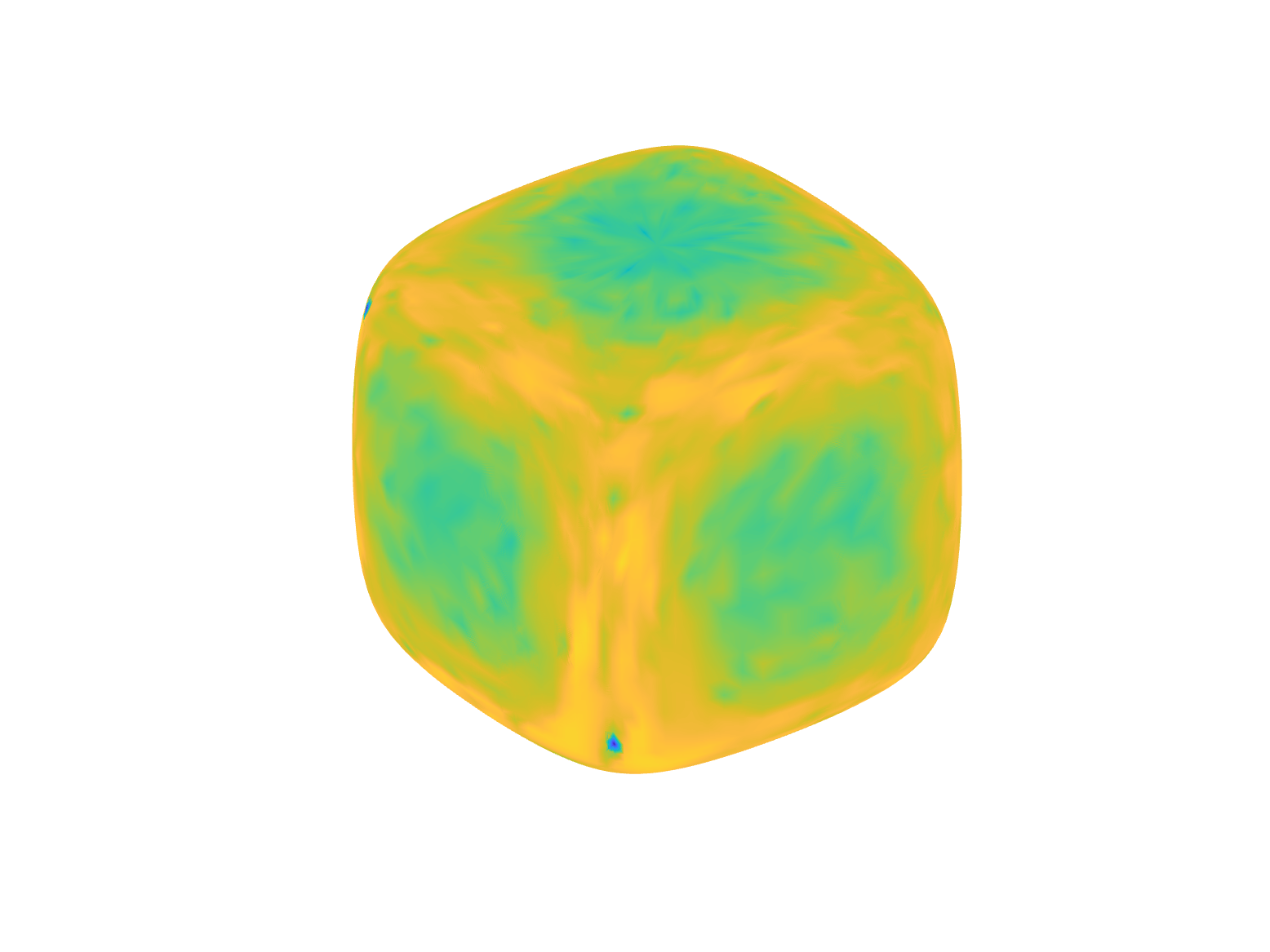}
    \put(-10,29){\rotatebox{90}{$s = 0.7$}}
    \put(10,87){$\| \mathcal{E}_{\n}\|_{\infty}=2.2\times10^{-1}$}
    \end{overpic}
    \hspace{-7mm}
    \begin{overpic}[width=0.33\textwidth, trim=140 60 50 40, clip=true,tics=10]{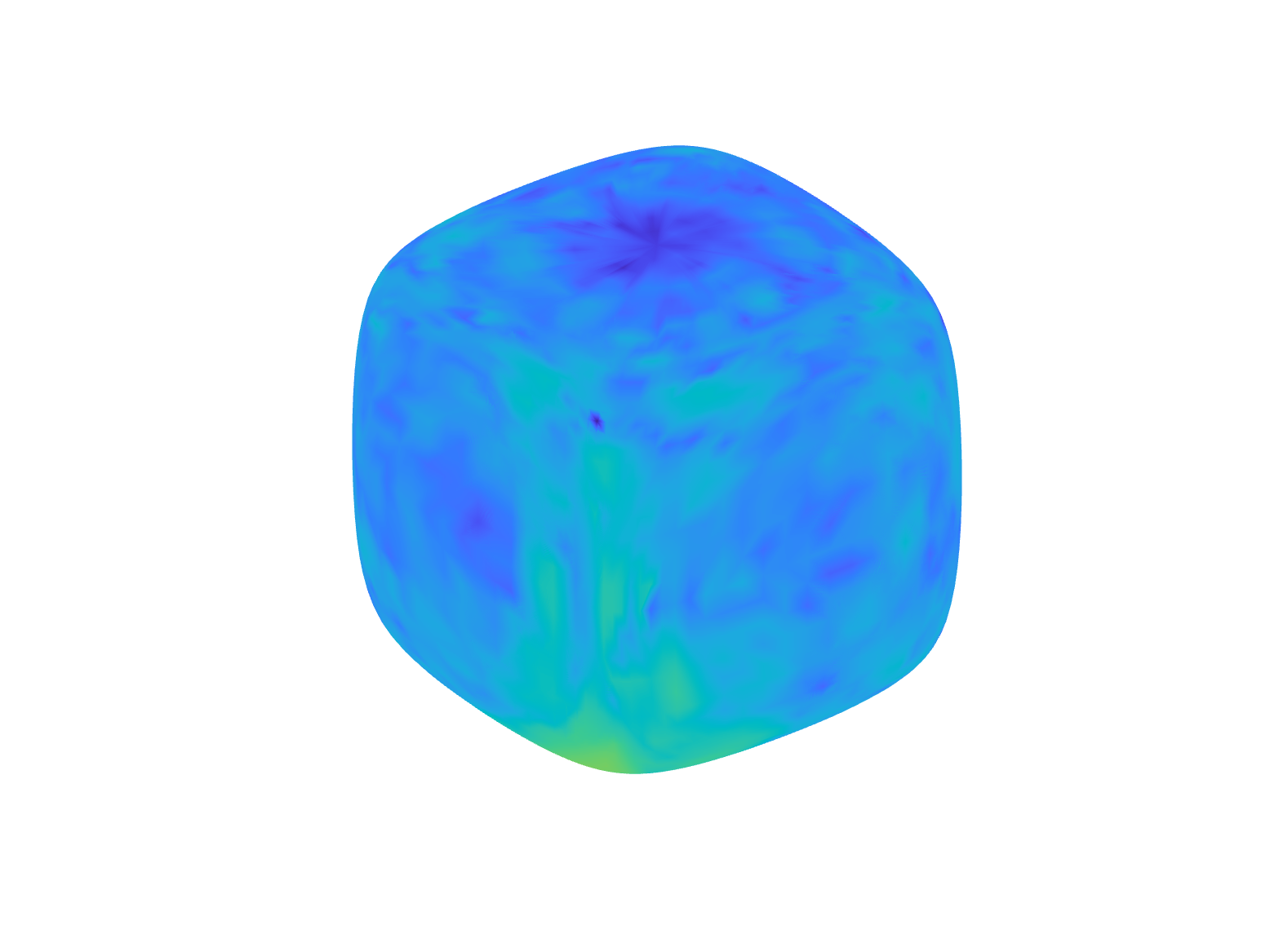}
    \put(11,87){$\| \mathcal{E}_{\n}\|_{\infty}=3.6\times10^{-3}$}
    \end{overpic}
    \hspace{-7mm}
    \begin{overpic}[width=0.33\textwidth, trim=140 60 50 40, clip=true,tics=10]{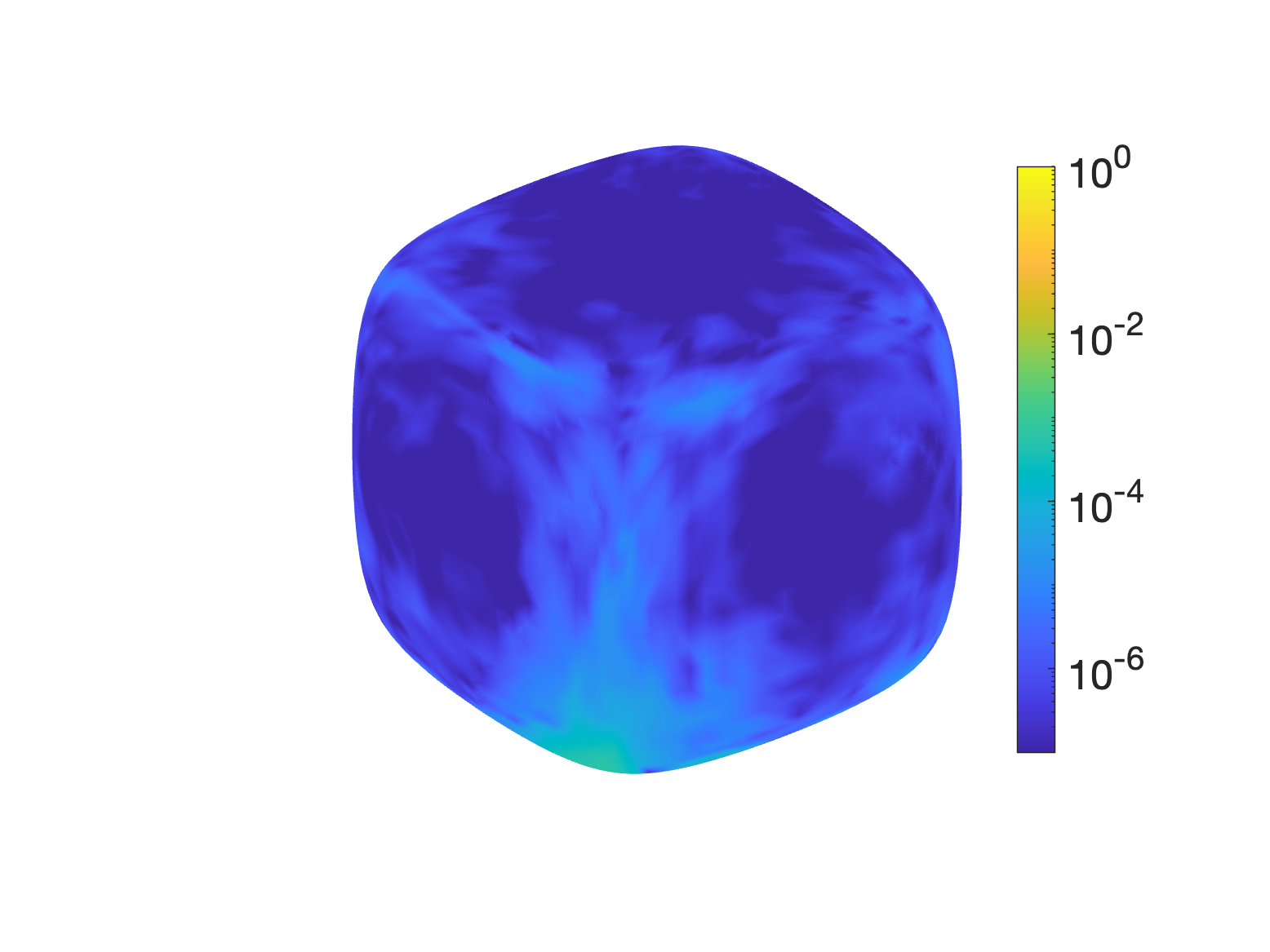}
    \put(12,87){$\| \mathcal{E}_{\n}\|_{\infty}=5.6\times10^{-4}$}
    \end{overpic}
    \\
    \vspace{4mm}
    \begin{overpic}[width=0.33\textwidth, trim=140 60 50 40, clip=true,tics=10]{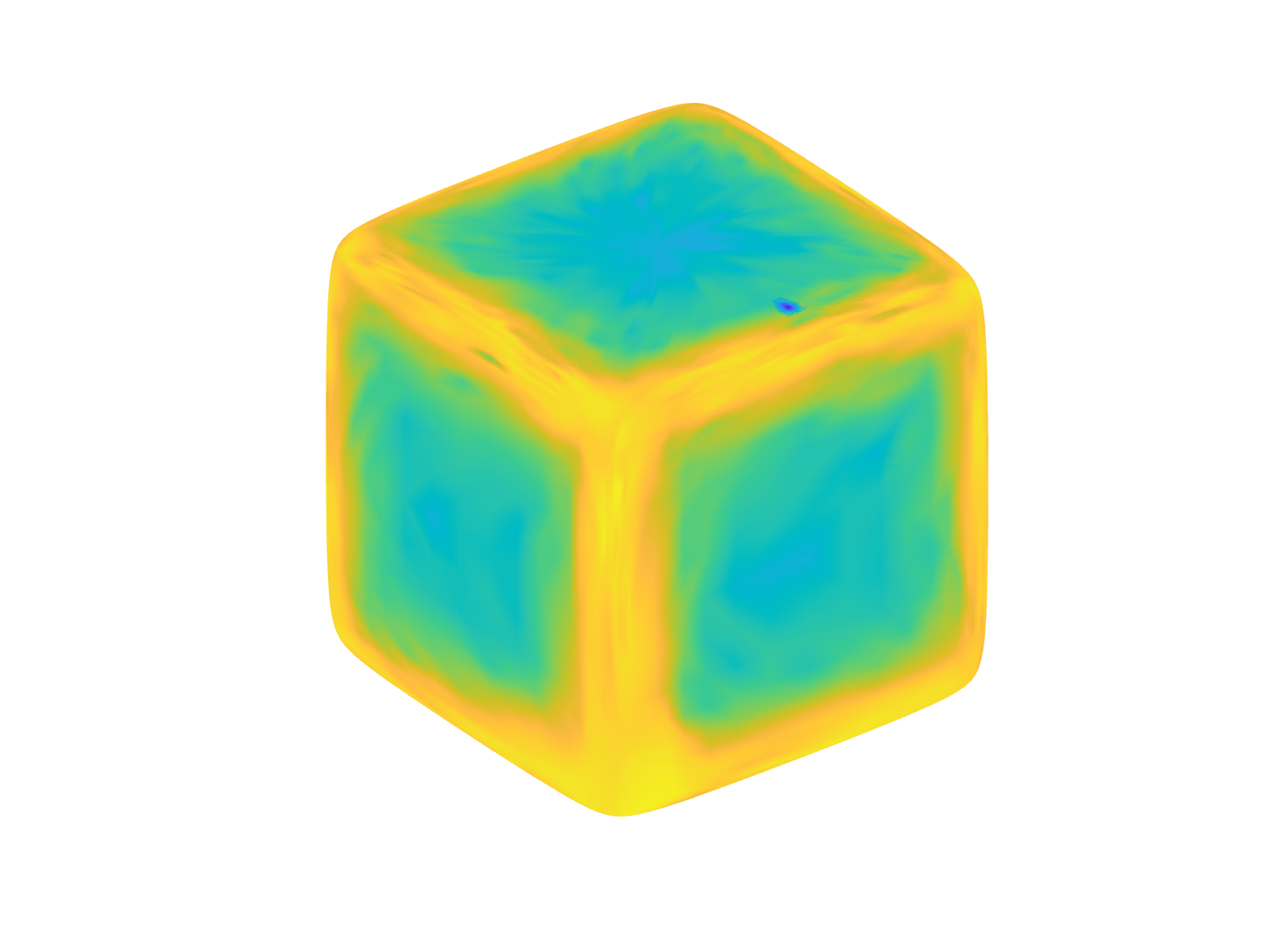}
    \put(-10,29){\rotatebox{90}{$s = 0.9$}}
    \put(10,89){$\| \mathcal{E}_{\n}\|_{\infty}=6.5\times10^{-1}$}
    \end{overpic}
    \hspace{-7mm} \vspace{4mm}
    \begin{overpic}[width=0.33\textwidth, trim=140 60 50 40, clip=true,tics=10]{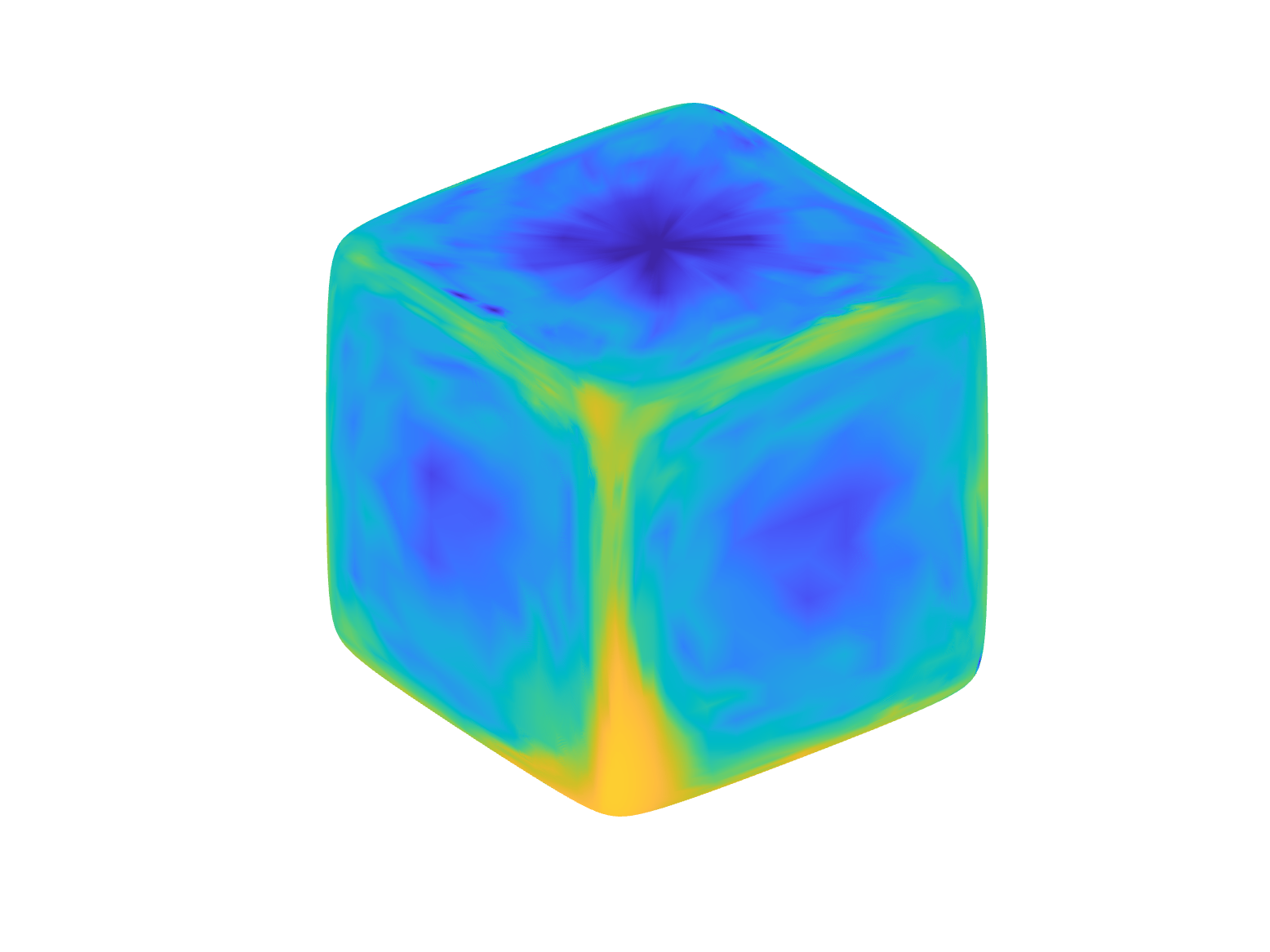}
    \put(11,89){$\| \mathcal{E}_{\n}\|_{\infty}=1.3\times10^{-1}$}
    \end{overpic}
    \hspace{-7mm} \vspace{4mm}
    \begin{overpic}[width=0.33\textwidth, trim=140 60 50 40, clip=true,tics=10]{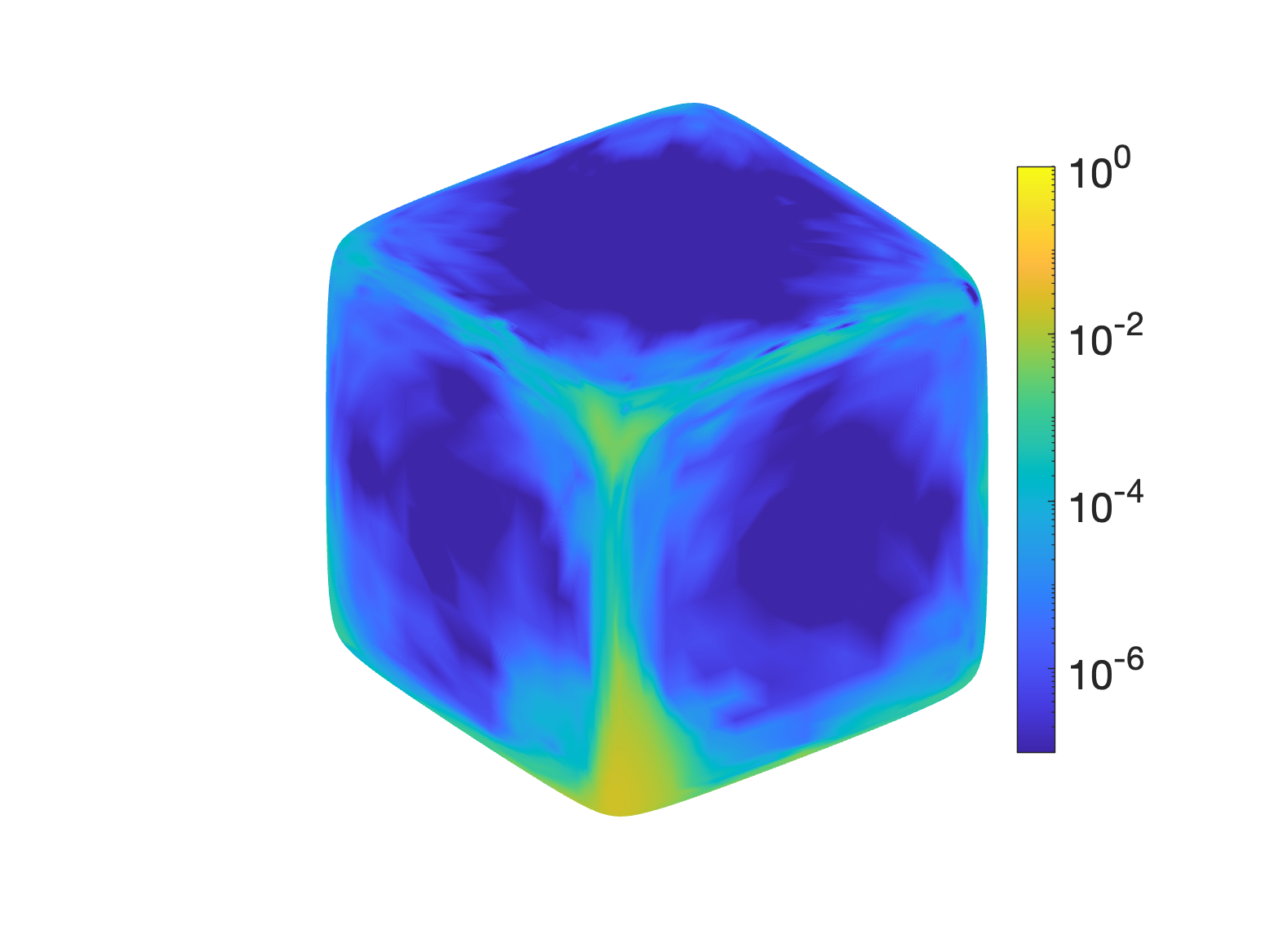}
    \put(12,89){$\| \mathcal{E}_{\n}\|_{\infty}=1.3\times10^{-1}$}
    \end{overpic}
    \vspace{-5mm}
    \caption{Error distribution of estimating surface normal of a Flattening sphube.}
    \label{fig:flattening_sphube}
\end{figure}
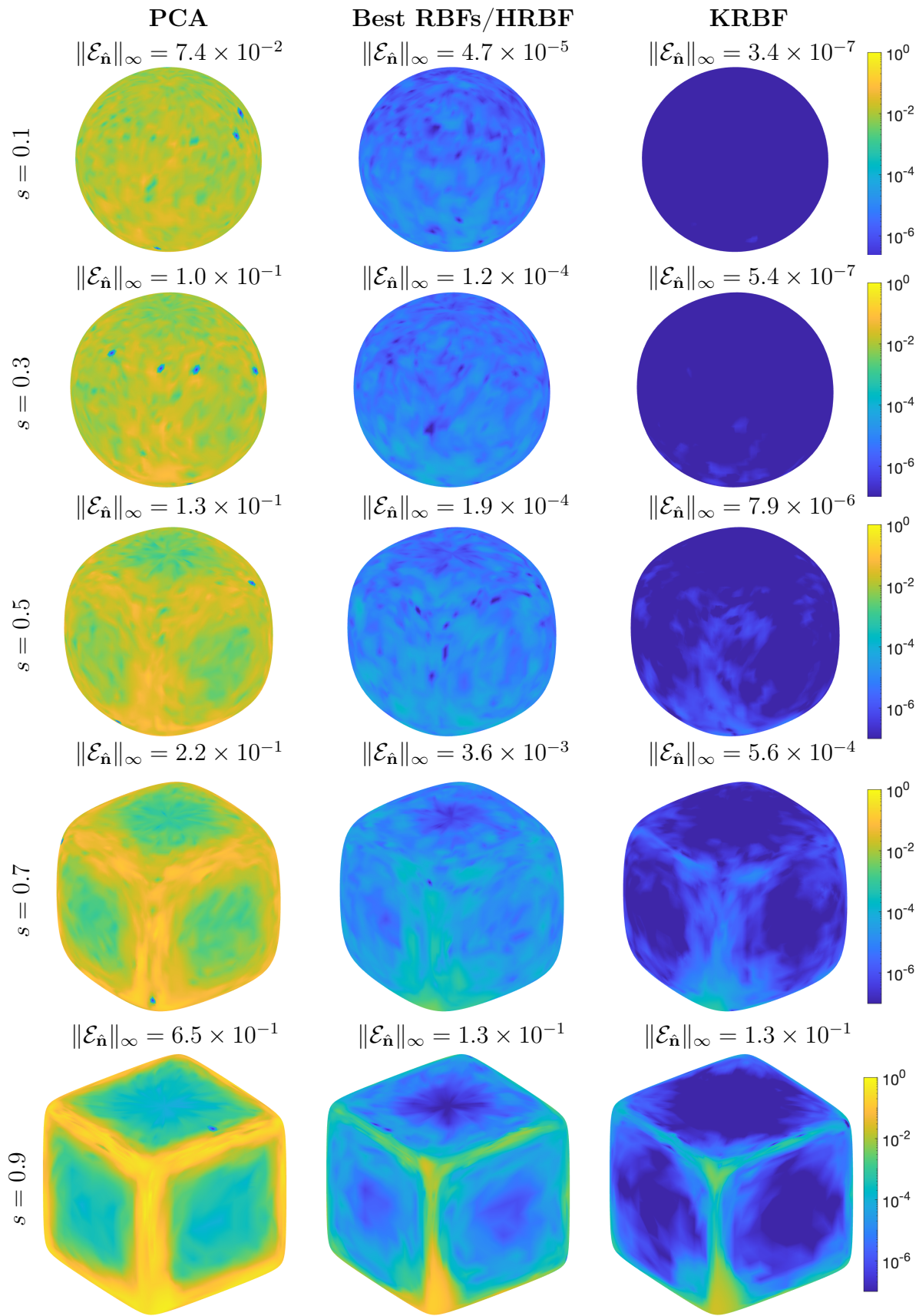

\autoref{fig:flattening_sphube_ptset} shows the distribution of \( N = 5000 \) data points on each surface \(\mathbb{S}_{s,1}\) with densities varying from coarse to fine. The colormap in this figure represents the approximated fill distance of each local stencil, with red indicating coarse regions and blue indicating fine regions. The normal estimation error computed by different methods on each surface is shown in \autoref{fig:flattening_sphube}. In this figure, the first column shows the normal estimation obtained via Principal Component Analysis (PCA), and the second column presents the output for the following approach:
\begin{itemize}
    \item \textbf{Best RBFs/HRBF:} Select the best normal approximation from the \textbf{RBF} and \textbf{HRBF} interpolants as defined in Example 3.
\end{itemize}
The final column is dedicated to KRBF. The colormap in this figure represents the normal estimation error defined in \eqref{eq:error_normal}, with blue indicating a more accurate approximation and yellow indicating a less accurate one.

Note that we are showing the face of $\mathbb{S}_{s,r}$ where the point with the maximum error $\| \mathcal{E}_{\n}\|_{\infty}$ is clearly visible at the bottom of each figure. Above each shape, we display the corresponding value of $\| \mathcal{E}_{\n}\|_{\infty}$, revealing that KRBF is approximately an order of magnitude more accurate than Best RBFs/HRBF in every scenario except the last row ($s=0.9$), where large grid spacing at the red corner induces significant error. As $s$ increases, $\| \mathcal{E}_{\n}\|_{\infty}$ grows where smaller errors are found clustering in the flat region and larger errors are predominantly found around the edges and corners of the near-cubic shapes, particularly at the coarse red corner (see \autoref{fig:flattening_sphube_ptset}). Overall, the error profile of KRBF converges faster than that of the Best RBFs/HRBF approaches as indicated by the deeper blue color found  in the flat regions.

Finally, the results in the third and fourth columns are significantly more accurate than the PCA estimation, with KRBF achieving the highest accuracy. This indicates that KRBF can serve as a superior pre-processing step compared to PCA for subsequent surface reconstruction techniques.

\subsection*{Example 5: Normal Estimation of Real Point Cloud Data}
\reviewB{This example considers several realistic point clouds of size $N=100,000$ from the PCPNet dataset \cite{PCPNet} to evaluate the performance of the proposed approach. We first test the RBF, HRBF, and KRBF approaches on two selected point clouds with high curvature and fine geometric details, namely the bunny and the dragon. We then extend the comparison to a pillar-shaped point cloud that has flat top and bottom bases and a striped pattern along its side. For this point cloud, we also include PCA and a recently developed deep learning method, HSurf \cite{HSurf}.} For the HSurf approach, no additional training was performed on the model; we simply used the pre-trained model provided by the authors along with the default parameters from the code package.

\begin{figure}
    \centering
    \begin{overpic}[width=0.44\textwidth, trim=90 100 70 110, clip=true,tics=10]{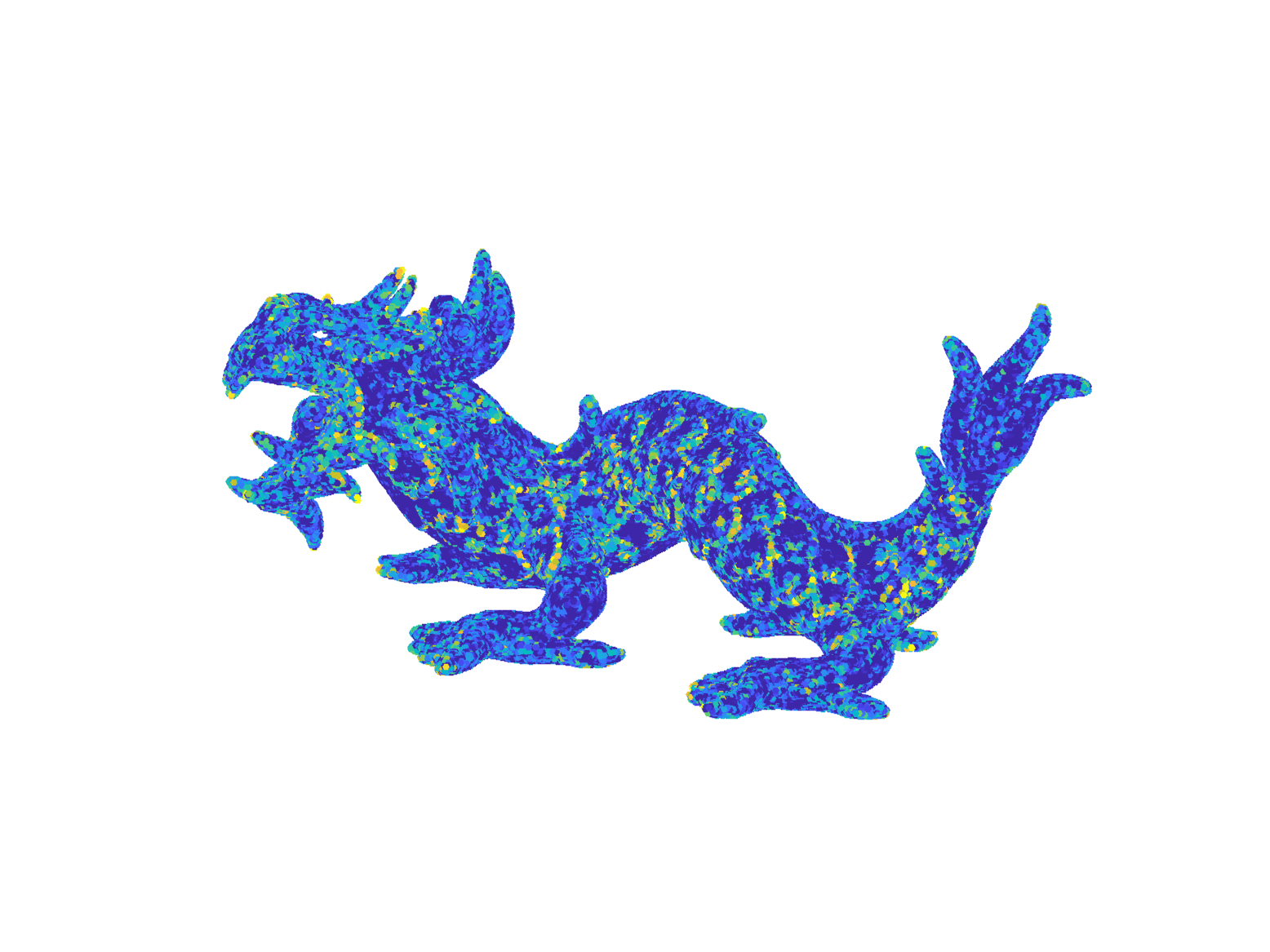}
    \put(8,66){\textbf{Best RBFs/HRBF}}
    \put(12,58){RMS = 0.19819}
    \end{overpic}
    \vspace{7mm}
    \begin{overpic}[width=0.44\textwidth, trim=90 100 70 110, clip=true,tics=10]{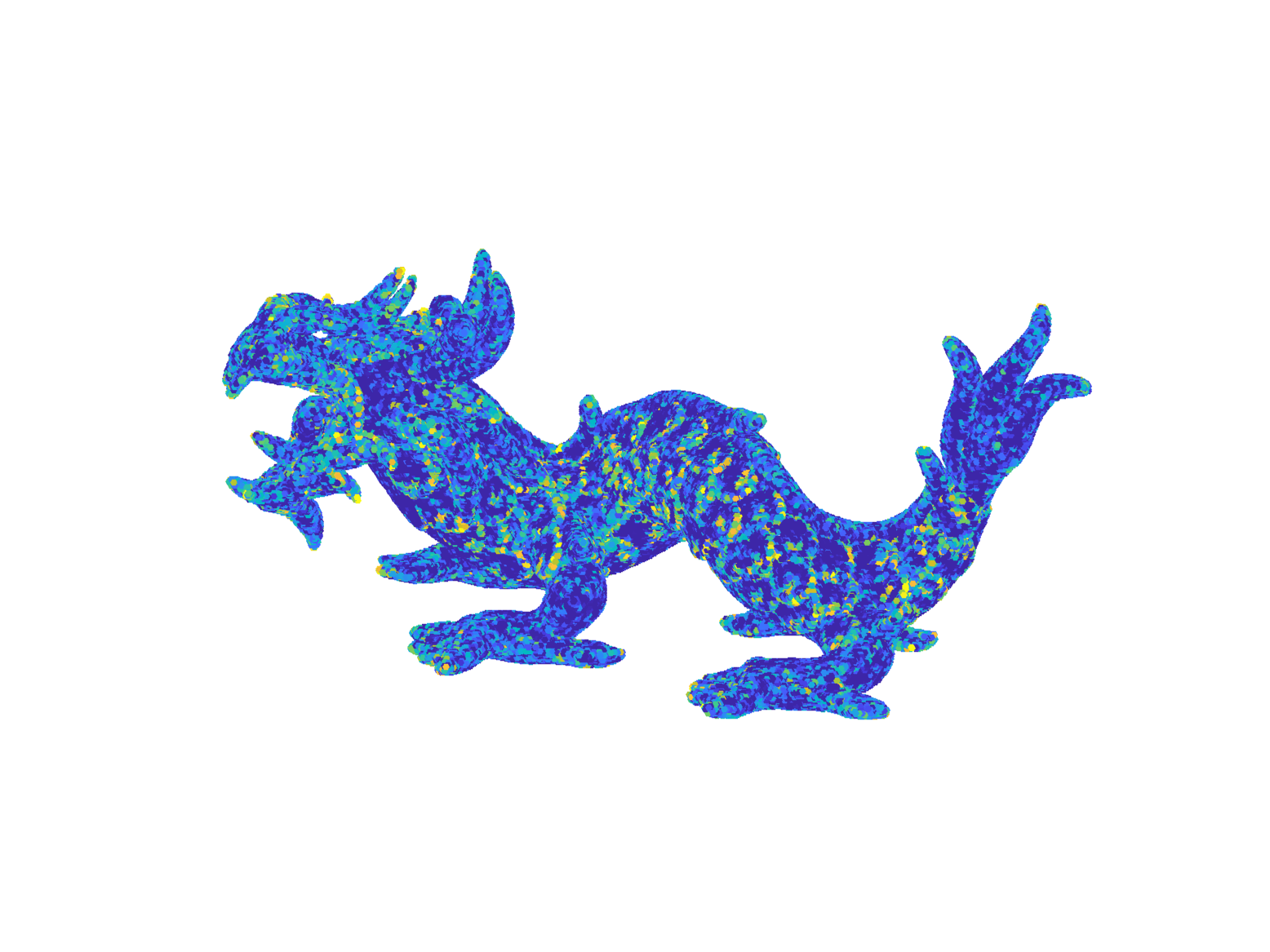}
    \put(32,66){\textbf{KRBF}}
    \put(24,58){RMS = 0.19706}
    \end{overpic}
    \\
    \begin{overpic}[width=0.44\textwidth, trim=100 60 100 50, clip=true,tics=10]{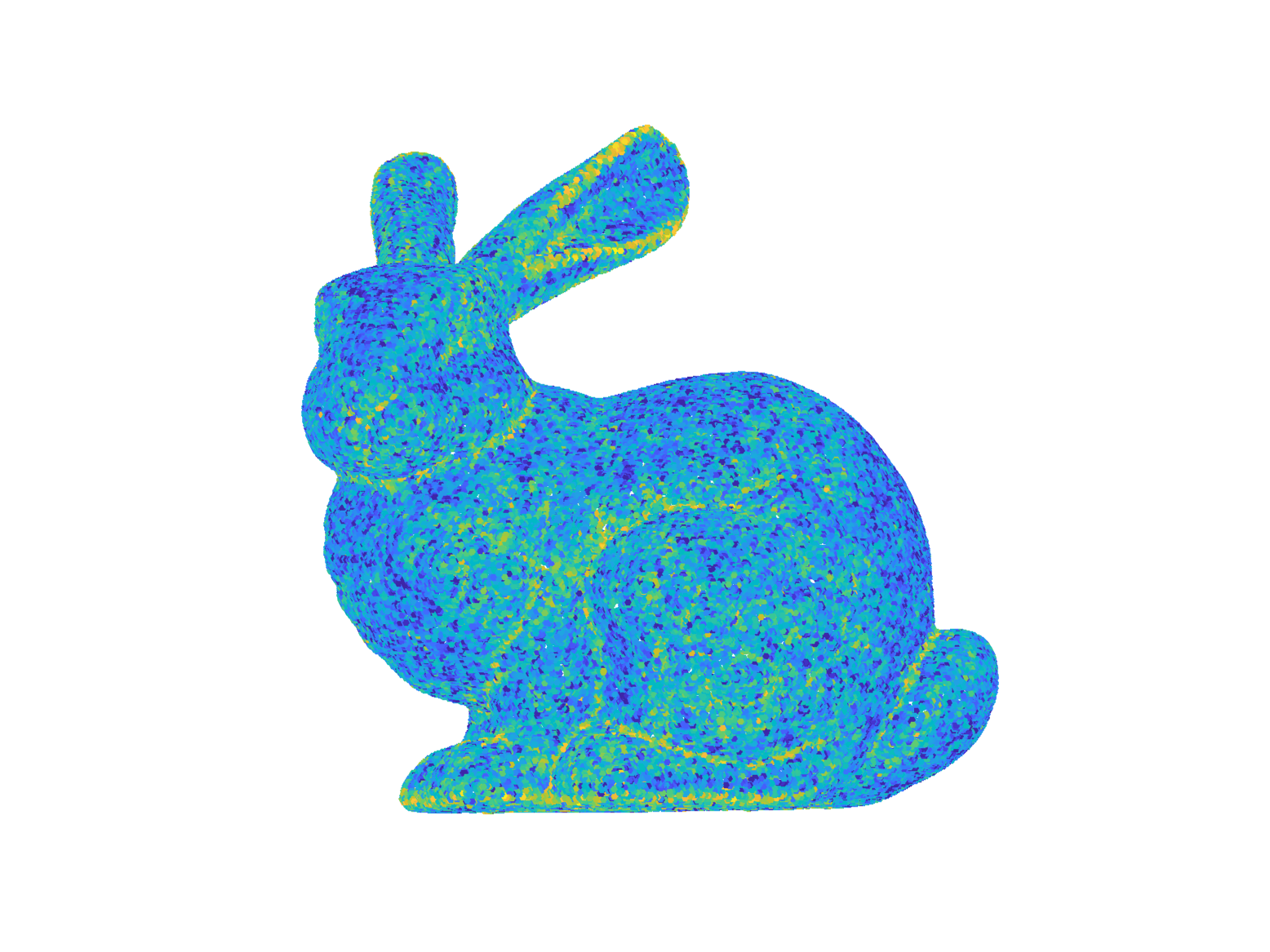}
    \put(12,87){RMS = 0.091331}
    \end{overpic}
    \vspace{7mm}
    \begin{overpic}[width=0.44\textwidth, trim=100 60 100 50, clip=true,tics=10]{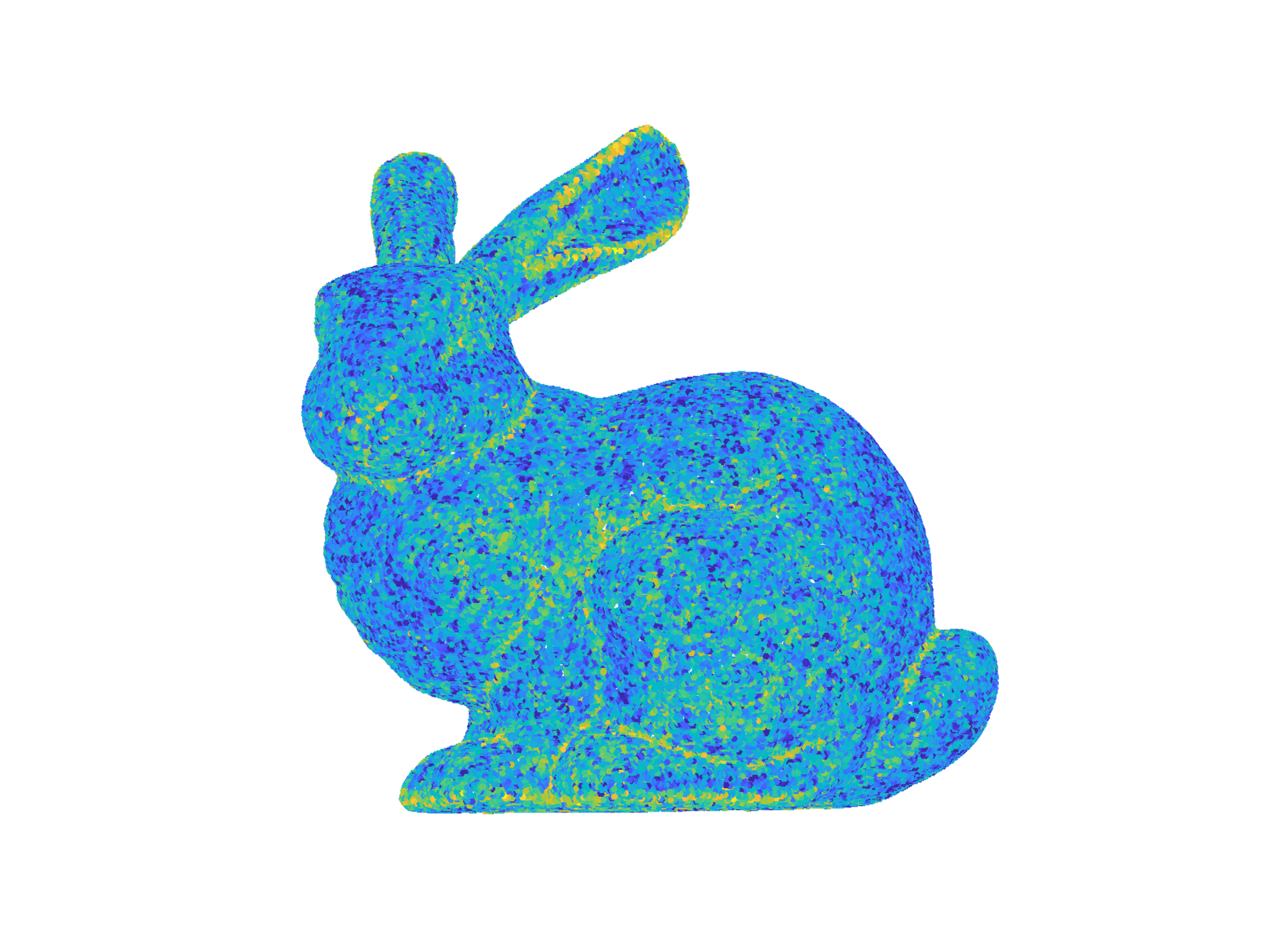}
    \put(24,87){RMS = 0.091202}
    \end{overpic}
    \caption{Root Mean Square (RMS) error of normal estimation on a bunny-shaped and a dragon-shaped point cloud of size $N=100,000$ points using Best RBFs/HRBF and KRBF approach.}
    \label{fig:bunny_dragon}
\end{figure}

\begin{figure}
    \centering
    \begin{minipage}[c]{0.49\textwidth}
        \begin{overpic}[width=\linewidth, trim=110 130 70 30, clip=true, tics=10]{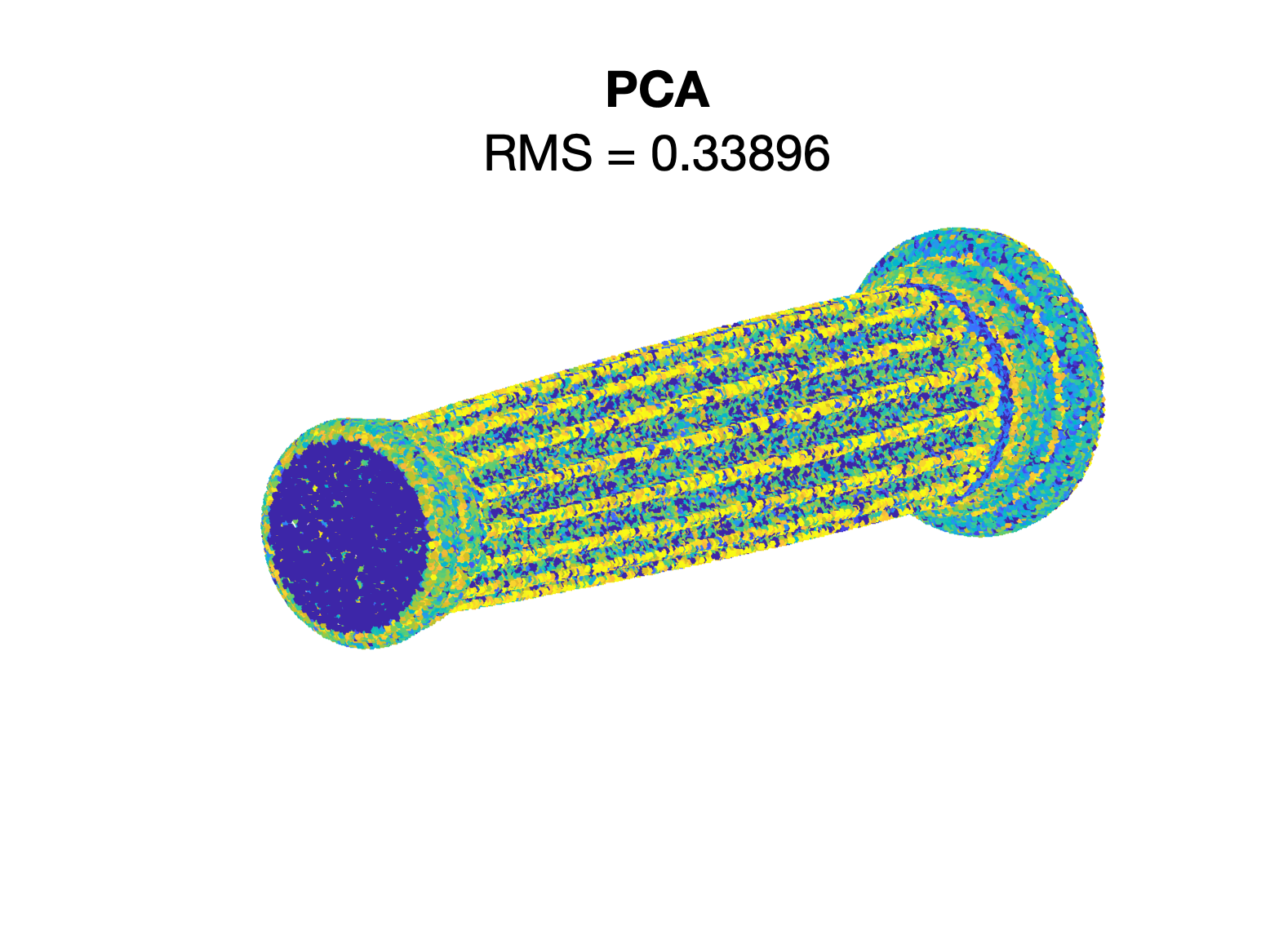}
            \put(40,15){
                \begin{tikzpicture}[overlay, x=\unitlength, y=\unitlength]
                    \draw[line width=1pt, black] (0,0) rectangle (30,20);
                \end{tikzpicture}
                }
            \put(70,35){\linethickness{0.25mm}\color{black}\line(1,0.47){35}}
            \put(70,15){\linethickness{0.25mm}\color{black}\line(1,0.05){35}}
        \end{overpic}
    \end{minipage}\hfill
    \begin{minipage}[c]{0.49\textwidth}
        \setlength{\fboxsep}{0pt}%
        \framebox{%
            \includegraphics[width=\linewidth, trim=60 140 40 120, clip=true]{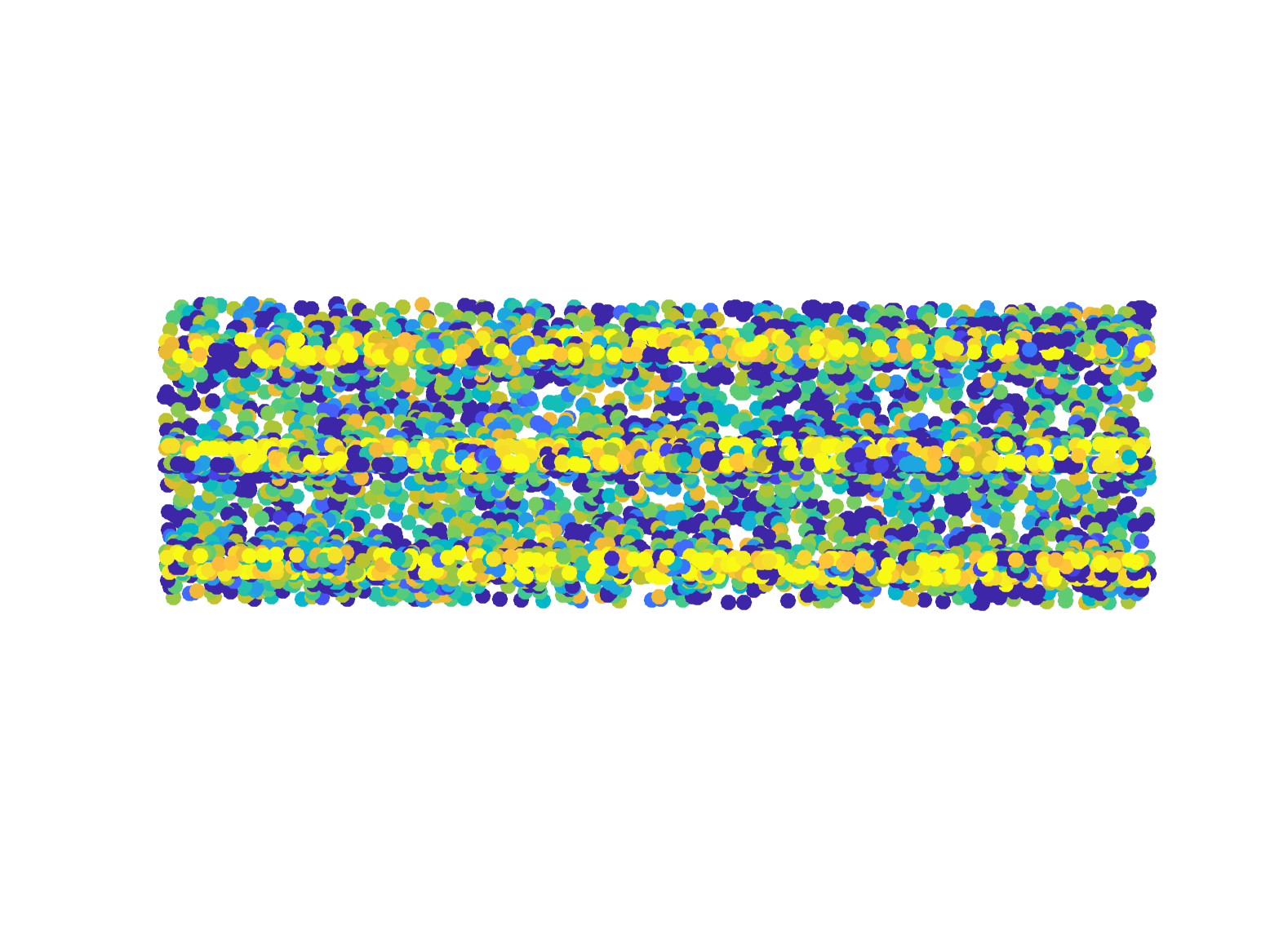}%
        }
    \end{minipage}
    \\
    \begin{minipage}[c]{0.49\textwidth}
        \begin{overpic}[width=\linewidth, trim=110 130 70 30, clip=true, tics=10]{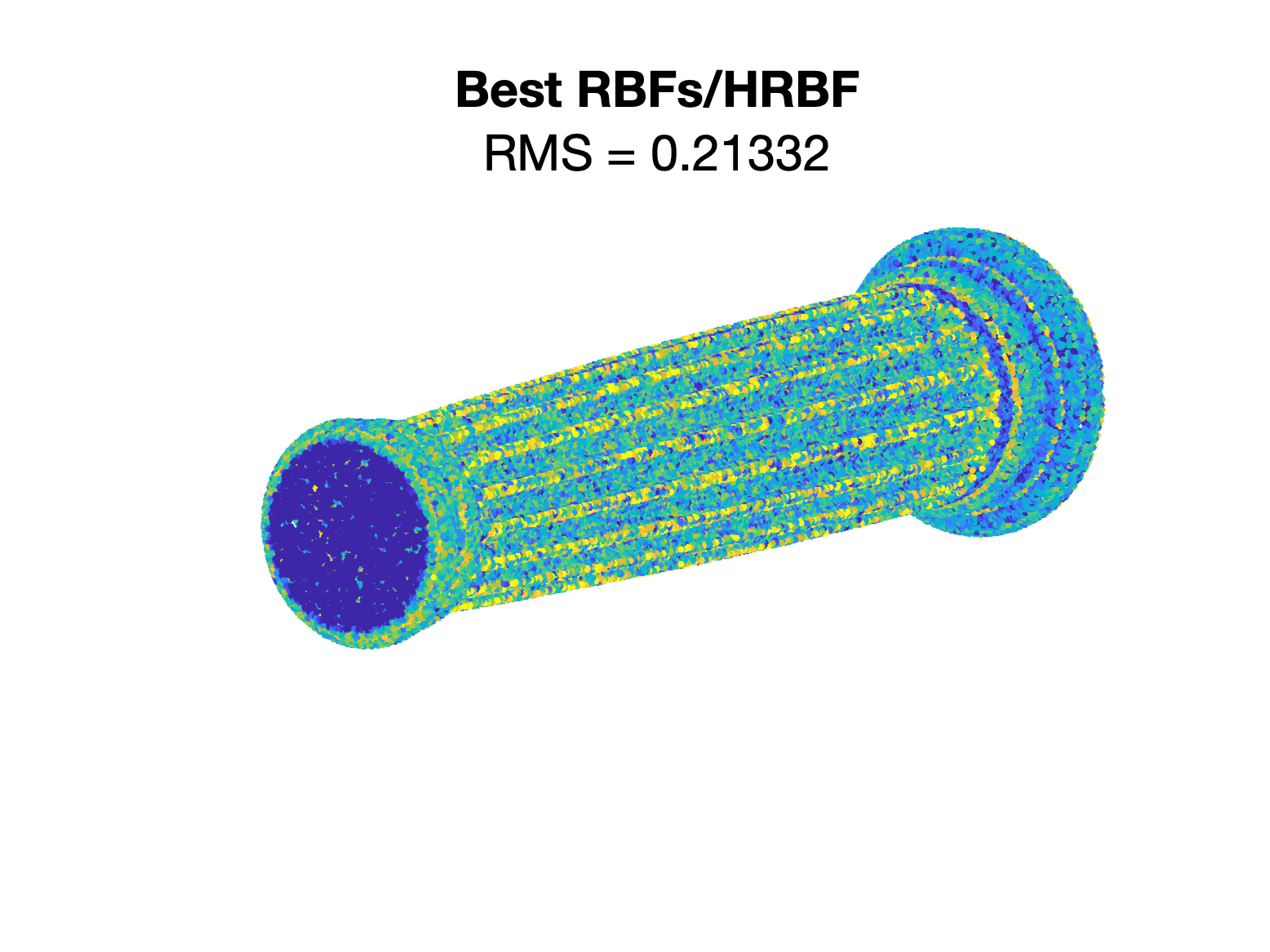}
            \put(40,15){
                \begin{tikzpicture}[overlay, x=\unitlength, y=\unitlength]
                    \draw[line width=1pt, black] (0,0) rectangle (30,20);
                \end{tikzpicture}
                }
            \put(70,35){\linethickness{0.25mm}\color{black}\line(1,0.47){35}}
            \put(70,15){\linethickness{0.25mm}\color{black}\line(1,0.05){35}}
        \end{overpic}
    \end{minipage}\hfill
    \begin{minipage}[c]{0.49\textwidth}
        \setlength{\fboxsep}{0pt}%
        \framebox{%
            \includegraphics[width=\linewidth, trim=60 140 40 120, clip=true]{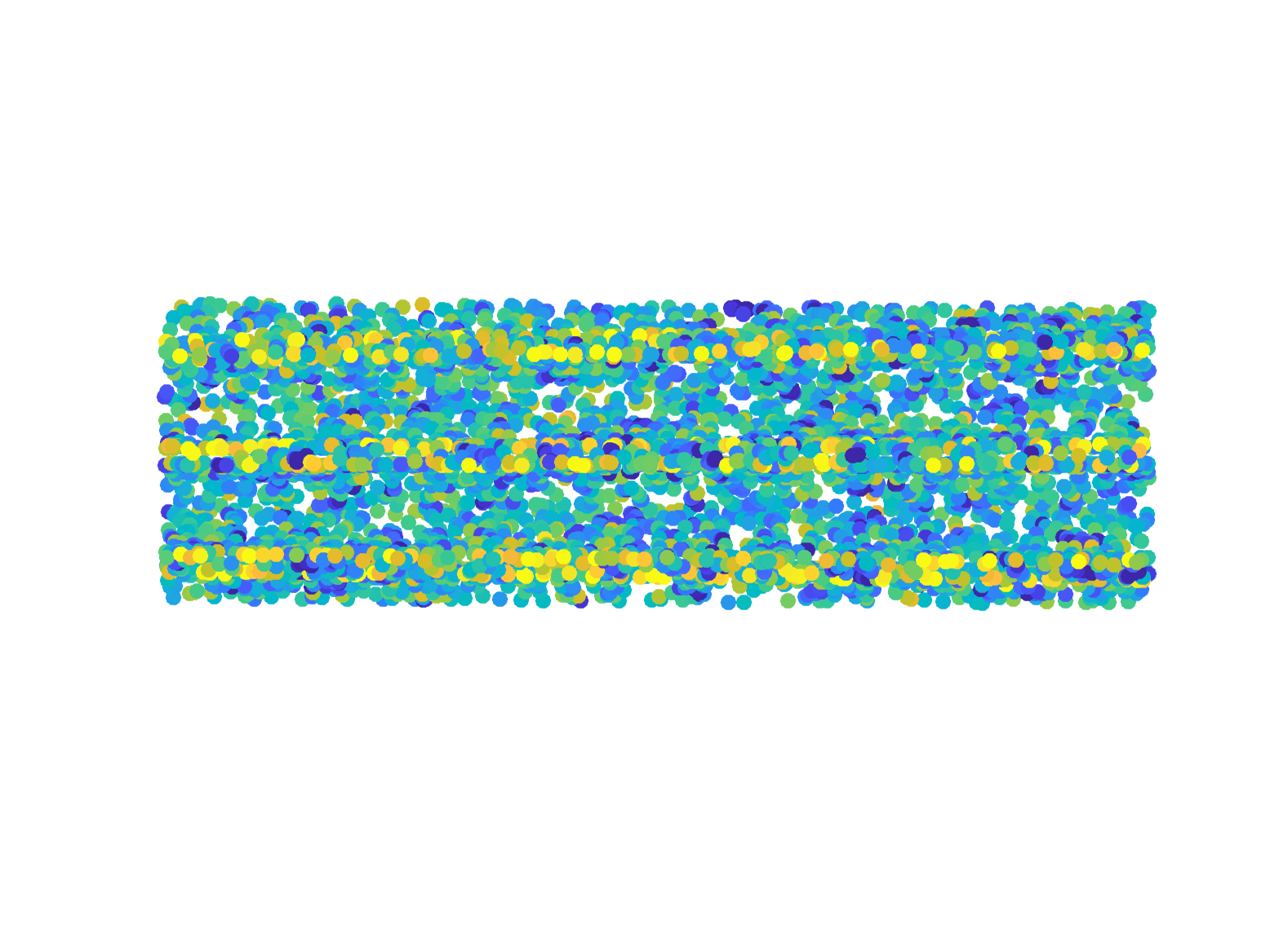}%
        }
    \end{minipage}
    \\
    \begin{minipage}[c]{0.49\textwidth}
        \begin{overpic}[width=\linewidth, trim=110 130 70 30, clip=true, tics=10]{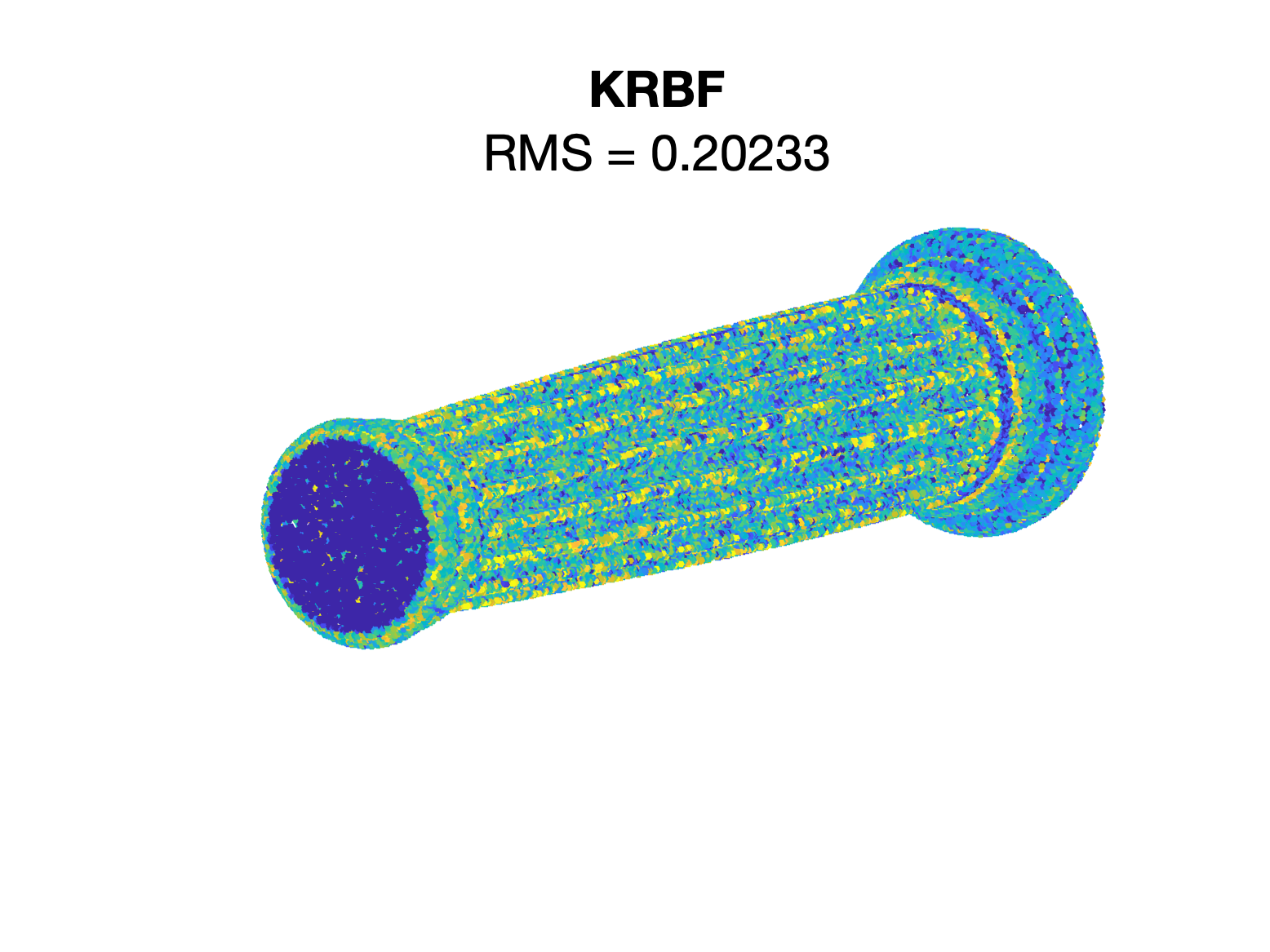}
            \put(40,15){
                \begin{tikzpicture}[overlay, x=\unitlength, y=\unitlength]
                    \draw[line width=1pt, black] (0,0) rectangle (30,20);
                \end{tikzpicture}
                }
            \put(70,35){\linethickness{0.25mm}\color{black}\line(1,0.47){35}}
            \put(70,15){\linethickness{0.25mm}\color{black}\line(1,0.05){35}}
        \end{overpic}
    \end{minipage}\hfill
    \begin{minipage}[c]{0.49\textwidth}
        \setlength{\fboxsep}{0pt}%
        \framebox{%
            \includegraphics[width=\linewidth, trim=60 140 40 120, clip=true]{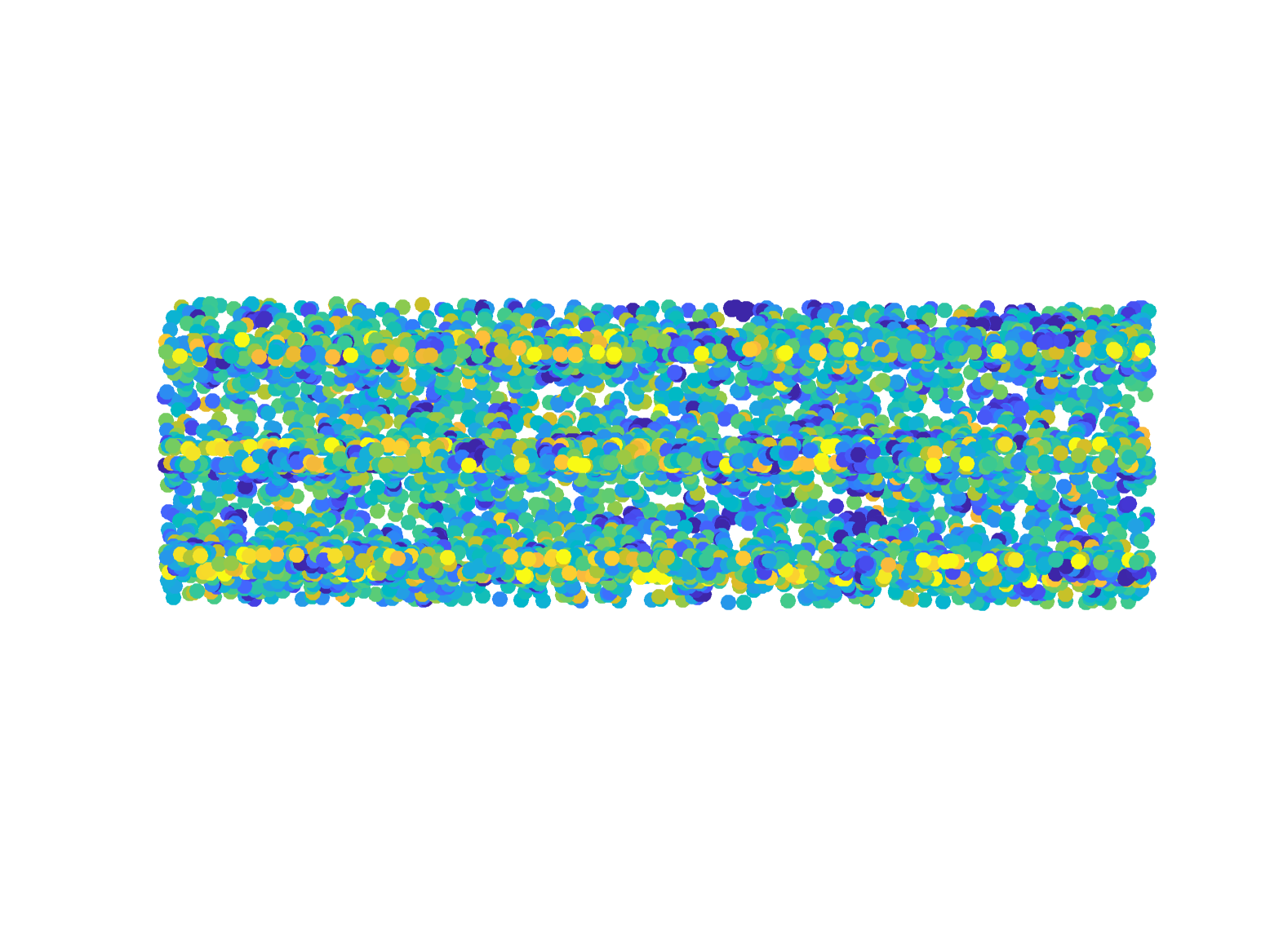}%
        }
    \end{minipage}
    \\
    \begin{minipage}[c]{0.49\textwidth}
        \begin{overpic}[width=\linewidth, trim=110 130 70 30, clip=true, tics=10]{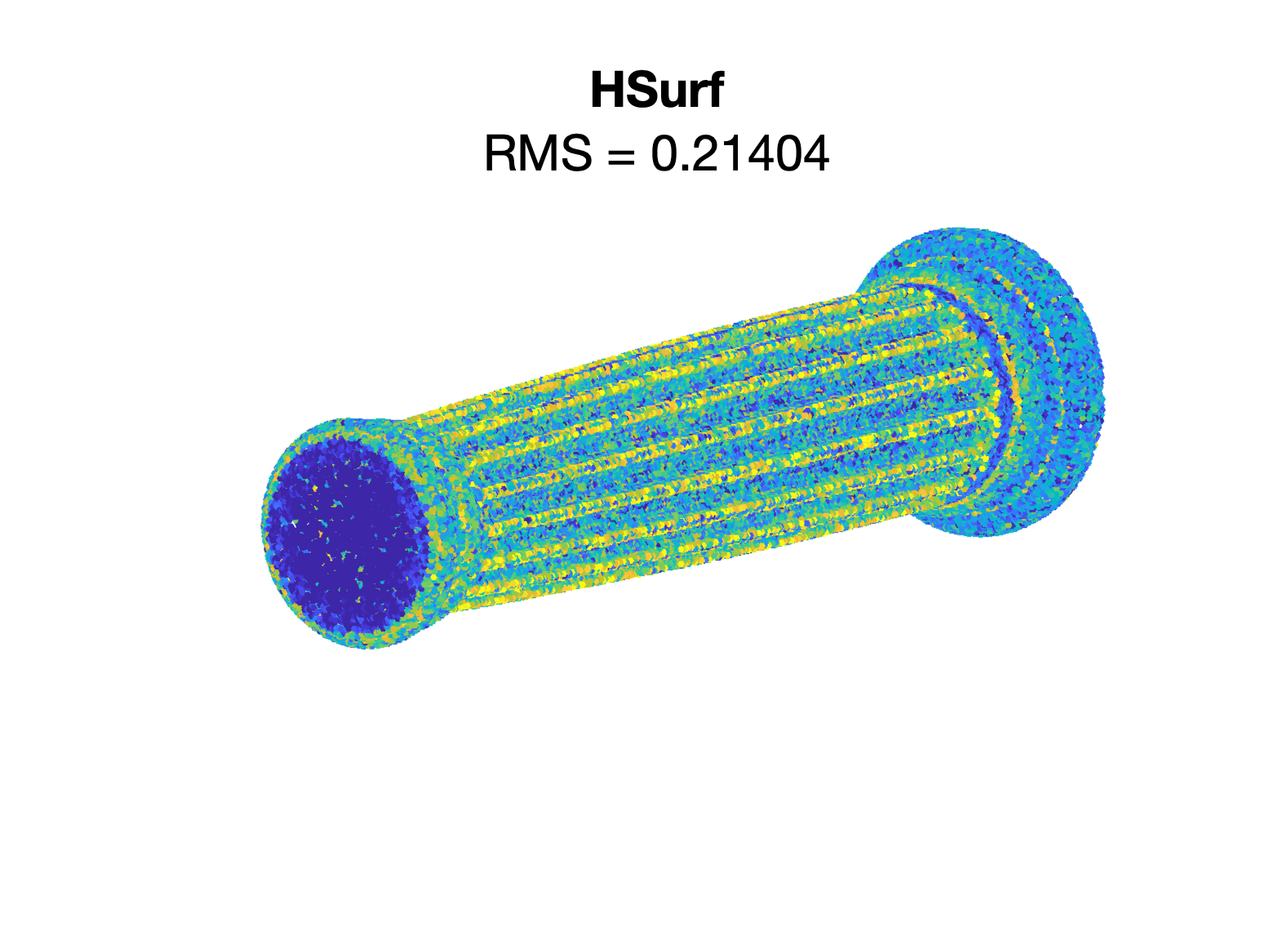}
            \put(40,15){
                \begin{tikzpicture}[overlay, x=\unitlength, y=\unitlength]
                    \draw[line width=1pt, black] (0,0) rectangle (30,20);
                \end{tikzpicture}
                }
            \put(70,35){\linethickness{0.25mm}\color{black}\line(1,0.47){35}}
            \put(70,15){\linethickness{0.25mm}\color{black}\line(1,0.05){35}}
        \end{overpic}
    \end{minipage}\hfill
    \begin{minipage}[c]{0.49\textwidth}
        \setlength{\fboxsep}{0pt}%
        \framebox{%
            \includegraphics[width=\linewidth, trim=60 140 40 120, clip=true]{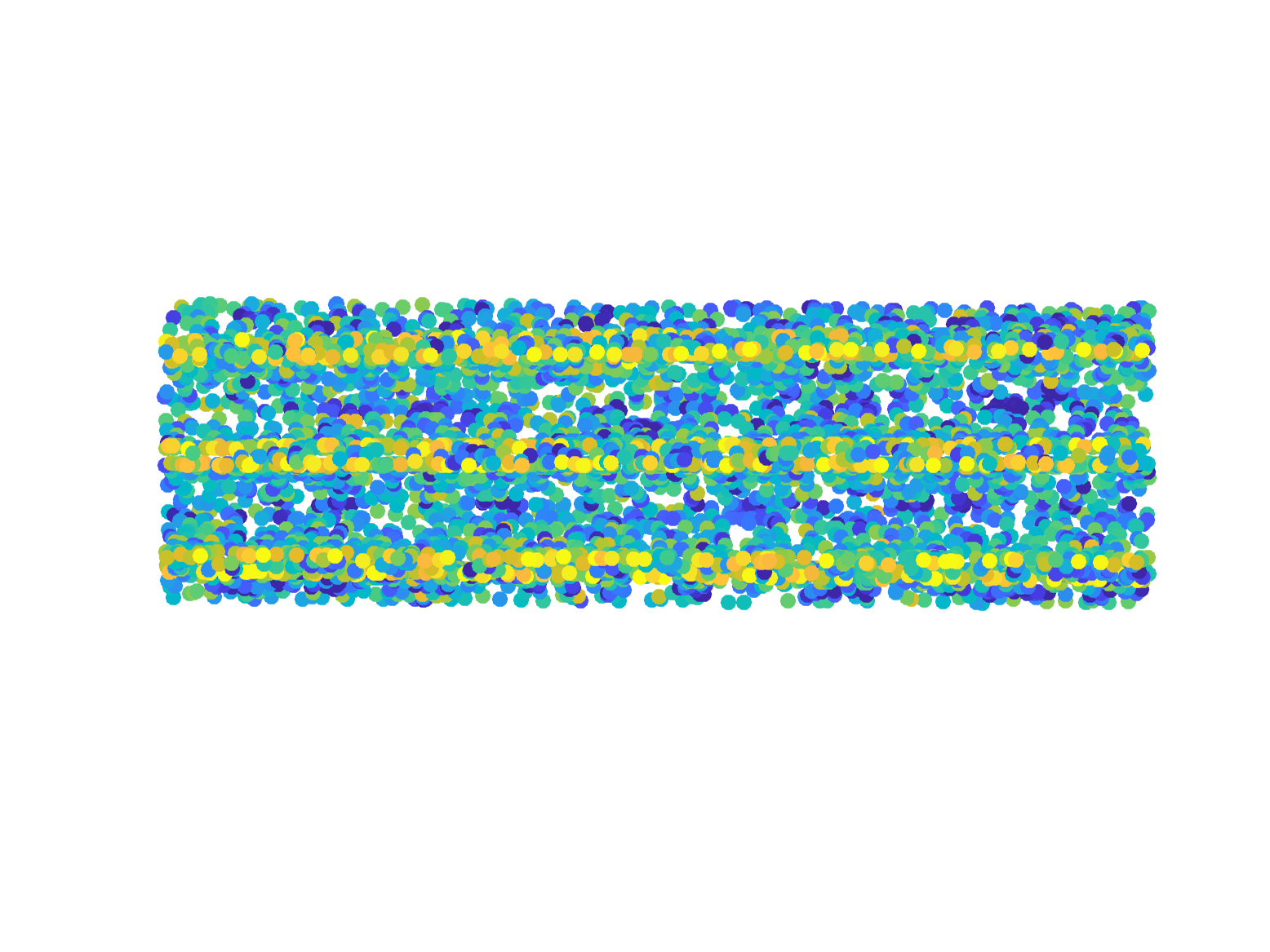}%
        }
    \end{minipage}
    \caption{Root Mean Square (RMS) error of normal estimation on a pillar-shaped point cloud of size $N=100,000$ points using various approaches.}
    \label{fig:real_ptcld}
\end{figure}

\reviewB{\autoref{fig:bunny_dragon} shows the best normal estimation results from RBF and HRBF alongside KRBF on the bunny and dragon, with KRBF achieving slightly better results than RBF and HRBF.}
The results from the different approaches on the pillar example are shown in \autoref{fig:real_ptcld}, with KRBF achieving the highest accuracy and PCA remaining the least accurate. Rather than using the maximum error $\| \mathcal{E}_{\n}\|_{\infty}$, we compare accuracy using the root mean squared error (RMS) because these point cloud data have complex surfaces and the maximum error is similar across all four methods. Although KRBF, RBF, and HSurf have the same order of accuracy, the RBF and KRBF approaches yield more accurate results in the flat region of the top base (indicated by solid blue), whereas HSurf displays a lighter blue ring around the base. Moreover, KRBF shows a more consistent blue pattern along the striped structure, which contributes to its lower RMS error compared to the other methods.

These examples demonstrate that KRBF can capture fine geometric details in complex, real-world point cloud data more effectively than traditional methods while maintaining comparable accuracy to deep learning method.

\subsection*{Example 6: Point Cloud Processing}
For the last example, we demonstrate the potential of the proposed KAN-inspired trial space and minimum-norm interpolation to perform point cloud processing. We consider a noisy point cloud sampled from the unit cubic point with Gaussian noise of standard deviation $0.02$ in each data point.
    
First, we use the proposed minimum-norm interpolation to estimate curvatures on the noisy point cloud for subsequent point-cloud processing. For comparison, we also employ standard Tikhonov regularization with native space norm regularization to find an approximant from \(\mathcal{U}_{\text{KAN}}\) for curvature estimation.
Subsequently, we use some anisotropic Laplacian operator based on the computed curvatures to perform fairing.
This procedure is similar to the bilateral filtering methods described in \cite{bilateral_denoise,anisotropic_search}, but with key differences: we use the maximum principal curvature \(\kappa_1\) to define the weight function instead of a {similarity weight function} that penalizes large variations in intensity. Additionally, instead of updating points along the normal direction, we update their positions using a weighted average location, allowing points to move freely in any direction during each iteration.

\begin{figure}
    \centering
    \begin{overpic}[width=0.24\textwidth, trim=100 40 100 20, clip=true,tics=10]{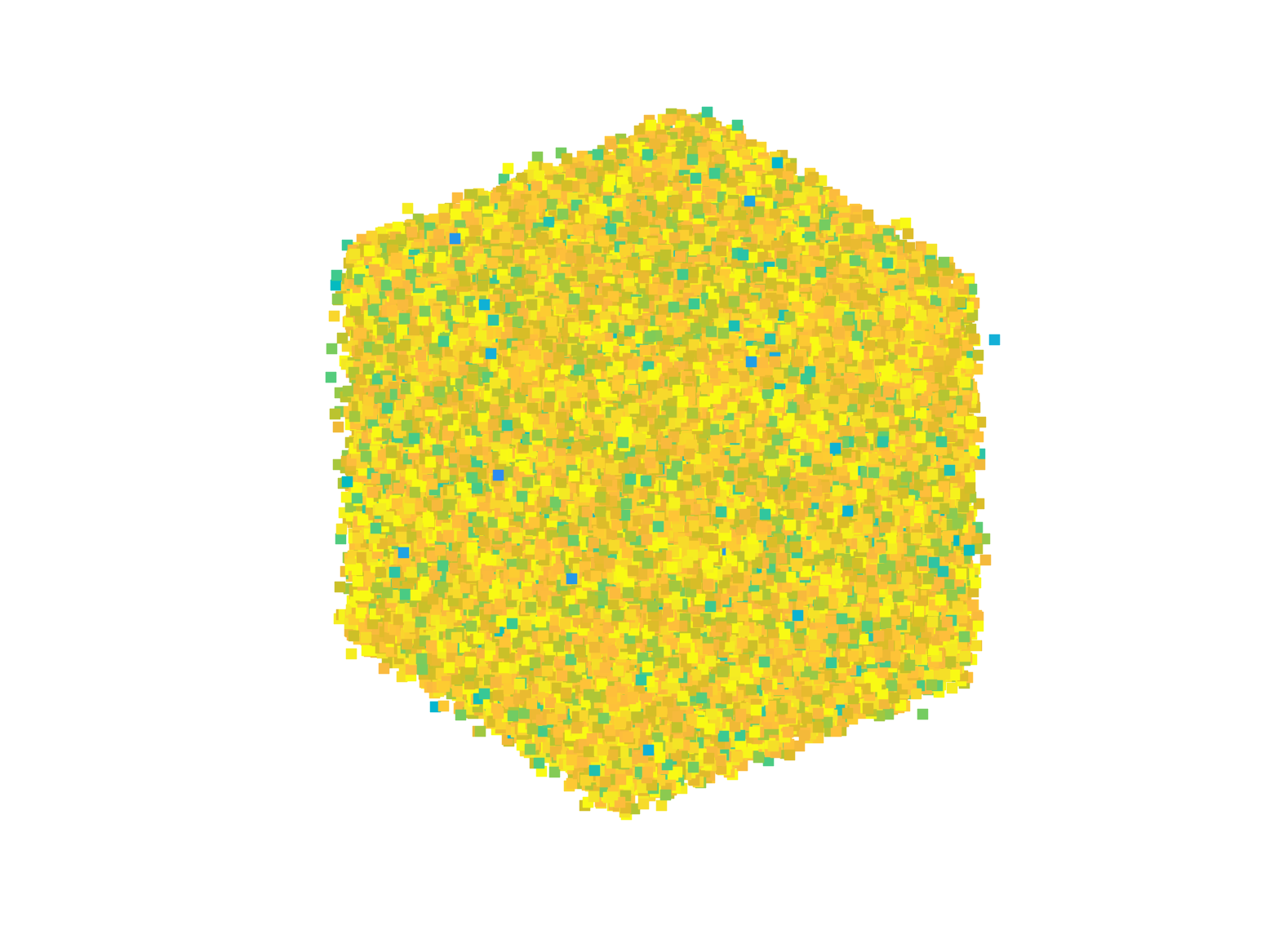}
        \put(30,120){\textbf{Minimum-norm interpolation}}
        \put(50,100){$\kappa_1$}
        \put(0,40){\rotatebox{90}{\small{Noisy data}}}
    \end{overpic}
    \begin{overpic}[width=0.24\textwidth, trim=100 40 100 20, clip=true,tics=10]{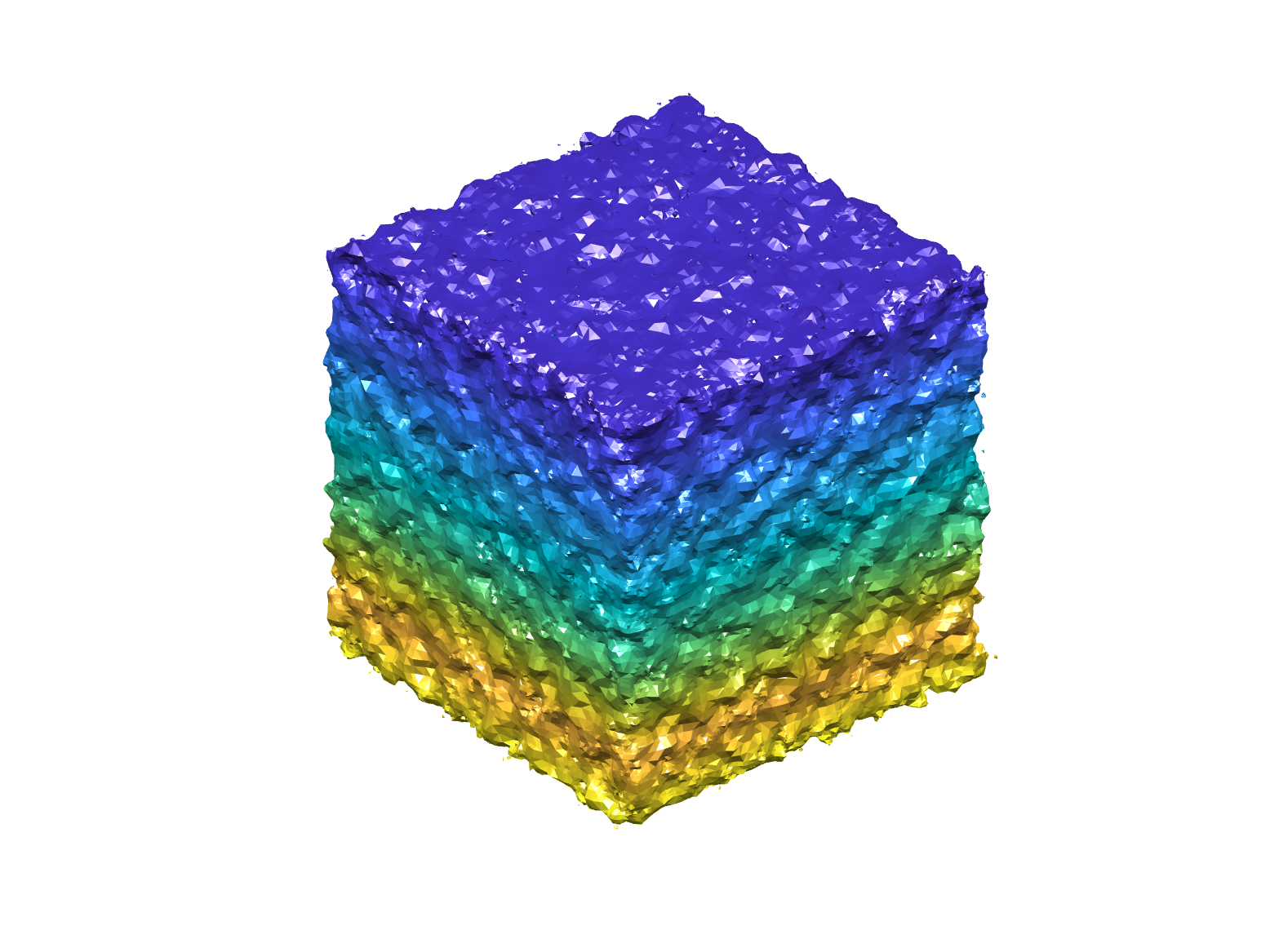}
        \put(35,100){surface}
    \end{overpic}
    \begin{overpic}[width=0.24\textwidth, trim=100 40 100 20, clip=true,tics=10]{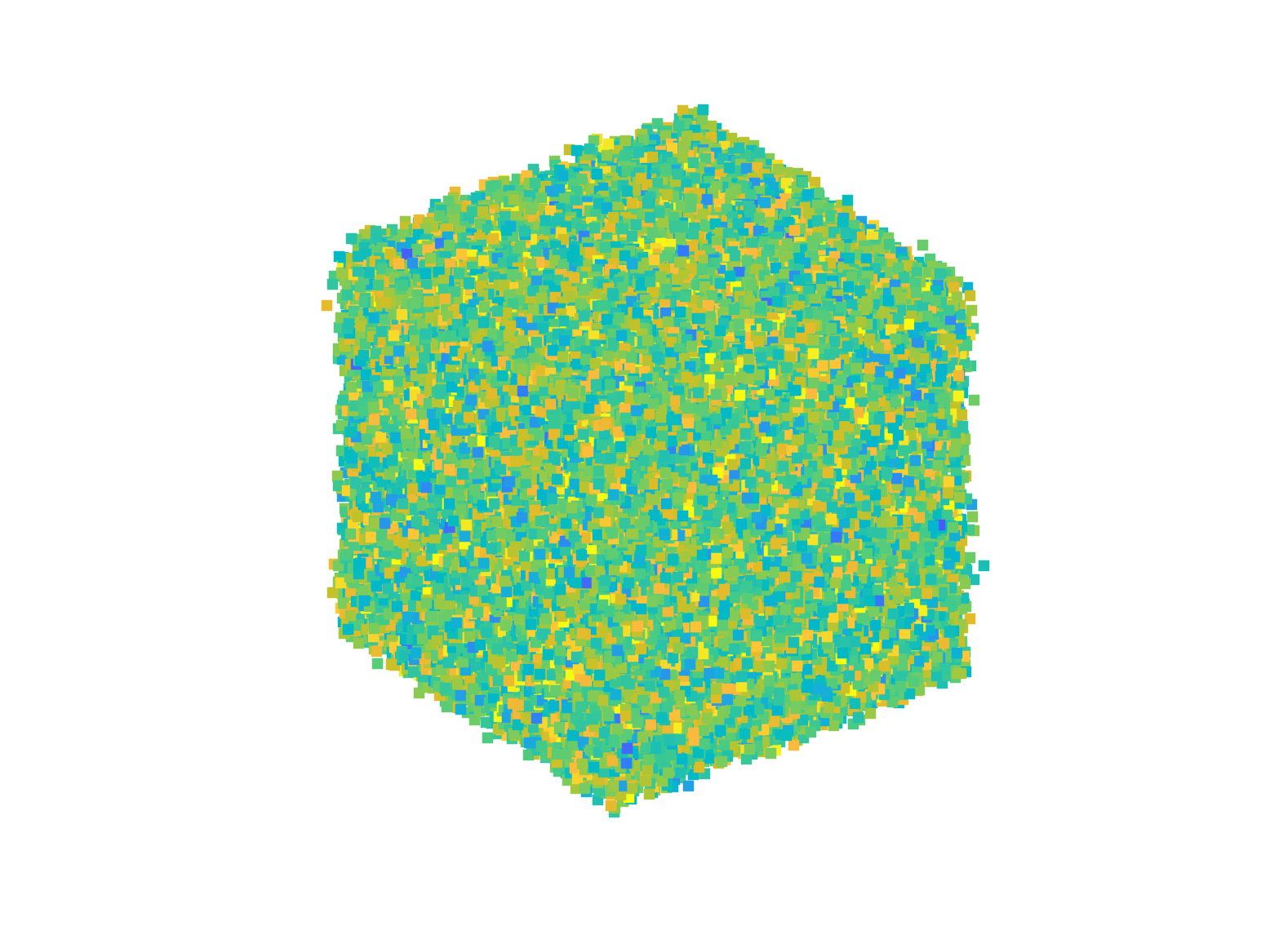}
        \put(55,120){\textbf{Tikhonov regularization }}
        \put(50,100){$\kappa_1$}
    \end{overpic}
    \begin{overpic}[width=0.24\textwidth, trim=100 40 100 20, clip=true,tics=10]{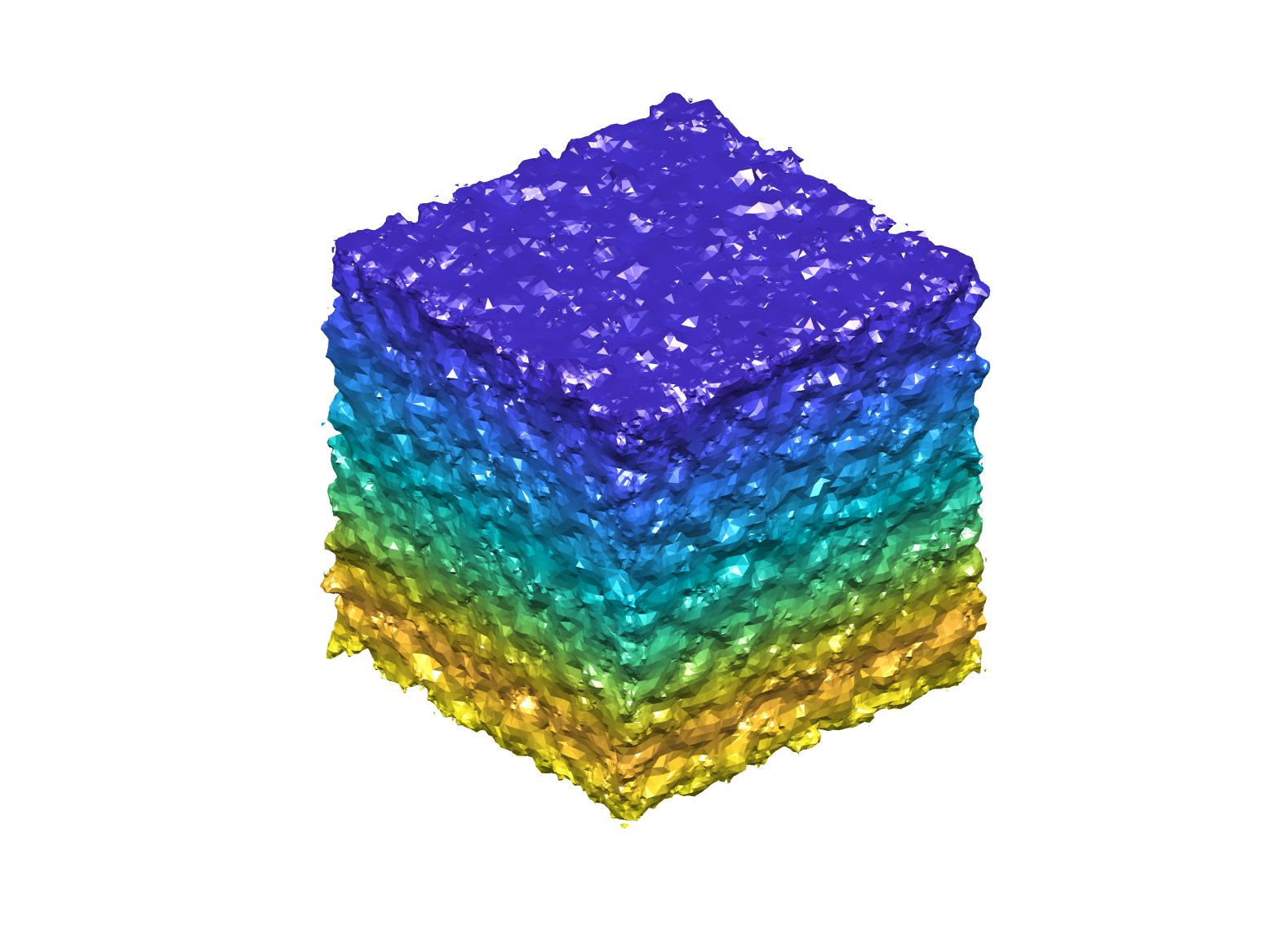}
        \put(35,100){surface}
    \end{overpic}
    \\
    \begin{overpic}[width=0.24\textwidth, trim=100 40 100 20, clip=true,tics=10]{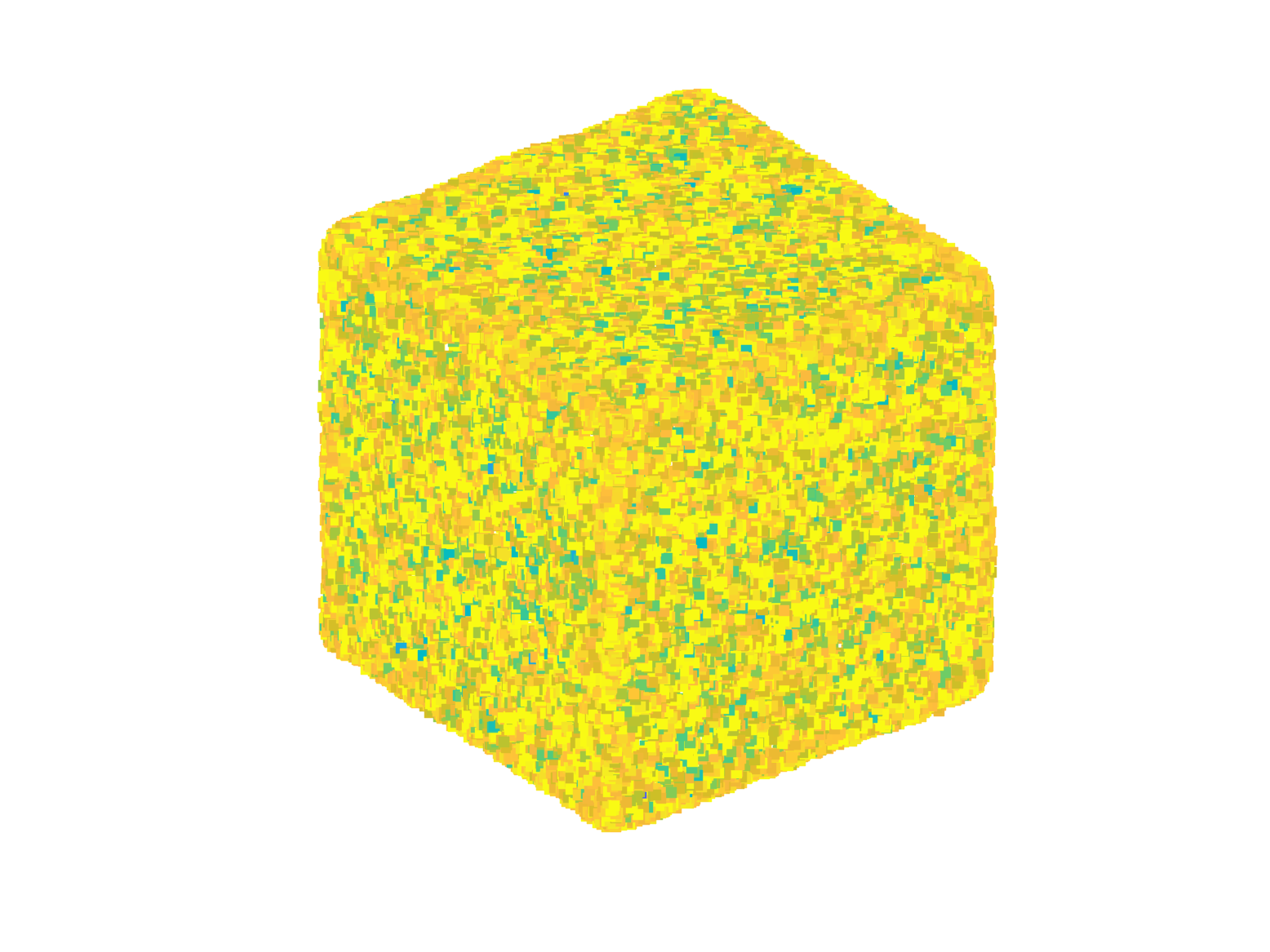}
        \put(0,40){\rotatebox{90}{\small{1 fairing}}}
    \end{overpic}
    \begin{overpic}[width=0.24\textwidth, trim=100 40 100 20, clip=true,tics=10]{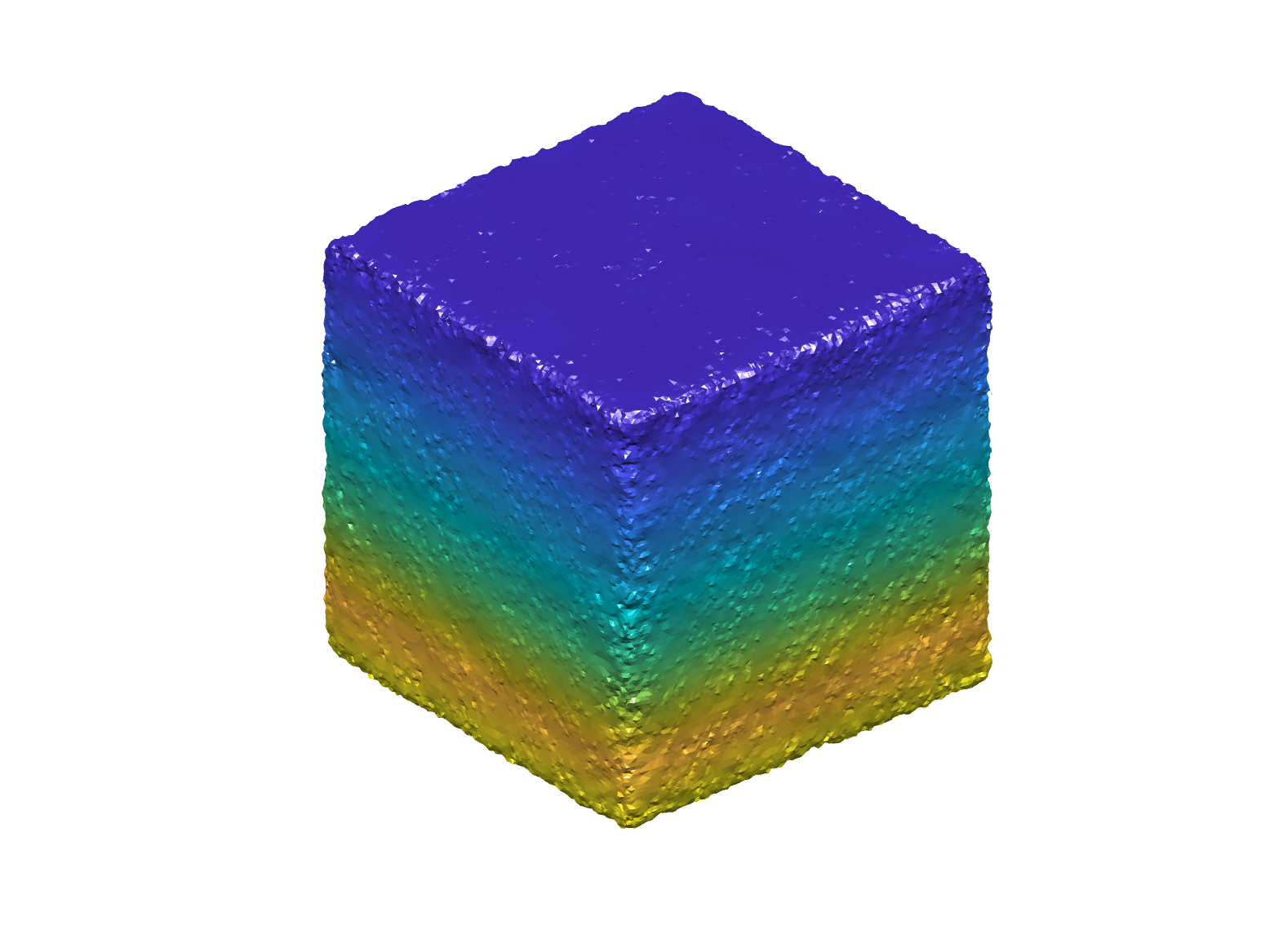}
    \end{overpic}
    \begin{overpic}[width=0.24\textwidth, trim=100 40 100 20, clip=true,tics=10]{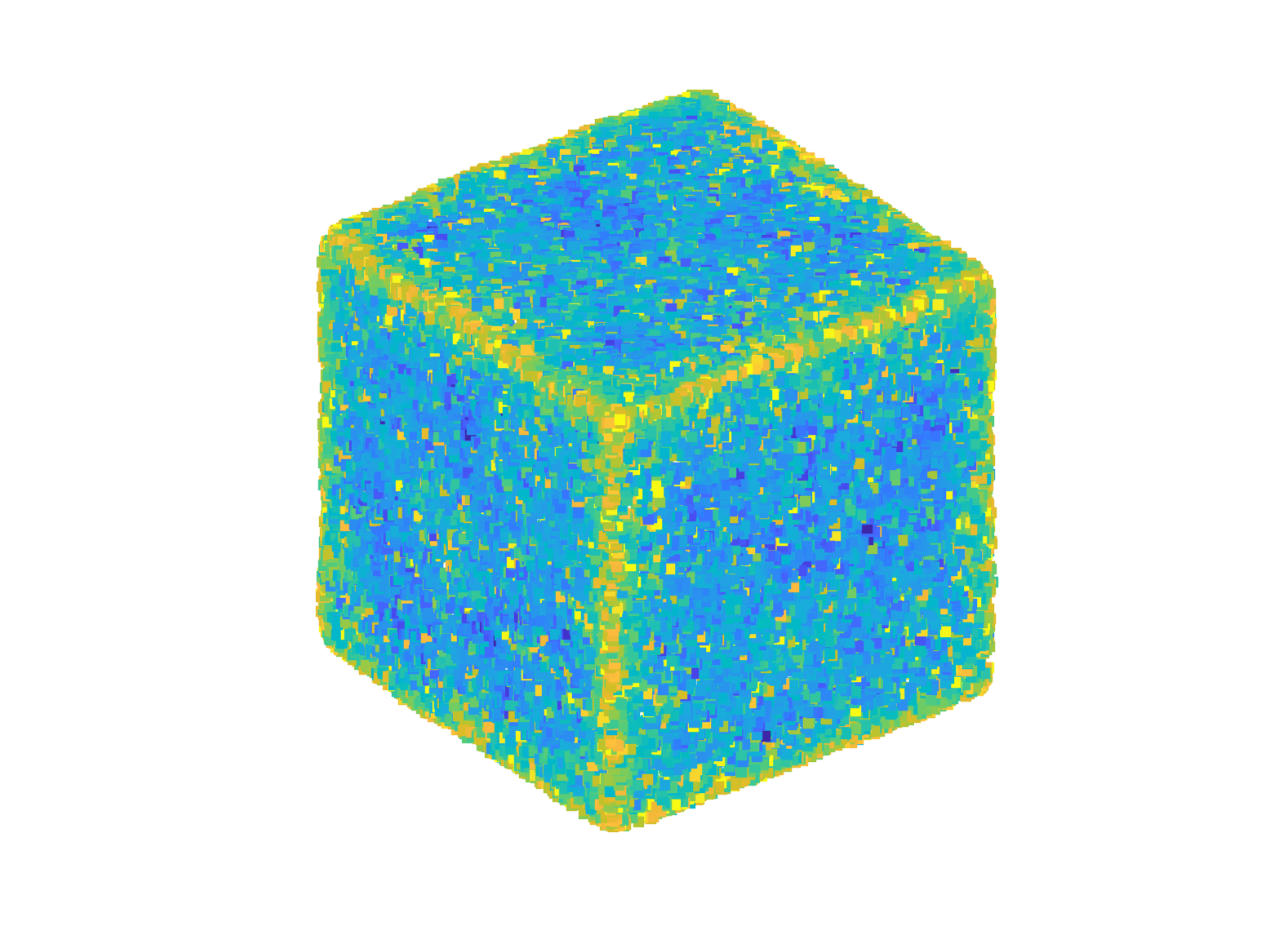}
    \end{overpic}
    \begin{overpic}[width=0.24\textwidth, trim=100 40 100 20, clip=true,tics=10]{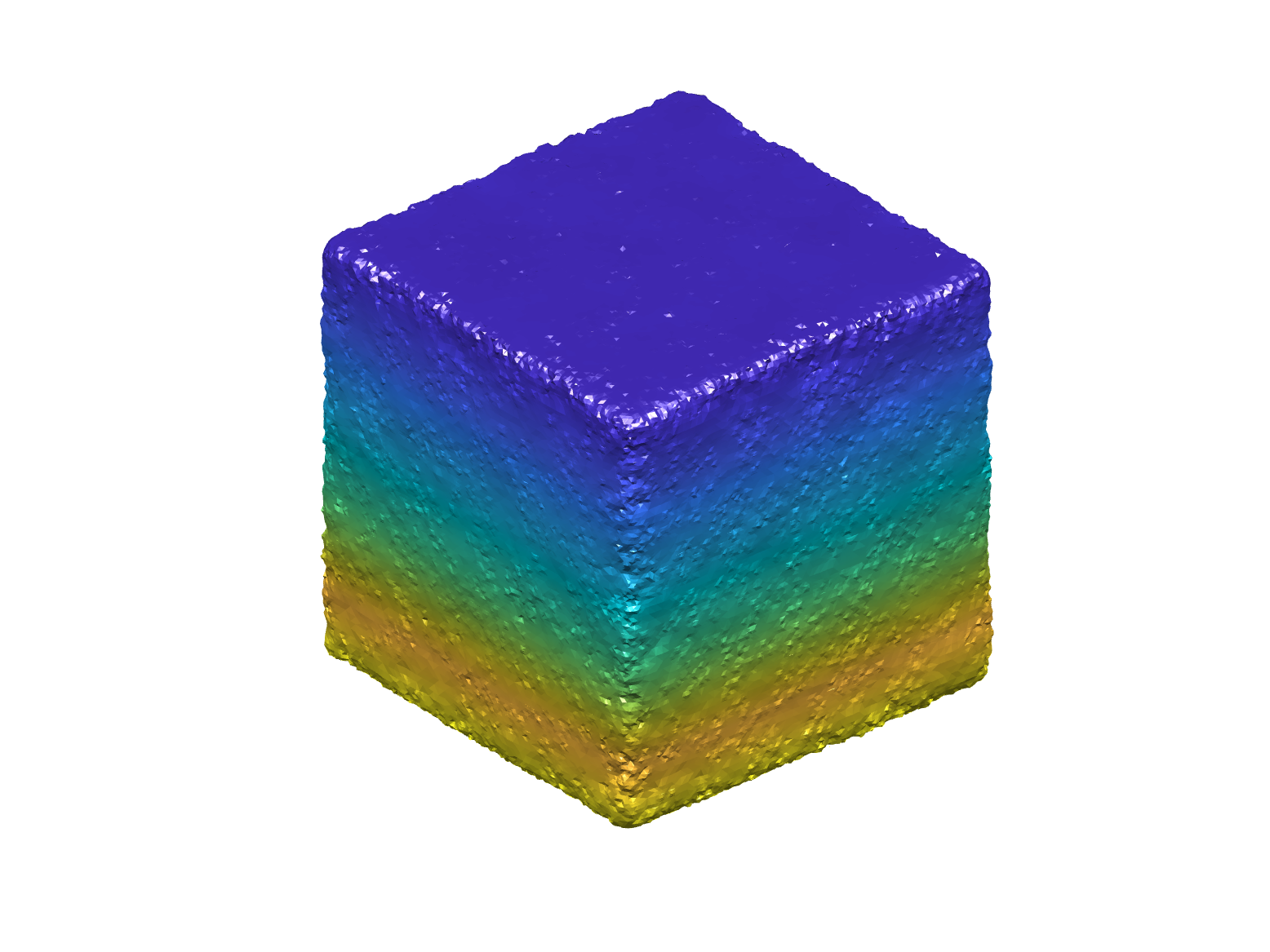}
    \end{overpic}
    \\
    \begin{overpic}[width=0.24\textwidth, trim=100 40 100 20, clip=true,tics=10]{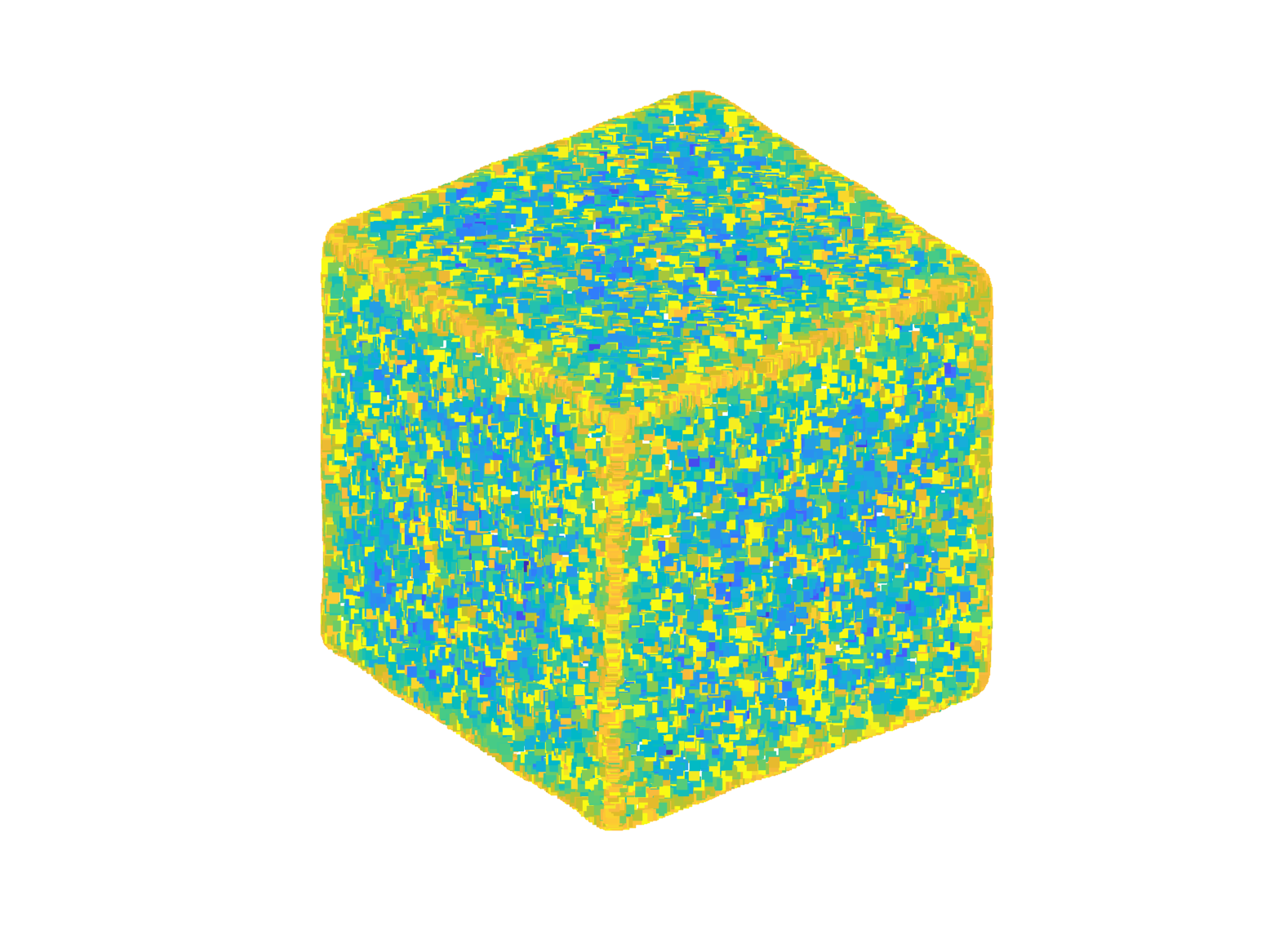}
        \put(0,40){\rotatebox{90}{\small{2 fairings}}}
    \end{overpic}
    \begin{overpic}[width=0.24\textwidth, trim=100 40 100 20, clip=true,tics=10]{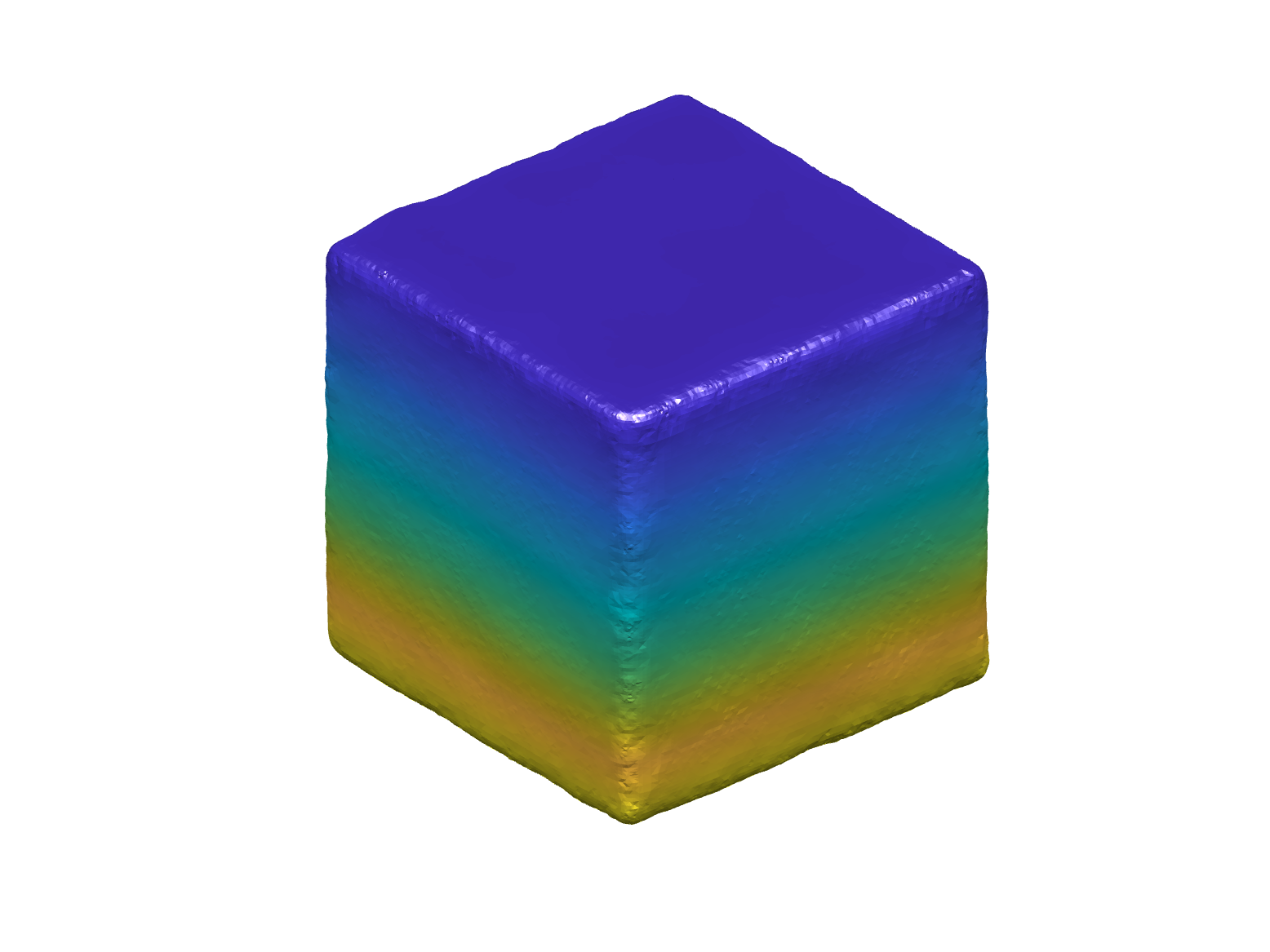}
    \end{overpic}
    \begin{overpic}[width=0.24\textwidth, trim=100 40 100 20, clip=true,tics=10]{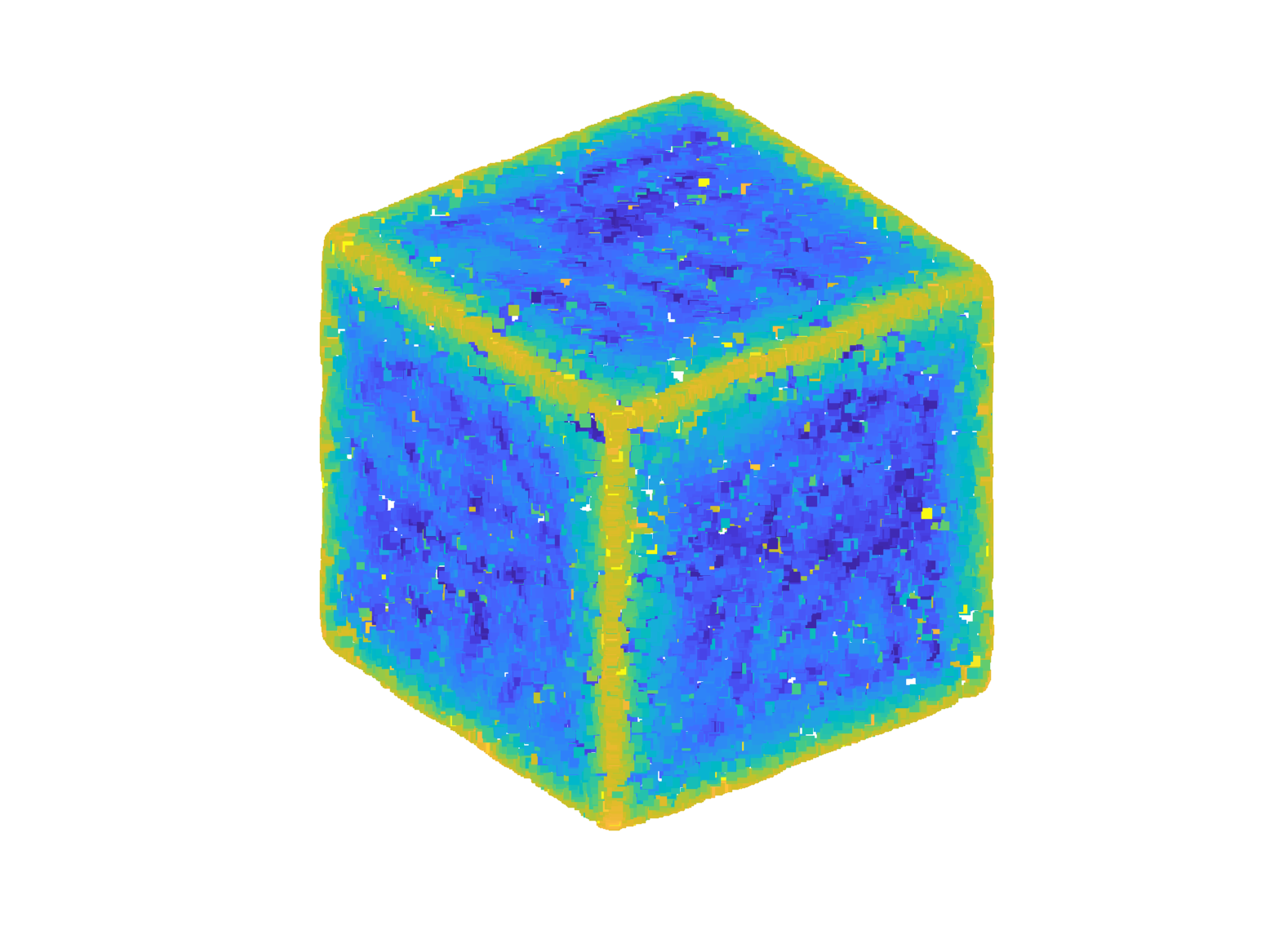}
    \end{overpic}
    \begin{overpic}[width=0.24\textwidth, trim=100 40 100 20, clip=true,tics=10]{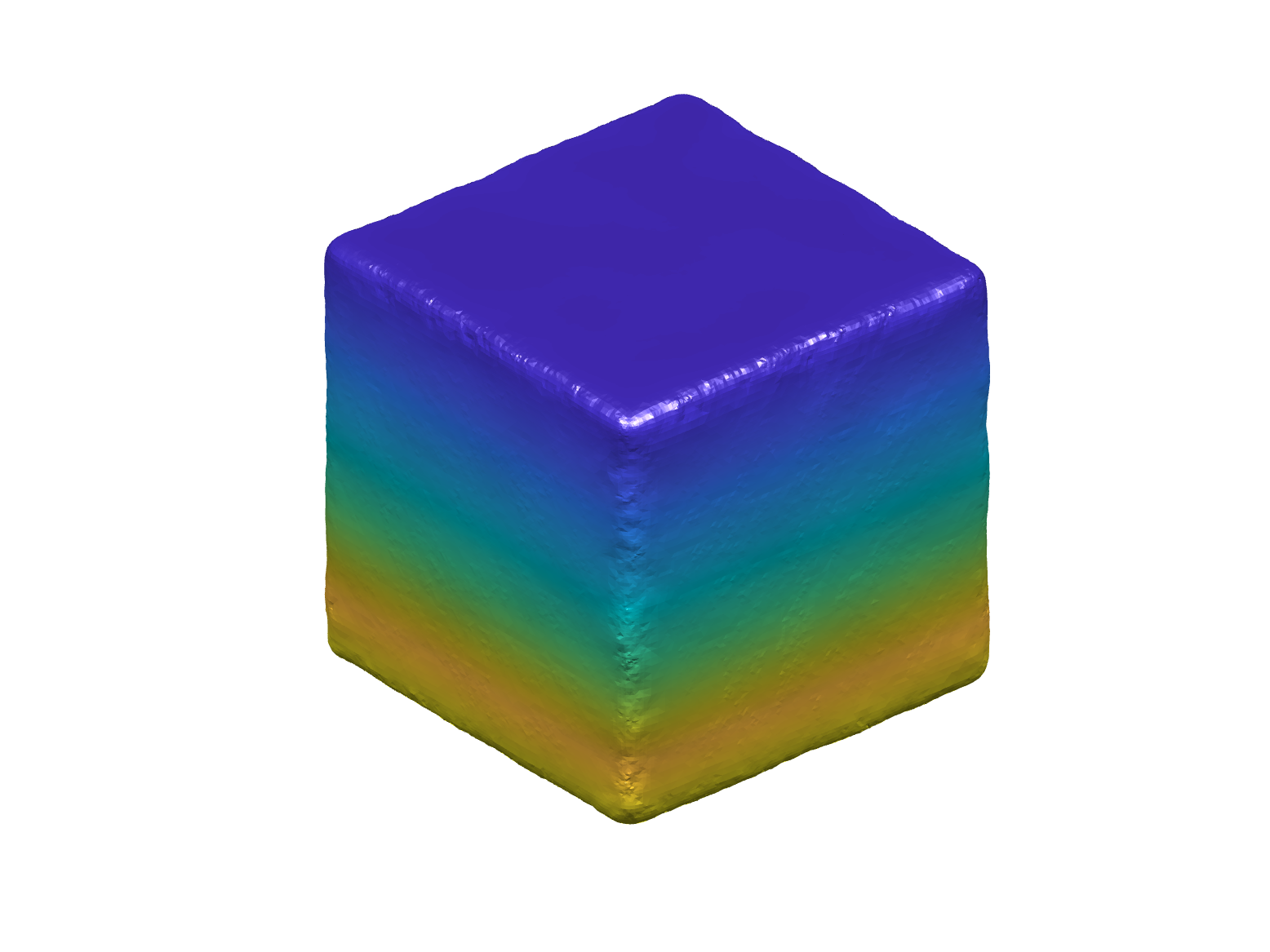}
    \end{overpic}
    \caption{Fairing of a point cloud consisting of $\approx30,000$ points forming a unit cube, subjected to Gaussian noise with a standard deviation of $0.02$. Fairing is performed using a Sobolev kernel with a smoothness order of $\tau=3$ and a local stencil size of $N_s=100$. The curvature is computed through two methods: (left) using a minimum-norm interpolant and (right) employing Tikhonov regularization.}
    \label{fig:denoise}
\end{figure}

The results are shown in \autoref{fig:denoise}. In the first row, the effect of Tikhonov regularization is evident, with its maximum principal curvature \(\kappa_1\) being much smaller than that computed by the minimum-norm interpolation. Comparing \(\kappa_1\) values, Tikhonov regularization exhibits greater smoothing power, making it easier to distinguish the corners from the plates after one fairing step. However, this may not be an ideal property for fairing. After two fairing steps, it becomes apparent (though subtle) that the edges of the cube are sharper when using minimum-norm interpolation.
We reconstruct surfaces in the second and fourth column by applying the Poisson reconstruction method \cite{Poisson-Surfrecon} to the point cloud in the first and third column, respectively.

\section{Conclusion} \label{sec6}
This paper introduces a novel optimization framework for surface reconstruction that uses a minimum-norm interpolation approach and extends the RBF trial space by incorporating a KAN-inspired trial space. Our numerical experiments demonstrate that this KAN-inspired trial space significantly outperforms traditional RBF approaches in accurately reconstructing surfaces and estimating surface normals from unorganized point clouds.
By enriching the approximation space with lower-dimensional kernels, the proposed method provides improved flexibility and accuracy, offering a significant contribution towards advancing the field of surface reconstruction.

\reviewB{Despite these promising results, the proposed method has several limitations. Its performance depends on the quality and distribution of the input point cloud data, and the advantages of the KAN-inspired trial space may become less consistent for sparse or highly irregular sampling patterns. In addition, the method introduces extra tuning choices such as the selection of configurations which may affect robustness and computational efficiency across different datasets. These aspects indicate that although the method is effective, its practical performance may vary depending on the sampling condition and parameter selection.}

\reviewB{For future work, it would be valuable to develop adaptive strategies for selecting configurations and stencil structures based on local geometric or sampling indicators. A deeper theoretical and computational study of the relationship between the KAN-inspired trial space and RBF neural networks is also a promising direction. Such an investigation may not only strengthen the theoretical foundation of the proposed framework but also open new directions for advances in neural network design and more robust surface reconstruction methods.}

\section*{Acknowledgments}
This work was supported by the General Research Fund (GRF No. 12300922, 12301824) of Hong Kong Research Grant Council.

\section*{Declaration}
During the preparation of this work, the author(s) used GPT-5.x to improve language and readability. After using this tool/service, the authors reviewed and edited the content as needed and take full responsibility for the content of the publication.

\bibliographystyle{elsarticle-num}
\bibliography{Reference}

\begin{table}
    \centering
    \resizebox{\textwidth}{!}{\begin{tabular}{ |c|c||*{4}{wc{2.5cm}|}}
        \hline
        \multicolumn{6}{|c|}{Trial centers configuration within KANs first layer} \\
        \hline
        $N$ \& fill distance & $\tau$ & Config 1. $\Xi$ & Config 2. $\mathcal{R}\Xi$ & Config 3. $\mathcal{S}\Xi$ & Config 4. $\mathcal{R}\mathcal{S}\Xi$ \\
        \hline
        \multirow{2}{5em}{$N = 100  $}  & 2 & 1.32e$-$1 & 1.30e$-$1 & 2.92e$-$1 & 1.63e$-$1 \\
        \multirow{2}{5em}{$h = 0.229$}  & 3 & 1.07e$-$2 & 2.58e$-$2 & 2.26e$-$2 & 2.14e$-$2 \\
                                        & 4 & 1.31e$-$3 & 4.55e$-$3 & 2.00e$-$3 & 3.11e$-$3 \\
                                        & 5 & 2.43e$-$4 & 1.11e$-$3 & 3.50e$-$4 & 1.07e$-$3 \\
        \hline
        \multirow{2}{5em}{$N = 200  $}  & 2 & 7.83e$-$2 & 1.05e$-$1 & 3.29e$-$1 & 1.58e$-$1 \\
        \multirow{2}{5em}{$h = 0.175$}  & 3 & 5.06e$-$3 & 8.71e$-$3 & 6.91e$-$3 & 7.05e$-$3 \\
                                        & 4 & 7.19e$-$4 & 1.58e$-$3 & 3.33e$-$4 & 5.16e$-$4 \\
                                        & 5 & 1.28e$-$4 & 1.52e$-$4 & 1.03e$-$4 & 1.06e$-$4 \\
        \hline
        \multirow{2}{5em}{$N = 500  $}  & 2 & 5.58e$-$2 & 4.52e$-$2 & 1.99e$-$1 & 4.44e$-$2 \\
        \multirow{2}{5em}{$h = 0.102$}  & 3 & 1.91e$-$3 & 2.28e$-$3 & 1.68e$-$3 & 1.10e$-$3 \\
                                        & 4 & 1.60e$-$4 & 1.48e$-$4 & 7.71e$-$5 & 7.41e$-$5 \\
                                        & 5 & 1.44e$-$5 & 1.45e$-$5 & 1.35e$-$5 & 8.68e$-$6 \\
        \hline
        \multirow{2}{5em}{$N = 700  $}  & 2 & 4.88e$-$2 & 3.03e$-$2 & 7.64e$-$2 & 3.13e$-$2 \\
        \multirow{2}{5em}{$h = 0.0900$} & 3 & 1.68e$-$3 & 1.86e$-$3 & 8.37e$-$4 & 7.33e$-$4 \\
                                        & 4 & 5.74e$-$5 & 1.06e$-$4 & 2.48e$-$5 & 1.23e$-$5 \\
                                        & 5 & 4.37e$-$6 & 4.50e$-$6 & 4.29e$-$6 & 2.80e$-$6 \\
        \hline
        \multirow{2}{5em}{$N = 1000 $}  & 2 & 1.50e$-$1 & 1.55e$-$1 & 1.23e$-$1 & 1.85e$-$1 \\
        \multirow{2}{5em}{$h = 0.0743$} & 3 & 1.18e$-$2 & 1.22e$-$2 & \textbf{3.51e$-$3} & \textbf{1.02e$-$3} \\
                                        & 4 & 1.61e$-$4 & 2.16e$-$4 & 3.36e$-$5 & 4.50e$-$5 \\
                                        & 5 & 6.96e$-$6 & 9.59e$-$6 & 1.91e$-$5 & 5.23e$-$6 \\
        \hline
        \multirow{2}{5em}{$N = 2000 $}  & 2 & 2.54e$-$2 & 1.46e$-$2 & 1.41e$-$2 & \textbf{8.42e$-$4} \\
        \multirow{2}{5em}{$h = 0.0627$} & 3 & 3.68e$-$4 & 3.71e$-$4 & 1.07e$-$4 & 1.92e$-$5 \\
                                        & 4 & 5.37e$-$6 & 9.45e$-$6 & 2.98e$-$6 & 1.51e$-$6 \\
                                        & 5 & 6.73e$-$7 & 1.05e$-$6 & 1.17e$-$6 & 2.10e$-$7 \\
        \hline
        \multirow{2}{5em}{$N = 5000 $}  & 2 & 1.43e$-$2 & 9.51e$-$3 & 3.30e$-$3 & \textbf{3.15e$-$4} \\
        \multirow{2}{5em}{$h = 0.0401$} & 3 & 7.17e$-$5 & 8.19e$-$5 & 6.60e$-$6 & 1.40e$-$6 \\
                                        & 4 & 1.17e$-$6 & 2.39e$-$6 & 6.07e$-$7 & \textbf{1.58e$-$7} \\
                                        & 5 & 5.71e$-$7 & 1.88e$-$6 & 9.98e$-$7 & 1.30e$-$7 \\
        \hline
    \end{tabular}}
    \caption{Example 1: Normal estimation accuracy for various configurations on an ellipsoid. Configurations are tested across a range of stencil sizes ($40 \leq N_s \leq 80$) and kernel smoothness orders ($2 \leq \tau \leq 5$). Notable results where configurations significantly outperform the benchmark  Config 1  are highlighted in bold.}
    \label{table:KnR_config}
\end{table}

\begin{table}
    \centering
    \resizebox{\textwidth}{!}{\begin{tabular}{ |c|c||*{4}{wc{3cm}|}}
        \hline
        & & \multicolumn{2}{c|}{min-norm interpolant in HRBF $\mathcal{H}_{\Xi,\Phi_{\tau,3}}$} & \multicolumn{2}{c|}{min-norm interpolant in KRBF $\mathcal{U}_{\text{KAN}}$} \\
        \hline
        $N$ \& fill distance & $\tau$ & $\mathcal{N}$-norm \eqref{eq:HRBF_K_norm}   & $\|\boldsymbol\lambda\|_{\ell^2}$  & $\mathcal{N}$-norm \eqref{eq:KA_RBF_K_norm} & $\|\boldsymbol\lambda\|_{\ell^2}$\\
        \hline
        \multirow{2}{5em}{$N = 100  $}  & 2 & 5.46e$-$1 & 5.14e$-$1 & 1.63e$-$1 & 1.15e$-$2 \\
        \multirow{2}{5em}{$h = 0.229$}  & 3 & 1.73e$-$1 & 2.48e$-$2 & 2.14e$-$2 & 1.49e$-$3 \\
                                        & 4 & 2.11e$-$2 & 4.20e$-$3 & 3.11e$-$3 & 7.18e$-$4 \\
                                        & 5 & 5.37e$-$3 & 8.05e$-$4 & 1.07e$-$3 & 4.26e$-$4 \\
        \hline
        \multirow{2}{5em}{$N = 200  $}  & 2 & 2.98e$-$1 & 3.10e$-$1 & 1.58e$-$1 & 1.39e$-$2 \\
        \multirow{2}{5em}{$h = 0.175$}  & 3 & 8.20e$-$2 & 2.08e$-$2 & 7.05e$-$3 & 1.03e$-$3 \\
                                        & 4 & 3.83e$-$3 & 5.23e$-$3 & 5.16e$-$4 & 1.78e$-$4 \\
                                        & 5 & 2.70e$-$3 & 1.51e$-$3 & 1.06e$-$4 & 3.57e$-$5 \\
        \hline
        \multirow{2}{5em}{$N = 500  $}  & 2 & 3.75e$-$1 & 4.10e$-$1 & 4.44e$-$2 & 4.23e$-$3 \\
        \multirow{2}{5em}{$h = 0.102$}  & 3 & 5.22e$-$2 & 1.15e$-$2 & 1.10e$-$3 & 1.19e$-$4 \\
                                        & 4 & 1.73e$-$3 & 2.60e$-$3 & 7.41e$-$5 & 1.85e$-$5 \\
                                        & 5 & 5.31e$-$4 & 3.25e$-$4 & 8.68e$-$6 & 2.52e$-$6 \\
        \hline
        \multirow{2}{5em}{$N =  700  $} & 2 & 9.35e$-$1 & 8.55e$-$1 & 3.13e$-$2 & 2.73e$-$3 \\
        \multirow{2}{5em}{$h = 0.0900$} & 3 & 1.66e$-$1 & 1.82e$-$2 & 7.33e$-$4 & 6.14e$-$5 \\
                                        & 4 & 4.06e$-$3 & 2.04e$-$3 & 1.23e$-$5 & 6.52e$-$6 \\
                                        & 5 & 3.18e$-$4 & 1.51e$-$4 & 2.80e$-$6 & 1.72e$-$6 \\
        \hline
        \multirow{2}{5em}{$N = 1000  $} & 2 & 1.39e$+$0 & 1.35e$+$0 & 1.85e$-$1 & 1.15e$-$2 \\
        \multirow{2}{5em}{$h = 0.0743$} & 3 & 1.02e$+$0 & 8.15e$-$2 & 1.02e$-$3 & 1.58e$-$4 \\
                                        & 4 & 7.60e$-$2 & 1.14e$-$2 & 4.50e$-$5 & 1.23e$-$5 \\
                                        & 5 & 4.84e$-$3 & 1.05e$-$4 & 5.23e$-$6 & 3.63e$-$6 \\
        \hline
        \multirow{2}{5em}{$N = 2000  $} & 2 & 3.30e$-$1 & 2.58e$-$1 & 8.42e$-$4 & 2.98e$-$4 \\
        \multirow{2}{5em}{$h = 0.0627$} & 3 & 2.08e$-$2 & 6.86e$-$3 & 1.92e$-$5 & 9.27e$-$6 \\
                                        & 4 & 6.35e$-$4 & 7.36e$-$4 & 1.51e$-$6 & 9.37e$-$7 \\
                                        & 5 & 1.08e$-$4 & 3.75e$-$5 & 2.10e$-$7 & 1.19e$-$6 \\
        \hline
        \multirow{2}{5em}{$N = 5000  $} & 2 & 2.74e$-$1 & 1.88e$-$1 & 3.15e$-$4 & 1.20e$-$4 \\
        \multirow{2}{5em}{$h = 0.0401$} & 3 & 9.89e$-$3 & 6.87e$-$3 & 1.40e$-$6 & 1.09e$-$6 \\
                                        & 4 & 3.33e$-$4 & 2.56e$-$4 & 1.58e$-$7 & 8.00e$-$7 \\
                                        & 5 & 1.12e$-$5 & 3.71e$-$6 & 1.30e$-$7 & 7.24e$-$7 \\
        \hline
    \end{tabular}}
    \caption{Example 2: Normal estimation accuracy comparison between KAN-inspired and HRBF trial spaces with configuration similar to Example 1.}
    \label{table:opt_norm_choice}
\end{table}
\end{document}

%% file: Packages_and_Command.tex
\usepackage{amsmath,amssymb}
\usepackage{amsthm}
\usepackage{graphicx}
\graphicspath{figure/}
\usepackage{subfig}
\usepackage[percent]{overpic}
\usepackage{pict2e} 
\usepackage{hyperref}
\usepackage[svgnames]{xcolor}
\usepackage{parskip}
\usepackage{tikz}
\usepackage{bm}
\usepackage{multirow}
\usepackage{array}
\usepackage{anysize}
\usepackage{soul}

\newtheorem{definition}{Definition}[section]
\newtheorem{remark}{Remark}[section]

\newcommand{\reviewA}[1]{{#1}}
\newcommand{\reviewB}[1]{{#1}}
\newcommand{\new}[1]{{#1}}

\NewDocumentCommand{\csq}{om}{%
  \IfNoValueTF{#1}{%
    \textcolor{#2}{\rule{1.8ex}{1ex}}%
  }{%
    \textcolor[#1]{#2}{\rule{1.8ex}{1ex}}%
  }%
}

\def\M{{\mathcal{M}}}
\def\R{{\mathbb{R}}}
\def\x{{\mathbf{x}}}

\def\p{{\mathbf{p}}}
\def\n{{\hat{\mathbf{n}}}}

%% file: sub_opt_KANs_1.tex
\tikzset{every picture/.style={line width=0.75pt}} 

\begin{tikzpicture}[x=0.75pt,y=0.75pt,yscale=-1,xscale=1]

\draw [color={rgb, 255:red, 0; green, 0; blue, 0 }  ,draw opacity=1 ]   (229.33,869.34) -- (169.33,840.01) -- (229.33,869.34) -- (202.67,840.67) -- (229.33,869.34) -- (272.67,841.34) ;
\draw [shift={(229.33,869.34)}, rotate = 206.05] [color={rgb, 255:red, 0; green, 0; blue, 0 }  ,draw opacity=1 ][fill={rgb, 255:red, 0; green, 0; blue, 0 }  ,fill opacity=1 ][line width=0.75]      (0, 0) circle [x radius= 3.35, y radius= 3.35]   ;
\draw [color={rgb, 255:red, 0; green, 0; blue, 0 }  ,draw opacity=1 ]   (200.67,780.01) -- (169,818.69) -- (200.67,780.01) -- (311.33,819.34) -- (200.67,780.01) -- (451.33,818.67) ;
\draw [shift={(200.67,780.01)}, rotate = 129.31] [color={rgb, 255:red, 0; green, 0; blue, 0 }  ,draw opacity=1 ][fill={rgb, 255:red, 0; green, 0; blue, 0 }  ,fill opacity=1 ][line width=0.75]      (0, 0) circle [x radius= 3.35, y radius= 3.35]   ;
\draw [color={rgb, 255:red, 0; green, 0; blue, 0 }  ,draw opacity=1 ]   (260.67,780.01) -- (202.67,819.34) -- (260.67,780.01) -- (344,819.34) -- (260.67,780.01) -- (485.33,819.34) ;
\draw [shift={(260.67,780.01)}, rotate = 145.86] [color={rgb, 255:red, 0; green, 0; blue, 0 }  ,draw opacity=1 ][fill={rgb, 255:red, 0; green, 0; blue, 0 }  ,fill opacity=1 ][line width=0.75]      (0, 0) circle [x radius= 3.35, y radius= 3.35]   ;
\draw [color={rgb, 255:red, 0; green, 0; blue, 0 }  ,draw opacity=1 ]   (503.33,782.01) -- (272,819.34) -- (503.33,782.01) -- (412.67,818.67) -- (503.33,782.01) -- (554.67,818.67) ;
\draw [shift={(503.33,782.01)}, rotate = 0] [color={rgb, 255:red, 0; green, 0; blue, 0 }  ,draw opacity=1 ][fill={rgb, 255:red, 0; green, 0; blue, 0 }  ,fill opacity=1 ][line width=0.75]      (0, 0) circle [x radius= 3.35, y radius= 3.35]   ;
\draw [color={rgb, 255:red, 0; green, 0; blue, 0 }  ,draw opacity=1 ]   (363.33,680.67) -- (201.33,729.34) -- (363.33,680.67) -- (264,728.01) -- (363.33,680.67) -- (504.67,728.67) ;
\draw [shift={(363.33,680.67)}, rotate = 163.28] [color={rgb, 255:red, 0; green, 0; blue, 0 }  ,draw opacity=1 ][fill={rgb, 255:red, 0; green, 0; blue, 0 }  ,fill opacity=1 ][line width=0.75]      (0, 0) circle [x radius= 3.35, y radius= 3.35]   ;
\draw    (200.67,750.01) -- (200.67,780.01) ;
\draw    (260.67,748.67) -- (260.67,780.01) ;
\draw    (503.33,750.01) -- (503.33,782.01) ;
\draw  [dash pattern={on 0.84pt off 2.51pt}] (152.71,829.67) .. controls (152.71,821.04) and (159.7,814.05) .. (168.33,814.05) .. controls (176.96,814.05) and (183.96,821.04) .. (183.96,829.67) .. controls (183.96,838.3) and (176.96,845.3) .. (168.33,845.3) .. controls (159.7,845.3) and (152.71,838.3) .. (152.71,829.67) -- cycle ;
\draw    (191,737.98) .. controls (206.65,708.16) and (196.03,770.1) .. (210,740.27) ;
\draw   (158,819.34) -- (178.67,819.34) -- (178.67,840.01) -- (158,840.01) -- cycle ;
\draw    (158,819.34) .. controls (162.67,844.67) and (174.67,844.01) .. (178.67,819.34) ;
\draw   (192,819.34) -- (212.67,819.34) -- (212.67,840.01) -- (192,840.01) -- cycle ;
\draw    (192,819.34) .. controls (196.67,844.67) and (208.67,844.01) .. (212.67,819.34) ;
\draw  [dash pattern={on 0.84pt off 2.51pt}]  (222,829.67) -- (249.33,830.01) ;
\draw   (260.67,819.34) -- (281.33,819.34) -- (281.33,840.01) -- (260.67,840.01) -- cycle ;
\draw    (260.67,819.34) .. controls (265.33,844.67) and (277.33,844.01) .. (281.33,819.34) ;
\draw [color={rgb, 255:red, 0; green, 0; blue, 0 }  ,draw opacity=1 ]   (371.33,869.34) -- (311.33,840.01) -- (371.33,869.34) -- (344.67,840.67) -- (371.33,869.34) -- (414.67,841.34) ;
\draw [shift={(371.33,869.34)}, rotate = 206.05] [color={rgb, 255:red, 0; green, 0; blue, 0 }  ,draw opacity=1 ][fill={rgb, 255:red, 0; green, 0; blue, 0 }  ,fill opacity=1 ][line width=0.75]      (0, 0) circle [x radius= 3.35, y radius= 3.35]   ;
\draw   (300,819.34) -- (320.67,819.34) -- (320.67,840.01) -- (300,840.01) -- cycle ;
\draw    (300,819.34) .. controls (304.67,844.67) and (316.67,844.01) .. (320.67,819.34) ;
\draw   (334,819.34) -- (354.67,819.34) -- (354.67,840.01) -- (334,840.01) -- cycle ;
\draw    (334,819.34) .. controls (338.67,844.67) and (350.67,844.01) .. (354.67,819.34) ;
\draw  [dash pattern={on 0.84pt off 2.51pt}]  (364,829.67) -- (391.33,830.01) ;
\draw   (402.67,819.34) -- (423.33,819.34) -- (423.33,840.01) -- (402.67,840.01) -- cycle ;
\draw    (402.67,819.34) .. controls (407.33,844.67) and (419.33,844.01) .. (423.33,819.34) ;
\draw [color={rgb, 255:red, 0; green, 0; blue, 0 }  ,draw opacity=1 ]   (512.67,869.34) -- (452.67,840.01) -- (512.67,869.34) -- (486,840.67) -- (512.67,869.34) -- (556,841.34) ;
\draw [shift={(512.67,869.34)}, rotate = 206.05] [color={rgb, 255:red, 0; green, 0; blue, 0 }  ,draw opacity=1 ][fill={rgb, 255:red, 0; green, 0; blue, 0 }  ,fill opacity=1 ][line width=0.75]      (0, 0) circle [x radius= 3.35, y radius= 3.35]   ;
\draw   (441.33,819.34) -- (462,819.34) -- (462,840.01) -- (441.33,840.01) -- cycle ;
\draw    (441.33,819.34) .. controls (446,844.67) and (458,844.01) .. (462,819.34) ;
\draw   (475.33,819.34) -- (496,819.34) -- (496,840.01) -- (475.33,840.01) -- cycle ;
\draw    (475.33,819.34) .. controls (480,844.67) and (492,844.01) .. (496,819.34) ;
\draw  [dash pattern={on 0.84pt off 2.51pt}]  (505.33,829.67) -- (532.67,830.01) ;
\draw   (544,819.34) -- (564.67,819.34) -- (564.67,840.01) -- (544,840.01) -- cycle ;
\draw    (544,819.34) .. controls (548.67,844.67) and (560.67,844.01) .. (564.67,819.34) ;
\draw   (190,728.67) -- (210.67,728.67) -- (210.67,749.34) -- (190,749.34) -- cycle ;
\draw   (250.67,728.67) -- (271.33,728.67) -- (271.33,749.34) -- (250.67,749.34) -- cycle ;
\draw   (493.33,729.34) -- (514,729.34) -- (514,750.01) -- (493.33,750.01) -- cycle ;
\draw  [dash pattern={on 0.84pt off 2.51pt}]  (282.33,740.01) -- (482,738.67) ;
\draw    (251,738.65) .. controls (266.65,708.83) and (256.03,770.76) .. (270,740.94) ;
\draw    (494.33,738.65) .. controls (509.98,708.83) and (499.36,770.76) .. (513.33,740.94) ;
\draw    (152.71,829.67) -- (113.02,804.41) ;
\draw [shift={(111.33,803.34)}, rotate = 32.47] [color={rgb, 255:red, 0; green, 0; blue, 0 }  ][line width=0.75]    (10.93,-3.29) .. controls (6.95,-1.4) and (3.31,-0.3) .. (0,0) .. controls (3.31,0.3) and (6.95,1.4) .. (10.93,3.29)   ;
\draw  [dash pattern={on 0.84pt off 2.51pt}] (184.71,739.01) .. controls (184.71,730.38) and (191.7,723.38) .. (200.33,723.38) .. controls (208.96,723.38) and (215.96,730.38) .. (215.96,739.01) .. controls (215.96,747.63) and (208.96,754.63) .. (200.33,754.63) .. controls (191.7,754.63) and (184.71,747.63) .. (184.71,739.01) -- cycle ;
\draw    (184.71,739.01) -- (169.4,723.43) ;
\draw [shift={(168,722.01)}, rotate = 45.5] [color={rgb, 255:red, 0; green, 0; blue, 0 }  ][line width=0.75]    (10.93,-3.29) .. controls (6.95,-1.4) and (3.31,-0.3) .. (0,0) .. controls (3.31,0.3) and (6.95,1.4) .. (10.93,3.29)   ;

\draw (229.33,869.34) node [anchor=north west][inner sep=0.75pt]   [align=left] {$\displaystyle \xi _{0}{}_{,}{}_{1}=x$};
\draw (115.37,695.17) node [anchor=south] [inner sep=0.75pt]   [align=left] {};
\draw (333.33,680.67) node [anchor=south west] [inner sep=0.75pt]   [align=left] {\begin{minipage}[lt]{54.09pt}\setlength\topsep{0pt}
\begin{center}
$\displaystyle  F \in \mathcal{K}_{\Xi }{}_{,}{}_{\Phi }{}_{_{\tau }}$
\end{center}

\end{minipage}};
\draw (111.33,803.34) node [anchor=south] [inner sep=0.75pt]   [align=left] 
{$\displaystyle \varphi_1 ( \xi ) \ =\ a^2(\xi-b)^2$};
\draw (371.33,869.34) node [anchor=north west][inner sep=0.75pt]   [align=left] {$\displaystyle \xi _{0}{}_{,}{}_{2}=y$};
\draw (512.67,869.34) node [anchor=north west][inner sep=0.75pt]   [align=left] {$\displaystyle \xi _{0}{}_{,}{}_{3}=z$};
\draw (200.67,780.01) node [anchor=south west] [inner sep=0.75pt]   [align=left] {$\displaystyle \xi _{1}{}_{,}{}_{1}$};
\draw (260.67,780.01) node [anchor=south west] [inner sep=0.75pt]   [align=left] {$\displaystyle \xi _{1}{}_{,}{}_{2}$};
\draw (503.33,782.01) node [anchor=south west] [inner sep=0.75pt]   [align=left] {$\displaystyle \xi _{1}{}_{,}{}_{N}$};
\draw (168,722.01) node [anchor=south east] [inner sep=0.75pt]   [align=left] 
{$\displaystyle \varphi_2 (\xi) \ = c\,\Phi_{\tau+1,3} (\sqrt{\xi}) \ $};

\end{tikzpicture} 